\documentclass[11pt]{elsarticle}

\usepackage[]{hyperref}

\journal{Journal of Computational Physics}
\usepackage[margin=2.5cm]{geometry}
\biboptions{sort&compress}

\newcommand{\figref}[1]{\hyperref[#1]{Fig.~\ref*{#1}}}
\newcommand{\secref}[1]{\hyperref[#1]{Section~\ref*{#1}}}
\newcommand{\tabref}[1]{\hyperref[#1]{Table~\ref*{#1}}}
\newcommand{\eqnref}[1]{\hyperref[#1]{Eq.~(\ref*{#1})}}

\makeatletter
\providecommand{\doi}[1]{%
	\begingroup
	\let\bibinfo\@secondoftwo
	\urlstyle{rm}%
	\href{http://dx.doi.org/#1}{%
		doi:\discretionary{}{}{}%
		\nolinkurl{#1}%
	}%
	\endgroup
}
\makeatother

\usepackage{graphicx}
\usepackage{dcolumn}
\usepackage{bm}
\usepackage{hyperref}
\usepackage{epstopdf}
\usepackage{amsmath,esint}
\usepackage{amsmath}
\usepackage{amsfonts}
\usepackage{amssymb}
\usepackage{graphicx}
\usepackage{amsmath}
\usepackage{amsthm}
\usepackage{eucal}
\usepackage{amssymb}
\usepackage{mathrsfs}
\usepackage{float}
\usepackage{hyperref}
\usepackage{lineno}
\usepackage{multirow}
\usepackage{booktabs}
\usepackage{graphicx}
\usepackage{subcaption}
\usepackage{dashrule}
\usepackage{adjustbox}
\usepackage[pdftex,dvipsnames]{xcolor}  
\usepackage{graphicx}
\usepackage{dcolumn}
\usepackage{bm}
\usepackage{natbib}
\usepackage{amssymb}
\usepackage{amsthm}
\usepackage{mathtools}
\usepackage{empheq}
\usepackage[most]{tcolorbox}
\usepackage{mathtools}
\usepackage{cleveref}
\usepackage{soul}
\usepackage{array}

\usepackage{tikz}
\usetikzlibrary{shapes,arrows,chains}
\usetikzlibrary{calc}
\usetikzlibrary{positioning}
\usepackage{pgfplots}
\usepackage{pgfplotstable}
\usepackage{grffile}
\usetikzlibrary{arrows,shapes,plotmarks}
\usepgfplotslibrary{groupplots}
\usetikzlibrary{matrix}
\usepackage{tikz,pgfplots}
\usetikzlibrary{snakes,arrows,shapes,trees}

\usepackage{tgpagella}
\usepackage{newpxmath}

\renewcommand{\vec}[1]{\mathbf{#1}}

\usepackage{xspace}
\let\storeBeta=\beta
\renewcommand\beta{\relax\ifmmode{\storeBeta}\else{$\storeBeta$}\fi\xspace}

\usepackage{lscape}

\let\storeAlpha=\alpha
\renewcommand\alpha{\relax\ifmmode{\storeAlpha}\else{$\storeAlpha$}\fi\xspace}

\newcommand{\pd}[2]{\frac{\partial #1}{\partial #2}} 
\newcommand{\pdd}[2]{\frac{\partial^2 #1}{\partial #2^2}} 
\let\baraccent=\= 
\renewcommand{\=}[1]{\stackrel{#1}{=}} 

\newcount\colveccount
\newcommand*\colvec[1]{
	\global\colveccount#1
	\begin{pmatrix}
		\colvecnext
	}
	\def\colvecnext#1{
		#1
		\global\advance\colveccount-1
		\ifnum\colveccount>0
		\\
		\expandafter\colvecnext
		\else
	\end{pmatrix}
	\fi
}

\def\tvi{\widetilde{v}_i}
\def\tvj{\widetilde{v}_j}

\hypersetup{
	colorlinks,
	linkcolor={blue!50!black},
	citecolor={blue!50!black},
	urlcolor={blue!80!black},
	anchorcolor = {blue!80!black},
	filecolor = {blue!80!black},
	menucolor = {blue!80!black},
	runcolor = {blue!80!black}
}

\usepackage{tcolorbox}
\tcbuselibrary{theorems}
\usepackage{varwidth}
\tcbuselibrary{skins}
\tcbuselibrary{breakable}
\tcbsetforeverylayer{shield externalize}
\newtcbtheorem[crefname={theorem}{theorems}]{theorem}{Theorem}%
{	enhanced,
	frame empty,
	interior empty,
	colframe=black,
	coltitle=black,
	fonttitle=\bfseries,
	colbacktitle=gray!10!white,
	borderline={0.5mm}{0mm}{gray!15!white},
	borderline={0.5mm}{0mm}{black!50!white,dashed},
	attach boxed title to top center={yshift=-2mm},
	boxed title style={boxrule=0.4pt},
	attach boxed title to top left={%
		yshift*=-\tcboxedtitleheight/2, 
		xshift=5mm},
	varwidth boxed title}{thrm}

\newtcbtheorem[crefname={lemma}{lemmas}]{lemma}{Lemma}%
{	enhanced,
	frame empty,
	interior empty,
	colframe=black,
	coltitle=black,
	fonttitle=\bfseries,
	colbacktitle=gray!10!white,
	borderline={0.5mm}{0mm}{gray!15!white},
	borderline={0.5mm}{0mm}{black!50!white,dashed},
	attach boxed title to top center={yshift=-2mm},
	boxed title style={boxrule=0.4pt},
	attach boxed title to top left={%
		yshift*=-\tcboxedtitleheight/2, 
		xshift=5mm},
	varwidth boxed title}{lem}

\newtcbtheorem[crefname={proposition}{propositions}]{proposition}{Proposition}%
{	enhanced,
	breakable,
	frame empty,
	interior empty,
	colframe=black,
	coltitle=black,
	fonttitle=\bfseries,
	colbacktitle=gray!10!white,
	borderline={0.5mm}{0mm}{gray!15!white},
	borderline={0.5mm}{0mm}{black!50!white,dashed},
	attach boxed title to top center={yshift=-2mm},
	boxed title style={boxrule=0.4pt},
	attach boxed title to top left={%
		yshift*=-\tcboxedtitleheight/2, 
		xshift=5mm},
	varwidth boxed title}{prop}

\newtcbtheorem[crefname={definition}{definitions}]{definition}{Definition}%
{	enhanced,
	frame empty,
	interior empty,
	colframe=black,
	coltitle=black,
	fonttitle=\bfseries,
	colbacktitle=gray!10!white,
	borderline={0.5mm}{0mm}{gray!15!white},
	borderline={0.5mm}{0mm}{black!50!white,dashed},
	attach boxed title to top center={yshift=-2mm},
	boxed title style={boxrule=0.4pt},
	attach boxed title to top left={%
		yshift*=-\tcboxedtitleheight/2, 
		xshift=5mm},
	varwidth boxed title}{defn}

\newtcbtheorem[crefname={claim}{claims}]{claim}{Claim}%
{	enhanced,
	frame empty,
	interior empty,
	colframe=black,
	coltitle=black,
	fonttitle=\bfseries,
	colbacktitle=gray!10!white,
	borderline={0.5mm}{0mm}{gray!15!white},
	borderline={0.5mm}{0mm}{black!50!white,dashed},
	attach boxed title to top center={yshift=-2mm},
	boxed title style={boxrule=0.4pt},
	attach boxed title to top left={%
		yshift*=-\tcboxedtitleheight/2, 
		xshift=5mm},
	varwidth boxed title}{claim}

\newtcbtheorem[crefname={corollary}{corollarys}]{corollary}{Corollary}%
{	enhanced,
	frame empty,
	interior empty,
	colframe=black,
	coltitle=black,
	fonttitle=\bfseries,
	colbacktitle=gray!10!white,
	borderline={0.5mm}{0mm}{gray!15!white},
	borderline={0.5mm}{0mm}{black!50!white,dashed},
	attach boxed title to top center={yshift=-2mm},
	boxed title style={boxrule=0.4pt},
	attach boxed title to top left={%
		yshift*=-\tcboxedtitleheight/2, 
		xshift=5mm},
	varwidth boxed title}{cor}

\tcolorboxenvironment{proof}{
	blanker,
	breakable,
	left=5mm,
	borderline west={0.5mm}{0pt}{black!50!white}
}

\newtheorem{remark}{Remark}

\usepackage{pgfplots}
\usepgfplotslibrary{colorbrewer}
\pgfplotsset{
	table/search path={Figures},
}
\usetikzlibrary{arrows,shapes}
\pgfplotsset{compat=1.3}
\pgfplotsset{
	discard if/.style 2 args={
		x filter/.code={
			\edef\tempa{\thisrow{#1}}
			\edef\tempb{#2}
			\ifx\tempa\tempb
			
			\fi
		}
	},
	discard if not/.style 2 args={
		x filter/.code={
			\edef\tempa{\thisrow{#1}}
			\edef\tempb{#2}
			\ifx\tempa\tempb
			\else
			
			\fi
		}
	}
}
\usetikzlibrary{plotmarks}
\usetikzlibrary{spy}
\usetikzlibrary{arrows,shapes,plotmarks}
\usepgfplotslibrary{groupplots}

\usetikzlibrary{calc}

\usetikzlibrary{matrix}

\usepackage{xargs}                      
\usepackage{siunitx}

\usepackage[colorinlistoftodos,prependcaption,textsize=normalsize]{todonotes}

\newcommand{\petsc}{\textsc{Petsc}}

\newcommand{\dkt}{\textsc{Dendro-KT}}
\newcommand{\In}{\textsc{InActive}}
\newcommand{\Out}{\textsc{Active}}
\newcommand{\Intercepted}{\textsc{Intercepted}}

\usepackage{listings}

\usepackage{xparse}

\NewDocumentCommand{\codeword}{v}{%
	\texttt{\textcolor{blue}{#1}}%
}

\definecolor{codegreen}{rgb}{0,0.6,0}
\definecolor{codegray}{rgb}{0.5,0.5,0.5}
\definecolor{codepurple}{rgb}{0.58,0,0.82}
\definecolor{backcolour}{rgb}{0.95,0.95,0.92}

\lstdefinestyle{mystyle}{
	backgroundcolor=\color{backcolour},   
	commentstyle=\color{codegreen},
	keywordstyle=\color{magenta},
	stringstyle=\color{codepurple},
	basicstyle=\ttfamily\footnotesize,
	breakatwhitespace=false,         
	breaklines=true,                 
	captionpos=b,                    
	keepspaces=true,                  
	showspaces=false,                
	showstringspaces=false,
	showtabs=false,                  
	tabsize=2
}

\usetikzlibrary{calc,matrix,decorations.markings,decorations.pathreplacing,shapes.geometric,backgrounds}

\definecolor{colone}{RGB}{209,220,204}
\definecolor{coltwo}{RGB}{207,215,210}
\definecolor{colthree}{RGB}{207,233,232}
\definecolor{colfour}{RGB}{248,243,214}
\definecolor{colfive}{RGB}{245,238,197}
\definecolor{colsix}{RGB}{243,235,179}
\definecolor{colseven}{RGB}{241,231,163}

\pgfplotsset{compat=newest,scaled y ticks=true} 
\usepgfplotslibrary{units} 
\usepgfplotslibrary{colorbrewer}

\definecolor{sq_b1}{RGB}{37,52,148}
\definecolor{sq_b2}{RGB}{44,127,184}
\definecolor{sq_b3}{RGB}{65,182,196}
\definecolor{sq_b4}{RGB}{127,205,187}
\definecolor{sq_b5}{RGB}{199,233,180}
\definecolor{sq_b6}{RGB}{255,255,204}

\definecolor{sq_r1}{RGB}{189,0,38}
\definecolor{sq_r2}{RGB}{240,59,32}
\definecolor{sq_r3}{RGB}{253,141,60}
\definecolor{sq_r4}{RGB}{254,178,76}
\definecolor{sq_r5}{RGB}{254,217,118}
\definecolor{sq_r6}{RGB}{255,255,178}

\definecolor{sq_g1}{RGB}{0,104,55}
\definecolor{sq_g2}{RGB}{49,163,84}
\definecolor{sq_g3}{RGB}{120,198,121}
\definecolor{sq_g4}{RGB}{173,221,142}
\definecolor{sq_g5}{RGB}{217,240,163}
\definecolor{sq_g6}{RGB}{255,255,204}

\definecolor{div_c1}{RGB}{230,171,2}
\definecolor{div_c2}{RGB}{102,166,30}
\definecolor{div_c3}{RGB}{231,41,138}
\definecolor{div_c4}{RGB}{117,112,179}
\definecolor{div_c5}{RGB}{217,95,2}
\definecolor{div_c6}{RGB}{27,158,119}
\definecolor{div_c7}{RGB}{215,48,39}

\definecolor{div_d1}{RGB}{215,25,28}
\definecolor{div_d2}{RGB}{253,174,97}
\definecolor{div_d3}{RGB}{255,255,191}
\definecolor{div_d4}{RGB}{171,217,233}
\definecolor{div_d5}{RGB}{44,123,182}

\allowdisplaybreaks
\usetikzlibrary{external}
\begin{document}
	
	\begin{frontmatter}
		
		\title{Direct numerical simulation of electrokinetic transport phenomena in fluids: variational multi-scale stabilization and octree-based mesh refinement\footnote{Here, by Direct numerical simulation, we mean that we resolve the physics of interest at appropriate time and length scales}}
		
		
		\author[isuMechEAddress,isuChemAddress]{Sungu Kim\corref{contrib}}
		%
		\author[isuMechEAddress]{Kumar Saurabh\corref{contrib}}
		%
		\author[stanfordaddress]{Makrand A. Khanwale\corref{contrib}}
		%
		\author[stanfordaddress]{Ali Mani}
		\author[isuChemAddress]{Robbyn K. Anand}
		%
		\author[isuMechEAddress]{Baskar~Ganapathysubramanian\corref{correspondingAuthor}}
		\ead{baskarg@iastate.edu}
		\cortext[contrib]{These three authors contributed equally}
		\cortext[correspondingAuthor]{Corresponding author}
		\address[isuMechEAddress]{Department of Mechanical Engineering, Iowa State University, Iowa, USA 50011}
		\address[isuChemAddress]{Department of Chemistry, Iowa State University, Iowa, USA 50011}
		\address[stanfordaddress]{Center for Turbulence Research, Stanford University, Stanford,  California, USA  94305 }

		\begin{abstract}
			Computational modeling of charged species transport has enabled the analysis, design, and optimization of a diverse array of electrochemical and electrokinetic devices. These systems are represented by the Poisson-Nernst-Planck (PNP) equations coupled with the Navier-Stokes (NS) equation. Direct numerical simulation (DNS) to accurately capture the spatio-temporal variation of ion concentration and current flux remains challenging due to the (a) small critical dimension of the diffuse charge layer (DCL)
			, (b) stiff coupling due to fast charge relaxation times, large advective effects, and steep gradients close to boundaries, and (c) complex geometries exhibited by electrochemical devices. 
			
			In the current study, we address these challenges by presenting a direct numerical simulation framework that incorporates (a) a variational multiscale (VMS) treatment, (b) a block-iterative strategy in conjunction with semi-implicit (for NS) and implicit (for PNP) time integrators, and (c) octree based adaptive mesh refinement. The VMS formulation provides numerical stabilization critical for capturing the electro-convective flows often observed in engineered devices. The block-iterative strategy decouples the difficulty of non-linear coupling between the NS and PNP equations and allows the use of tailored numerical schemes separately for NS and PNP equations. The carefully designed second-order, hybrid implicit methods circumvent the harsh timestep requirements of explicit time steppers, thus enabling simulations over longer time horizons. Finally, the octree-based meshing allows efficient and targeted spatial resolution of the DCL. These features are incorporated into a massively parallel computational framework, enabling the simulation of realistic engineering electrochemical devices. The numerical framework is illustrated using several challenging canonical examples.  
		\end{abstract}
		
		
	\end{frontmatter}
	
	
	\section{Introduction}
	\label{sec:introduction}
	
	A diverse set of energy conversion, energy storage, manufacturing, and healthcare processes involve the transport and interaction of charged species in aqueous environments. These interactions span multiple time and length scales, often involving complex geometries ~\citep{rubinstein1997electric,fleury1993coupling,fleury1994mechanism,huth1995role,kwak2013shear,wessling2014nanometer}. Examples of such processes include electrolysis for hydrogen production, electrodialysis for water desalination, gas diffusion electrodes for chemical conversion, carbon capture, as well as several other processes that have a direct impact on sustainability and climate resilience~\citep{chi2018water, kim2010direct, kwak2016enhanced, nikonenko2014desalination, lees2022gas, weng2018modeling}.
	
	Quantitative analyses of transport processes in devices are performed using continuum mechanics-derived partial differential equations (PDEs), specifically the Navier-Stokes (NS) governing the flow and the Poisson-Nernst-Plank (PNP) equations governing chemical species transport and electrostatic interactions. Detailed, high-fidelity simulations- i.e., direct numerical simulations- are critical for comprehensive understanding, design, and quantitative prediction of system-level responses of electrokinetic systems. The coupling between the NS and PNP equation is responsible for producing a variety of very interesting phenomena~\citep{chi2018water,kim2010direct,kwak2016enhanced,nikonenko2014desalination,lees2022gas,weng2018modeling,kim2020concentration,anand2022concentration,devarakonda2022designing}. For instance, even in systems with small hydrodynamic velocities, i.e., low Reynolds (Re) numbers, a strong applied electric field triggers strongly non-linear and multi-scale electrokinetic instabilities~\citep{druzgalski2013direct}. 
	However, the complexities arising from the tight coupling of multiscale and multiphysics phenomena in these processes make simulating these systems very challenging. A comprehensive simulation framework can significantly accelerate the exploration and understanding of complex electrochemical processes and enable the design and optimization of societally relevant applications like more efficient water electrolysis and cheaper and lightweight dialysis machines.

	One can broadly identify three challenges to performing direct numerical simulations of NS-PNP equations in complex geometries: First, electrokinetic systems involve coupled mass transport, fluid flow, and electrodynamics. Additionally, the coupled non-linear NS-PNP system involves a variety of stiffness.  Therefore, such coupling and stiffnesses in the NS-PNP system call for numerical methods that can robustly and efficiently the multiscale physics. Second, such systems often exhibit thin regions with a high electric field and high gradients in ion concentration in the diffused charge layer. Accurately resolving these near-boundary gradients is essential, as they determine not only the bulk behavior but also impact the boundary fluxes, which are often the key quantities of interest from an experimental point of view ~\citep{kim2017ion,dydek2011overlimiting,han2014over,nam2015experimental}. The challenge of accurately resolving these regions is further exacerbated by the complex geometries exhibited by various electrochemical devices (packed beds, dendritic surfaces). This calls for approaches that are endowed with the ability to spatially resolve complex geometries in a scalable fashion. Third, simulating practical systems of engineering interest requires performing simulations across a long time horizon. This calls for time discretizations that are second-order (or higher) and have good stabilizing properties, thus offering the ability to successfully capture long time horizon simulations affordably~\citep{mani2020electroconvection}. 
	
	Resolving these challenges requires careful design of spatial stabilization and temporal discretization as well as appropriate strategies of (de)coupled solution of the equations and adaptive mesh refinement. To resolve those challenges, \citet{druzgalski2013direct} performed pioneering direct numerical simulations of electroconvective instabilities. They used a non-uniform mesh in the wall-normal direction to resolve the DCLs for small Debye lengths and a semi-implicit time scheme.  The framework was then pushed for more complex surfaces with patterned electrochemical properties in~\citep{Davidson2014, Davidson2016}. One of the main challenges in simulating these systems is that if one resolves the thin DCLs with a fine mesh resolution, one often requires very small timesteps for stability reasons. The problem of resolving thin DCLs is exacerbated when the problem length scale is larger, i.e., small non-dimensional Debye lengths, requiring finer meshes and smaller timesteps. \citet{druzgalski2013direct}~resolved the challenge of small timesteps using implicit time-stepping. 
	Specifically, the framework in \citep{druzgalski2013direct} utilized Fourier methods to efficiently handle linear equations arising from stiffnesses in NS and PNP independently. However, such strategies limit simulations to uniform meshes and periodic boundary conditions in the transverse directions, preventing their extendability to complex geometries. Additionally, structured meshes enforced fine resolution in the wall-tangential direction, even in regions away from the walls, where no physical small-scale structure exists. 
	
	There have been computational studies using commercial software to analyze the non-linear behavior of electrokinetics~\citep{jia2014multiphysics,gong2018direct,seo2022non}.
	It is important to note that most commercial softwares are unsuitable for such electrochemical systems, either because the numerical methods are not designed for these phenomena or are not scalable enough to accurately resolve the spatial and temporal scales required. \citet{karatay2015simulation} showed that with a carefully designed numerical method, their code was 20$\times$ faster than the commercial alternative (COMSOL in this case) on a single processor. This advantage amplified to about 160$\times$ speed up when both the codes were deployed in parallel but on a single node. Therefore, careful design of numerical methods and efficient parallel implementation are essential for affordable high-fidelity simulations of electrochemical systems. While there has been some progress in developing flexible frameworks~\citep{roberts2014framework,roberts2016insights}, numerical challenges involved in solving the PNP equations still remain. 
	
	In this work, we seek to resolve all three challenges detailed above for the direct numerical simulation of electrokinetic phenomena. We build upon our previous work on a finite elements based framework for solving the NS and PNP equations~\citep{kim2022computational}.  In contrast to~\citep{kim2020concentration}, we deploy a novel framework that utilizes a Variational multiscale method-based finite element approach and fast octree-based meshing.  Finite element approaches have been successfully applied for high-resolution simulations of NS~(see \citep{Volker2016} for an overview) and PNP~\citep{lim2007transient, lu2010poisson, sun2016error, gao2017linearized, he2018mixed, liu2021efficient,liu2022positivity,dehghan2022optimal,chen2022electrokinetic} equations. They are well suited for the NS-PNP system of equations due to (a) the use of a variational formulation that allows the weakening of derivatives, (b) the natural incorporation of heterogeneous and mixed boundary conditions often seen in these systems, and (c) ability to construct rigorous \textit{a posteriori} error estimates for mesh adaptivity that enable substantial computational gains.  
	
	The main contributions of this paper are as follows: 
	\begin{enumerate}
		\item We decouple the solution procedure of NS and PNP equations in a block-solve strategy. The PNP equations are solved using a fully implicit Newton method-based non-linear solver, ensuring the choice of large time steps even for cases with fast charge relaxation times and very fine mesh resolution used for thin DCLs. On the other hand, the Navier-Stokes equations are discretized using a semi-implicit linear discretization, which also ensures stability for larger time steps. In conjunction with this strategy, we utilize a second-order time discretization for both NS and PNP to ensure accuracy across long time horizons. Along with carefully chosen preconditioners, this block-iterative strategy allows robust and efficient DNS computations.
		
		\item We use a variational multi-scale (VMS) based stabilization for the pressure that allows for a pressure-coupled solver in conjunction with continuous Galerkin finite elements. Additionally, the VMS method models fine-scale velocities and ion concentrations not resolved by the grid size. This allows us to use adaptive mesh discretization without compromising the accuracy of the field variables.
		
		\item We use a massively parallel octree-based meshing framework to resolve the thin DCLs for small Debye lengths and mesh complex geometries in 3D. We use \dkt{}~\citep{MasadoDKT}, a highly scalable parallel octree library, to generate, manage and traverse adaptive octree-based meshes in distributed memory. Octrees are widely used in computational sciences~\cite{SundarSampathBiros08,BursteddeWilcoxGhattas11,Fernando2018_GR,Fernando:2017, khanwale2020simulating,khanwale2022fully,saurabh2021industrial,xu2021octree,khanwale2023projection}, due to its conceptual simplicity and ability to scale across a large number of processors.
		We leverage the highly optimized and scalable implementations of general finite element kernels in~\textsc{Dendrite-KT}~\citep{Ishii2019, saurabh2021industrial, saurabh2022scalable, khanwale2022fully, khanwale2022breakup, khanwale2023projection, saurabh2023scalable} built on top of the scalable \dkt{} framework.
		In this work, we use the octree-based refinement to achieve highly localized refinement. We would like to note that the framework also supports adaptive mesh refinement for extensions to velocity/vorticity-based refinement or refinement based on interfaces when the framework is extended to handle multiphase fluid flow.
		
		\item We illustrate the framework with a thorough series of test cases with increasing complexity, including electro-convective instabilities and complex geometries.
	\end{enumerate}

	The rest of the paper is structured as follows: Section 2 details the equations in non-dimensional form. Section 3 lays out the spatial and temporal discretization and associated properties. Section 4 describes the implementation and algorithmic details of the framework. Section 5 reports on the numerical results, starting from convergence tests using manufactured solutions to a full 3D simulation of a bio-diagnostic device. We conclude in~\cref{sec:conclusions}.

	
	\section{Governing equations}
	\label{sec:governing_equations}
	\subsection{Dimensional form:}
	Consider a bounded domain $\Omega \subset \mathbb{R}^d$ (for $d = 2,3$) and the time interval $[0, T]$. The domain contains a fluid (usually aqueous) with $N$ ionic species, and we seek to model the behavior of the $N$ species concentration under the imposition of electric fields. We consider the case of strong electrolytes with no reactions. The concentration distribution of the $N$ species is represented using Nernst-Plank transport equations and the Poisson equation for the electric potential, together making up the Poisson-Nernst-Plank (PNP) equations. The pressure and fluid velocity are represented using the Navier Stokes (NS) equations. There is a two-way coupling between the species concentration and the underlying fluid motion, given by an advection term in the PNP equation and an electrostatic body force term in the NS equation. 
	
	We consider $N > 1$ number of charged species with subscript $s$ indicating the species index. The species flux, $j^{s,*}_{i}$, which is driven by diffusion, migration, and convection
	~is written as:
	\begin{align}\label{eq:NP1}
		j^{s,*}_{i} = -D^s \pd{c^{s,*}}{x^*_i}  
		- D^s\frac{z^{s}F}{RT}c^{s,*}\;\pd{\phi^{s,*}}{x^*_i} 
		+  v^{*}_{i} c^{s,*} \quad\quad\text{\text{for $s^{th}$ species,}}  
	\end{align}
	$D^s$ is diffusivity of $s^{th}$ species, $c^{s,*}$ is concentration of $s^{th}$ species, $z^s$ is valance of $s^{th}$ species, $F$ is Faraday constant 96,485.33 \si{C/\mole}, $R$ is gas constant 8.314 \si{\joule/\mole.\kelvin}, $T$ is temperature, $\phi^*$ is electric potential, and $v_i^*$ is fluid velocity. 
	\begin{remark}
		Note that we use Einstein notation throughout this work; in this notation, $v_i$ represents the $i^{\text{th}}$ component of the vector $\vec{v}$, and any repeated index is implicitly summed over.  Superscript $s$ is the species index.  To avoid confusion powers of the terms are represented by writing the power outside the braces.  For example square of $s^{th}$ species of concentration is written as $\left(c^{s,*}\right)^2$.  For non-dimensional quantities, $*$ in the superscript will be dropped. 
	\end{remark}
	The time variation of concentration can be obtained by taking the divergence of species flux,
	\begin{align}
		\label{eq:NP2}
		\pd{c^{s,*}}{t^*} = - \pd{j^{s,*}_{i}}{x^*_i}, \quad\quad\text{\text{for the $s^{th}$ species.}}   
	\end{align}
	Combining \cref{eq:NP1,eq:NP2} 
	produces the Nernst-Planck equation,
	\begin{equation}\label{eq:NPF}
		\pd{c^{s,*}}{t^*} 
		+  \pd{\left(v^{*}_{i} c^{s,*}\right)}{x^*_i}
		= 
		\pd{}{x^*_i} \left(D^s\pd{c^{s,*}}{x^*_i} + D^s\frac{z^{s}F}{RT}c^{s,*}\;\pd{\phi^{s,*}}{x^*_i}\right) . 
	\end{equation}
	The second term on the left-hand side of Eqn.~\ref{eq:NPF} represents the advective transport of species, which is the coupling between the species concentration and the flow field. The first term on the right-hand side of Eqn.~\ref{eq:NPF} represents the diffusion of species, while the second term represents the electric migration of species. 
	
	The electric potential, $\phi^*$ in the domain, depends on the concentration distribution of the $N$ species and the imposed electric field. The potential is the solution to the Poisson equation~\cite{Kirby2010}  
	\begin{equation}\label{eq:Poisson}
		\pd{}{x^*_i}\left(-\varepsilon \pd{\phi^{*}}{x^*_i}\right) = \rho_e^*.
	\end{equation}
	$\varepsilon$ is the electric permittivity of fluid, and $\rho_e$ is charge density. Charge density can be expressed with species concentration.
	\begin{equation}
		\label{eq:ChargeDensity}
		\rho_e^*\ = F \sum c^{s,*} z^{s}
	\end{equation}
	Then, equation (\(\ref{eq:Poisson}\)) becomes (assuming that the permittivity, $\varepsilon$, is constant)
	\begin{equation}\label{eq:PoissonF}
		-\varepsilon\pd{}{x^*_i}\left(\pd{\phi^{*}}{x^*_i}\right) = F \sum c^{s,*} z^{s}.
	\end{equation}
	
	Finally, the Navier-Stokes equation is written as,
	\begin{equation}\label{eq:NS}
		\rho\pd{v^*_i}{t^*} + \rho\pd{\left(v^*_jv^*_i\right)}{x^*_j} + 
		\pd{p^*}{x^*_i} - \eta\pd{}{x^*_j}\left({\pd{v^*_i}{x^*_j}}\right) - \rho_e^*E_i^* = 0
	\end{equation}
	$\rho$ is the density of the fluid, $p$ is pressure, $\eta$ is dynamic viscosity, and $E_i$ is $i^{th}$ component of the electric field. The last term in the above equation is the body force due to the electric field resulting in NS being coupled with the NS and the PNP. Using the relationship between electric potential and electric field,
	\begin{equation}\label{eq:ElectricField}
		-\pd{\phi^{*}}{x^*_i} = E^*_i,
	\end{equation}
	and from the expression of charge density (equation (\ref{eq:ChargeDensity})), equation (\ref{eq:NS}) can be written as
	\begin{equation}\label{eq:NS2}
		\rho\pd{v^*_i}{t^*} + \rho\pd{\left(v^*_jv^*_i\right)}{x^*_j} + 
		\pd{p^*}{x^*_i} - \eta\pd{}{x^*_j}\left({\pd{v^*_i}{x^*_j}}\right) + F \sum c^{s,*} z^{s} \pd{\phi^{*}}{x^*_i} = 0.
	\end{equation}
	\begin{remark}
		Note that some papers in the literature use elementary charge $e$ and Boltzmann constant $k_B$ instead of Faraday constant $F$ and gas constant $R$~\citep{Mani2013}. The choice of these parameters is simply dependent on how to define concentration. Here, we choose the Faraday constant, $F$, and gas constant $R$ since molarity $M$ was used for concentration.
	\end{remark}
	We write all the governing equations together as, 
	\begin{align}
		\begin{split}
			\text{Momentum Eqns:} & \quad \rho\pd{v^*_i}{t^*} + \rho\pd{\left(v^*_jv^*_i\right)}{x^*_j} + 
			\pd{p^*}{x^*_i} - \eta\pd{}{x^*_j}\left({\pd{v^*_i}{x^*_j}}\right) 
			+ F \left(\sum c^{s,*} z^{s}\right) \pd{\phi^{*}}{x^*_i} = 0, \label{eqn:nav_stokes} 
		\end{split} \\
		\text{Solenoidality:} & \quad \pd{v^{*}_i}{x^*_i} = 0, \label{eqn:dim_cont} \\
		\text{Poisson:} & \quad -\varepsilon\pd{}{x^*_i}\left(\pd{\phi^{*}}{x^*_i}\right) = F \sum c^{s,*} z^{s},\label{eqn:dim_poisson}\\
		\text{Nernst-Planck:} & \quad \pd{c^{s,*}}{t^*} 
		+ v^{*}_{i} \pd{c^{s,*}}{x^*_i}
		=  \pd{}{x^*_i} \left(D^s\pd{c^{s,*}}{x^*_i} + D^s\frac{z^{s}F}{RT}c^{s,*}\;\pd{\phi^{*}}{x^*_i}\right). \label{eqn:dim_np}
	\end{align}
	\Cref{tab:DParameters} shows all the dimensional parameters involved. 
	\begin{table}[H]
		\centering
		\begin{tabular}{|c||c|c|}
			\hline
			& & \\ [-10pt]
			\textbf{Parameters} & \textbf{Units} & \textbf{Value/Formula}\\ 
			& & \\ [-10pt]
			\hline
			\hline
			Species diffusivity, $D^s$ & \(m^2/s\) & \\
			\hline
			Species valence, $z^s$ &  & \\
			\hline
			Faraday constant, $F$ & $C/mol$ & 96,485.33 $C/mol$\\
			\hline
			Gas constant, $R$ & $J/K\cdot mol$ & 8.314 $J/K\cdot mol$\\
			\hline
			Temperature, $T$ & $K$ & \\
			\hline
			Vacuum permittivity, $\varepsilon_0$ & $F/m, \: C/V\cdot m, \: N/V^2$ & 8.854$\times10^{-12}  F/m$\\
			\hline
			Relative permittivity, $\varepsilon_r$ &  & \\
			\hline
			Permittivity, $\varepsilon$ & $F/m, \: C/V\cdot m, \: N/V^2$ & $\varepsilon = \varepsilon_r\varepsilon_0$\\
			\hline
			Charge density, $\rho_e$ & $C$ & \\
			\hline
			Density, $\rho$ & $kg/m^3$ & \\
			\hline
			Viscosity, $\eta$ & $Pa \cdot s, \: cP$ & 8.90 $\times 10^{-4} P\cdot s$ for water at \SI{25}{\celsius}\\
			\hline
			Ionic strength of bulk, $I_b$ & $M, \: mol/m^3$ & $I_b=\frac{1}{2}\Sigma z_i^2 c_i$ \\
			\hline
			Thermal voltage, $V_T$ & $V$ & $V_T=RT/F$\\
			\hline\end{tabular}
		\caption{Dimensional parameters in the governing equations.}
		\label{tab:DParameters}
	\end{table}
	
	\subsection{Non-dimensional form:}
	Each variable is scaled with its associated characteristic quantity to define its non-dimensional counterpart.
	
	\begin{equation} \label{eq:scaling1}
		x_i = \frac{x_i^*}{L_{ch}}, \quad t=\frac{t^*}{t_{ch}}, \quad
		v_i=\frac{v_i^*}{U_{ch}}, 
		\quad p=\frac{p^*}{p_{ch}}, \quad c_i=\frac{c_i^*}{c_{ch}}, \quad \phi=\frac{\phi^*}{\phi_{ch}}.
	\end{equation}
	And for the derivative terms,
	\begin{equation} \label{eq:scaling2}
		\pd{}{x_i}=L_{ch}\pd{}{x^*_i}, \quad \pd{}{x_i}\left(\pd{}{x_i}\right)=L_{ch}^2\pd{}{x^*_i}\left(\pd{}{x^*_i}\right)
	\end{equation}
	To non-dimensionalize the governing equations, it is important to select appropriate reference quantities in \cref{eq:scaling1}. 
	The reference concentration and reference potential are selected as follows:
	\begin{equation} \label{eq:refIPhi}
		c_{ch}=I_b=\frac{1}{2}\sum \left(z^s\right)^2 c^s_{initial} \quad \text{and} \quad \phi_{ch} = V_T=\frac{RT}{F} 
	\end{equation}
	$I_b$ is the ionic strength of bulk fluid at a reference time (initial), and $V_T$ is the thermal voltage. The rest of the variables have multiple options for the reference quantity.  We choose the following reference quantities: The characteristic velocity $U_{ch}$ is given by $D/L_{ch}$, characteristic time $t_{ch}$ is given by $L_{ch}^2/D$ and characteristic pressure given by $\eta D/ L_{ch}^2$. Here, $D$ is the mean species diffusivity ($D = 1/N\sum^N_{s=1} D^s$), and $\eta$ is the dynamic viscosity of the fluid medium. For the cases in this paper, we assume $D^s$ is the same for all species; therefore, the non-dimensional $D^s$ equals 1. 
	
	With these characteristic quantities, the non-dimensional equations are given as follows, 
	
	\begin{align}
		\begin{split}
			\text{Momentum Eqns:} & \quad \frac{1}{Sc}\left(\pd{v_i}{t} + \pd{\left(v_jv_i\right)}{x_j}\right) + 
			\pd{p}{x_i} - \pdd{v_i}{x_j} 
			+ \frac{\kappa}{2 \Lambda^2} \left(\sum c^{s} z^{s}\right) \pd{\phi}{x_i} = 0, \label{eqn:nav_stokes_nondim} 
		\end{split} \\
		\text{Solenoidality:} & \quad \pd{v_i}{x_i} = 0, \label{eqn:nondim_cont} \\
		\text{Poisson:} & \quad -2 \Lambda^2 \pdd{\phi}{x_i} = \sum c^{s} z^{s},\label{eqn:nondim_poisson}\\
		\text{Nernst-Planck:} & \quad \pd{c^{s}}{t} 
		+ v_{i} \pd{c^{s}}{x_i}
		=  \pd{}{x_i} \left(\pd{c^{s}}{x_i} + z^{s}c^{s}\;\pd{\phi}{x_i}\right). \label{eqn:nondim_np}
	\end{align}
	Where the Schmidt number, $Sc = \frac{\eta}{\rho D}$, the electrohydrodynamic (EHD) coupling constant, $\kappa = \frac{\varepsilon}{\eta D} \left(\frac{RT}{F}\right)^2$, and non-dimensionalized Debye length, $\Lambda = \lambda_D/L_{ch}$, where $\lambda_D$ is the Debye length,~$
	\lambda_{D} = \sqrt{\frac{1}{2}\frac{\varepsilon RT}{F^{2}I_{b}}}$.

	\section{Numerical method and its properties}
	\label{sec:numerical_tecniques}
	
	We use a second-order backward difference formula (BDF) scheme as a time-stepping strategy for~\crefrange{eqn:nav_stokes_nondim}{eqn:nondim_np}.  We solve momentum equations and Poisson Nernst-Planck (PNP) equations in two separate blocks.   
	
	Let $\delta t$ be a time-step; let $t^k:= k \delta t$; let $n+1$ be the current time-point, and $n$ be the previous time-point. We define an intermediate velocity at time $n+1$ as follows: 
	\begin{equation}
		\label{eqn:averaging1}
		\widetilde{\vec{v}}^{k+1} := 2\vec{v}^{k} - \vec{v}^{k-1}, 
	\end{equation}
	We define the time-discretized weak form of
	the Navier-Stokes-Poisson-Nernst-Plank (NS-PNP) equations are as follows: 
	
	
	
	%
	
	\begin{definition}{}{variational_form_sem_disc}
		Let $(\cdot,\cdot)$ be the standard $L^2$ inner product. We state the time-discrete variational problem as follows: find $\vec{v}^{k+1}(\vec{x}) \in \vec{H}^r(\Omega)$, $p^{k+1}(\vec{x})$, $\phi^{k+1}(\vec{x})$, $c^{s,k+1}(\vec{x})$ $\in {H}^r(\Omega)$ such that
		\begin{align}
			\begin{split}
				\text{Momentum Eqns:} & \quad \frac{1}{Sc}\left(w_i, \, \frac{\beta_0v^{k+1}_i + \beta_1v^k_i + \beta_2v^{k-1}_i}{\delta t}\right)_\Omega 
				+ \frac{1}{Sc}\left(w_i, \, \, \tvj^{k+1} \, \pd{v_i^{k+1}}{x_j}\right)_\Omega \\ & \quad 
				+ \left(w_i, \, \pd{p^{k+1}}{x_i} \right)_\Omega 
				+ \left(\pd{w_i}{x_j}, \, {\pd{v_i^{k+1}}{x_j}} \right)_\Omega
				+ \frac{\kappa}{2 \Lambda^2} \left(w_i,\left(\sum c^{s, k+1} z^{s}\right) \pd{\phi^{k+1}}{x_i}\right)_{\Omega}\\ & \quad
				- \left(w_i, \, \hat{n}_j {\pd{v_i^{k+1}}{x_j}} \right)_{\partial\Omega} 
				= 0, 
				\label{eqn:nav_stokes_var_semi_disc}
			\end{split} \\ 
			\text{Solenoidality:} & \quad 
			\quad \left(q, \, \pd{v^{k+1}_{i}}{x_i}\right)_\Omega = 0,
			\label{eqn:cont_var_semi_disc} \\
			\text{Poisson:} & \quad 
			2\Lambda^2 \left(\pd{q}{x_i}, \pd{\phi^{k+1}}{x_i}\right)_\Omega 
			- 2\Lambda^2 \left(q, n_i\pd{\phi^{k+1}}{x_i}\right)_{\partial\Omega}
			-  \left(q,  \sum c^{s,k+1} z^{s}\right)_\Omega   = 0, \label{eqn:poisson_eqn_var_semi_disc}\\
			\begin{split}
				\text{Nernst-Planck:} & \quad \left(q, \frac{\beta_0 c^{s,k+1} + \beta_1 c^{s,k} + \beta_2 c^{s,k-1}}{\delta t} \right)_\Omega 
				+ \left(q, \, \tvi^{k+1} \pd{c^{s, k+1}}{x_i} \right)_\Omega \\  & \quad 
				+ \left(\pd{q}{x_i}, \, \pd{c^{s, k+1}}{x_i} + z^{s}c^{s, k+1}\;\pd{\phi^{k+1}}{x_i} \right)_\Omega \\  & \quad
				- \left(q, \, \hat{n}_i \left(\pd{c^{s, k+1}}{x_i} + z^{s}c^{s, k+1}\;\pd{\phi^{k+1}}{x_i}\right) \right)_{\partial\Omega} = 0,
				\label{eqn:np_eqn_var_semi_disc}
			\end{split} 
		\end{align}
		$\forall \vec{w} \in \vec{H}^r(\Omega)$, $\forall q \in H^r(\Omega)$, given $\vec{v}^{k},\;\vec{v}^{k-1} \in \vec{H}^r(\Omega)$, and $\phi^{k},\phi^{k-1}, c^{s,k}, c^{s,k} \in H^r(\Omega)$.
		$\hat{\vec{n}}$ is the outward pointing normal to the boundary $\partial \Omega$
	\end{definition}

	\begin{remark}
		Note that the variational forms above are for general Hilbert spaces with regularity $r$ (order of allowed derivatives in the space). However, for the cases in this paper, we restrict ourselves to $r=1$, which corresponds to bilinear basis functions. 
	\end{remark}
	
	The semi-implicit discretization for~\cref{eqn:nav_stokes_var_semi_disc} is useful as it naturally linearizes the momentum equations.  This discretization is particularly useful for chaotic systems where the timesteps are naturally small to capture small-scale flow structures.  
	Solving this linearized semi-implicit system helps with performance as opposed to the fully implicit discretization, which requires a Newton iteration and multiple sub-iterations (each iteration being a linear solve).

	\subsection{Spatial discretization and the variational multiscale approach}
	\label{subsec:space_scheme}
	We use continuous Galerkin finite elements for spatial discretization of the time-discrete system presented in \cref{defn:variational_form_sem_disc}.  However, imposing the d'Alembert constraint (\cref{eqn:cont_var_semi_disc}) when solving $\vec{v}, \, p$ together presents a challenge.  An equal order polynomial approximation of velocity and pressure leads to pressure instability because of the violation of the discrete inf-sup condition (Ladyzhenskaya-Babuska-Brezzi condition, e.g., see \citet[page 31]{ Volker2016}).  There are two methods to circumvent this issue: 1) add pressure stabilization to~\cref{eqn:cont_var_semi_disc}, which converts the inf-sup stability restriction to a V-elliptic coercivity problem~\citep{article:TezMitRayShi92}, thereby curbing the artificial pressure instabilities; 2) The second method is to solve a separate pressure Poisson equation using the Helmholtz Hodge decomposition of solenoidal velocity fields.  For our specific problem of NS-PNP equations, we use the first method.  This method allows us to use the Variational multi-scale method (VMS)~\citep{ hughes2018multiscale} to stabilize pressure.  In addition to providing pressure stabilization, VMS also has the following useful features. It allows us to naturally perform Large Eddy simulation type decomposition and appropriate fine scale modeling~\citep{Hughes1995}.  The VMS enrichment for modeling fine-scale velocities, pressures, and ion concentrations allows grid coarsening away from the extremely fine resolutions near the boundary needed for the EDLs.  Note that providing targeted resolution near the boundary for resolving EDLs is relatively simpler; however, resolving emerging small-scale structures due to electrochemical instabilities for high electric fields is extremely difficult without doing adaptive meshing~\footnote{Adaptive mesh refinement to capture small-scale instabilities also presents its challenges, primarily of constructing \textit{a}posteriori error estimates for solutions to Navier-Stokes which is an open research field}.  Therefore, the subgrid-scale modeling away from the fine resolutions of the EDLs is valuable for accurately capturing electroconvective instabilities with relatively coarser resolutions. However, in this manuscript, we focus on much finer resolutions, where the primary contribution of VMS is to circumvent the inf-sup condition and provide pressure stabilization. 
	
	The VMS approach uses a direct-sum decomposition of the function spaces as follows.  If $\vec{v} \in \vec{V}$, $p \in Q$, and $\phi \in Q$ then we decompose these spaces as: 
	\begin{eqnarray}
		\vec{V} = \vec{V}^{c} \oplus \vec{V}^{f} \qquad \text{and} \qquad Q = Q^{c} \oplus Q^{f},
	\end{eqnarray}
	where $\overline{\vec{V}}$ and $\overline{Q}$ are the cG(r) subspaces of $\vec{V}$ and $Q$, respectively, and the primed versions are the complements of the cG(r) subspaces in $\vec{V}$ and $Q$, respectively.
	We decompose the velocity and pressure as follows: 
	\begin{equation}
		\vec{v} = \vec{v}^{c} + \vec{v}^{f}, \quad
		\phi = \phi^{c} + \phi^{f}, \quad
		c^{s} = c^{s,c} + c^{s,f},
		\quad \text{and}
		\quad p = p^{c} + p^{f},
	\end{equation} 
	where the {\it coarse scale} components are
	$\vec{v}^{c} \in \vec{V}^{c}$ and  $p^{c}, \phi^{c}, c^{s,c} \in Q^{c}$, and the {\it fine scale} components are
	$\vec{v}^{f} \in \vec{V}^{f}$ and $p^{f},\phi^{f} \in Q^{f}$.  We define a projection operator, $\mathscr{P}:\vec{V} \rightarrow \vec{V}^{c}$, such that
	$\vec{v}^{c} = \mathscr{P}\{\vec{v}\}$ and $\vec{v}^{f} = \vec{v} - \mathscr{P}\{\vec{v}\}$.  Similar operators decompose the other variables, $p$, $\phi$, $c^{s}$.  
	
	Substituting these decompositions in the original variational form in~\cref{defn:variational_form_sem_disc} yields:
	\begin{align}	
		\begin{split}	
			\text{Momentum Eqns:}& \quad \frac{1}{Sc}\left(w_i,\; \frac{\beta_0v^{c,k+1}_i + \beta_1v^{c,k}_i + \beta_2v^{c,k-1}_i}{\delta t} \right)_{\Omega}\\ 
			& 
			+  \frac{1}{Sc}\left(w_i, \tvj^{k+1}\pd{v_i^{c,k+1}}{x_j} \right)_{\Omega} 
			- \frac{1}{Sc}\left(\pd{w_i}{x_j}, \tvj^{k+1}v_i^{f,k+1} \right)_{\Omega} \\
			& +  \left(w_i,\pd{\left(p^{c,k+1} + p^{f,k+1}\right)}{x_i}\right)_{\Omega} 
			+ \left(\pd{w_i}{x_j}, \pd{v_i^{c,k+1}}{x_j}\right)_{\Omega}
			+ \boxed{\left(\pd{w_i}{x_j}, \pd{v_i^{f,k+1}}{x_j}\right)_{\Omega}} \\
			&+ \frac{\kappa}{2 \Lambda^2} \left(w_i, \left(\sum \left(c^{s, c, k+1} + c^{s, f, k+1} \right)z^{s}\right) \pd{\phi^{c, k+1}}{x_i} \right)_{\Omega} \\
			&\boxed{+ \frac{\kappa}{2 \Lambda^2} \left(w_i, \left(\sum \left(c^{s, c, k+1} + c^{s, f, k+1} \right)z^{s}\right) \pd{\phi^{f, k+1}}{x_i}\right)_{\Omega}} = 0,
			\label{eqn:weak_VMS_ns}
		\end{split}\\
		\begin{split}
			\text{Solenoidality:}
			& \quad  \left(q,\pd{v_i^{c,k+1}}{x_i}\right)_{\Omega} + \left(q,\pd{v_i^{f,k+1}}{x_i}\right)_{\Omega} = 0,\label{eqn:cont_eqn_var__vms1}
		\end{split}\\
		\begin{split}
			\text{Poisson:} & \quad 
			2\Lambda^2 \left(\pd{q}{x_i}, \pd{ \left(\phi^{c,k+1} + \phi^{f,k+1}\right) }{x_i}\right)_\Omega 
			+ 2\Lambda^2 \left(q, n_i\pd{\left(\phi^{c,k+1} + \phi^{f,k+1}\right)}{x_i}\right)_{\partial\Omega} \\ 
			&  
			- \left(q,  \sum c^{s,c,k+1} z^{s}\right)_\Omega  
			- \left(q,  \sum c^{s,f,k+1} z^{s}\right)_\Omega = 0, \label{eqn:poisson_eqn_var__vms1}
		\end{split}\\
		\begin{split}
			\text{Nernst-Planck:} & \quad \left(q, \frac{\beta_0 c^{s,c,k+1} + \beta_1 c^{s,c,k} + \beta_2 c^{s,c,k-1}}{\delta t} \right)_\Omega 
			+ \left(q, \, \tvi^k \pd{\left(c^{s, c, k+1} + c^{s, f, k+1}\right)}{x_i} \right)_\Omega \\  & \quad 
			+ \left(\pd{q}{x_i}, \, \pd{\left(c^{s, c, k+1} + c^{s, f, k+1}\right)}{x_i} 
			+ z^{s}\left(c^{s, c, k+1} + c^{s, f, k+1}\right)\;\pd{\left(\phi^{c, k+1} + \phi^{f, k+1}\right)}{x_i} \right)_\Omega \\  & \quad
			- \left(q, \, \hat{n}_i \left(\pd{\left(c^{s, c, k+1} + c^{s, f, k+1}\right)}{x_i} 
			+ z^{s}\left(c^{s, c, k+1} + c^{s, f, k+1}\right)\;\pd{\left(\phi^{c, k+1} + \phi^{f, k+1}\right)}{x_i} \right) \right)_{\partial\Omega} = 0,
			\label{eqn:np_eqn_var_vms1}
		\end{split}
	\end{align} 
	where $\vec{w},\vec{v}^{c}, \in \mathscr{P}\vec{H}^{r}(\Omega)$, $p^{c},\;\phi^{c} \in \mathscr{P}H^r(\Omega), \vec{v}^{f} \in (\mathscr{I} - \mathscr{P})\vec{H}^{r}(\Omega)$,
	$\phi^{f}, c^{s,f}, p^{f} \in (\mathscr{I} - \mathscr{P})H^r(\Omega)$, and 
	$c^{s, c}, q \in \mathscr{P}H^r(\Omega)$.
	Here $\mathscr{I}$ is the identity operator, and $\mathscr{P}$ is the projection operator.  
	
	We can further simplify \cref{eqn:np_eqn_var_vms1}, 
	
	\begin{align}
		\begin{split}
			& \left(q, \frac{\beta_0 c^{s,c,k+1} + \beta_1 c^{s,c,k} + \beta_2 c^{s,c,k-1}}{\delta t} \right)_\Omega 
			+ \left(q, \, \tvi^k \pd{\left(c^{s, c, k+1} + c^{s, f, k+1}\right)}{x_i} \right)_\Omega \\  & \quad 	
			+ \left(\pd{q}{x_i}, \, \pd{\left(c^{s, c, k+1} + c^{s, f, k+1}\right)}{x_i} 
			+ z^{s}\left(c^{s, c, k+1} + c^{s, f, k+1}\right)\;\pd{\left(\phi^{c, k+1} + \phi^{f, k+1}\right)}{x_i} \right)_\Omega \\  & \quad
			- \left(q, \, \hat{n}_i \left(\pd{\left(c^{s, c, k+1} + c^{s, f, k+1}\right)}{x_i} 
			+ z^{s}\left(c^{s, c, k+1} + c^{s, f, k+1}\right)\;\pd{\left(\phi^{c, k+1} + \phi^{f, k+1}\right)}{x_i} \right) \right)_{\partial\Omega} = 0,
			\label{eqn:np_eqn_var_vms_s1}
		\end{split}\\
		\begin{split}
			\implies 
			& \quad 
			\left(q, \frac{\beta_0 c^{s,c,k+1} + \beta_1 c^{s,c,k} + \beta_2 c^{s,c,k-1}}{\delta t} \right)_\Omega 
			+ \left(q, \, \tvi^k \pd{c^{s, c}}{x_i} \right)_\Omega 
			+ \left(q, \, \tvi^k \pd{c^{s, f}}{x_i} \right)_\Omega \\  & \quad
			+ \left(\pd{q}{x_i}, \, \pd{c^{s, c, k+1}}{x_i} \right)_\Omega
			+ \boxed{\left(\pd{q}{x_i}, \, \pd{c^{s, f, k+1}}{x_i} \right)_\Omega} \\  & \quad
			+ \left(\pd{q}{x_i}, \, z^{s}c^{s, c, k+1}\;\pd{\phi^{c, k+1}}{x_i} \right)_\Omega
			+ \left(\pd{q}{x_i}, \, z^{s}c^{s, f, k+1}\;\pd{\phi^{c, k+1}}{x_i} \right)_\Omega \\  & \quad
			+ \boxed{\left(\pd{q}{x_i}, \, z^{s}c^{s, c, k+1}\;\pd{\phi^{f, k+1}}{x_i} \right)_\Omega}
			+ \boxed{\left(\pd{q}{x_i}, \, z^{s}c^{s, f, k+1}\;\pd{\phi^{f, k+1}}{x_i} \right)_\Omega} 
			\\  & \quad
			- \left(q, \, n_i \pd{c^{s, c, k+1}}{x_i} \right)_{\partial\Omega}
			- \boxed{\left(q, \, n_i \pd{c^{s, f, k+1}}{x_i} \right)_{\partial\Omega}} \\  & \quad
			- \left(q, \, n_i z^{s}c^{s, c, k+1}\;\pd{\phi^{c, k+1}}{x_i} \right)_{\partial\Omega}
			- \boxed{\left(q, \, n_i z^{s}c^{s, f, k+1}\;\pd{\phi^{c, k+1}}{x_i} \right)_{\partial\Omega}}\\  & \quad
			- \boxed{\left(q, \, n_i z^{s}c^{s, c, k+1}\;\pd{\phi^{f, k+1}}{x_i} \right)_{\partial\Omega}}
			- \boxed{\left(q, \, n_i z^{s}c^{s, f, k+1}\;\pd{\phi^{f, k+1}}{x_i} \right)_{\partial\Omega}}\\
			& \quad = 0,
			\label{eqn:np_eqn_var_vms_s2}
		\end{split}
	\end{align}
	We set the boxed terms to zero (see next paragraph, and also, Table.~\ref{tab:models}). These assumptions reduce the Nernst-Plank equations to the following, 
	\begin{align}
		\begin{split}
			\text{Nernst-Planck:} 
			& \quad 
			\left(q, \frac{\beta_0 c^{s,c,k+1} + \beta_1 c^{s,c,k} + \beta_2 c^{s,c,k-1}}{\delta t} \right)_\Omega 
			+ \left(q, \, \tvi^k \pd{c^{s, c, k+1}}{x_i} \right)_\Omega 
			+ \left(q, \, \tvi^k \pd{c^{s, f, k+1}}{x_i} \right)_\Omega\\  & \quad
			+ \left(\pd{q}{x_i}, \, \pd{c^{s, c, k+1}}{x_i} \right)_\Omega
			+ \left(\pd{q}{x_i}, \, z^{s}c^{s, c, k+1}\;\pd{\phi^{c, k+1}}{x_i} \right)_\Omega   
			+ \left(\pd{q}{x_i}, \, z^{s}c^{s, f, k+1}\;\pd{\phi^{c, k+1}}{x_i} \right)_\Omega 
			\\  & \quad
			- \left(q, \, n_i \pd{c^{s, c, k+1}}{x_i} \right)_{\partial\Omega}
			- \left(q, \, n_i z^{s}c^{s, c, k+1}\;\pd{\phi^{c, k+1}}{x_i} \right)_{\partial\Omega} = 0
			\label{eqn:np_eqn_var_vms_app_s1}
		\end{split} \\
		\begin{split}
			\implies
			& \quad 
			\left(q, \frac{\beta_0 c^{s,c,k+1} + \beta_1 c^{s,c,k} + \beta_2 c^{s,c,k-1}}{\delta t} \right)_\Omega 
			+ \left(q, \, \tvi^k \pd{c^{s, c, k+1}}{x_i} \right)_\Omega 
			- \left(\pd{q}{x_i}, \, \tvi^k c^{s, f, k+1} \right)_\Omega \\  & \quad
			+ \left(\pd{q}{x_i}, \, \pd{c^{s, c, k+1}}{x_i} \right)_\Omega
			+ \left(\pd{q}{x_i}, \, z^{s}c^{s, c, k+1}\;\pd{\phi^{c, k+1}}{x_i} \right)_\Omega   
			+ \left(\pd{q}{x_i}, \, z^{s}c^{s, f, k+1}\;\pd{\phi^{c, k+1}}{x_i} \right)_\Omega
			\\  & \quad
			+  \boxed{\left(q, \, n_i \tvi^k c^{s, f, k+1} \right)_{\partial\Omega}}
			- \left(q, \, n_i \pd{c^{s, c, k+1}}{x_i} \right)_{\partial\Omega}
			- \left(q, \, n_i z^{s}c^{s, c, k+1}\;\pd{\phi^{c, k+1}}{x_i} \right)_{\partial\Omega} = 0.
			\label{eqn:np_eqn_var_vms_app_s2}
		\end{split}
	\end{align}
	
	Table ~\ref{tab:models} lists the various boxed terms in the coupled set of NS-PNP equations that are set to zero, along with the rationale for setting them to zero.  
	
	
	\begin{table}[H]
		\centering
		\centering\footnotesize\setlength\tabcolsep{4pt}
		\begin{tabular}{|>{\centering\arraybackslash}m{0.2\linewidth}|>{\centering\arraybackslash}m{0.4\linewidth}|>{\centering\arraybackslash}m{0.4\linewidth}|}
			\toprule
			\toprule
			Equation & Terms & Rationale for assumption \\
			\midrule
			\midrule
			Momentum equations & $\left(\pd{w_i}{x_j}, \pd{v_i^{f,k+1}}{x_j}\right)_\Omega$ &  \textit{Projection:} Orthogonality  conditions from the projection decomposition, as $w_i \in \mathscr{P}\vec{H}^{r}(\Omega)$ and $v_i^{f,k+1} \in (\mathscr{I} - \mathscr{P})\vec{H}^{r}(\Omega)$~\citep{Bazilevs2007, Hughes2000, hughes2018multiscale}\\
			\midrule
			Momentum equations & $\frac{\kappa}{2 \Lambda^2} \left(w_i, \left(\sum \left(c^{s, c, k+1} + c^{s, f, k+1} \right)z^{s}\right) \pd{\phi^{f, k+1}}{x_i}\right)_\Omega$ & $\phi^{f, k+1}$ assumed to be zero\\
			\midrule
			Nernst-Planck equations & $\left(q, \, n_i \tvi^k c^{s, f, k+1} \right)_{\partial\Omega}$ & Depends on the boundary condition, could be zero by either a Dirichlet condition of velocity or a Dirichlet condition of concentration \\
			\midrule		
			Nernst-Planck equations & $\left(\pd{q}{x_i}, \, \pd{c^{s, f, k+1}}{x_i} \right)_\Omega$ &  \textit{Projection:} Orthogonality  conditions from the projection decomposition, as $q \in \mathscr{P}H^{r}(\Omega)$ and $c^{s,f,k+1} \in (\mathscr{I} - \mathscr{P})H^{r}(\Omega)$~\citep{Bazilevs2007, Hughes2000, hughes2018multiscale}\\
			\midrule
			Nernst-Planck equations & $\left(\pd{q}{x_i}, \, z^{s}c^{s, c, k+1}\;\pd{\phi^{f, k+1}}{x_i} \right)_\Omega$ & $\phi^{f, k+1}$ assumed to be zero \\
			\midrule
			Nernst-Planck equations & $\left(\pd{q}{x_i}, \, z^{s}c^{s, f, k+1}\;\pd{\phi^{f, k+1}}{x_i} \right)_\Omega$ & $\phi^{f, k+1}$ assumed to be zero \\
			\midrule
			Nernst-Planck equations & $\left(q, \, n_i \pd{c^{s, f, k+1}}{x_i} \right)_{\partial\Omega}$ & flux of $c^{s, f, k+1}$  assumed to be zero on boundaries   \\
			\midrule
			Nernst-Planck equations & $\left(q, \, n_i z^{s}c^{s, f, k+1}\;\pd{\phi^{c, k+1}}{x_i} \right)_{\partial\Omega}$ & $c^{s,f, k+1}$ assumed to be zero at the boundary \\
			\midrule
			Nernst-Planck equations & $\left(q, \, n_i z^{s}c^{s, c, k+1}\;\pd{\phi^{f, k+1}}{x_i} \right)_{\partial\Omega}$ & $\phi^{f, k+1}$ assumed to be zero \\
			\midrule
			Nernst-Planck equations & $\left(q, \, n_i z^{s}c^{s, f, k+1}\;\pd{\phi^{f, k+1}}{x_i} \right)_{\partial\Omega}$ & $c^{s,f, k+1}$ assumed to be zero at the boundary \\
			\midrule
			Poisson equations & $2\Lambda^2 \left(\pd{q}{x_i}, \pd{ \phi^{f,k+1} }{x_i}\right)_\Omega$ & $\phi^{f, k+1}$ assumed to be zero \\
			\midrule
			Poisson equations & $2\Lambda^2 \left(q, n_i\pd{\phi^{f,k+1}}{x_i}\right)_{\partial\Omega}$ & $\phi^{f, k+1}$ assumed to be zero \\
			\bottomrule
			\bottomrule
		\end{tabular}
		\caption{Fine scale terms assumed to be zero and along with the rationale.}
		\label{tab:models}
	\end{table}
	With all the simplifications listed in Table~\ref{tab:models} we can now write the fully discrete variational form for the NS-PNP equations. 
	\begin{definition}{}{variational_form_sem_disc_vms1}
		Let $(\cdot,\cdot)$ be the standard $L^2$ inner product. We state the time-discrete variational problem as follows: find 
		$\vec{v}^{c, k+1}, \in \mathscr{P}\vec{H}^{r}(\Omega)$, and $p^{c, k+1},\;c^{s,c, k+1}(\vec{x}), \phi^{c, k+1} \in \mathscr{P}H^r(\Omega)$
		such that
		\begin{align}
			\begin{split}	
				\text{Momentum Eqns:}& \quad \frac{1}{Sc}\left(w_i,\; \frac{\beta_0v^{c,k+1}_i + \beta_1v^{c,k}_i + \beta_2v^{c,k-1}_i}{\delta t} \right)_{\Omega}\\ 
				& 
				+  \frac{1}{Sc}\left(w_i, \tvj^{k+1}\pd{v_i^{c,k+1}}{x_j} \right)_{\Omega} 
				- \frac{1}{Sc}\left(\pd{w_i}{x_j}, \tvj^{k+1}v_i^{f,k+1} \right)_{\Omega} \\
				& +  \left(w_i,\pd{\left(p^{c,k+1} + p^{f,k+1}\right)}{x_i}\right)_{\Omega} 
				+ \left(\pd{w_i}{x_j}, \pd{v_i^{c,k+1}}{x_j}\right)_{\Omega} \\
				&+ \frac{\kappa}{2 \Lambda^2} \left(w_i, \left(\sum c^{s, k+1}z^{s}\right) \pd{\phi^{c, k+1}}{x_i}\right)_{\Omega}= 0, 
				\label{eqn:nav_stokes_var_semi_disc_vms1}
			\end{split}\\
			\text{Solenoidality:} & \quad 
			\quad \left(q,\pd{v_i^{c,k+1}}{x_i}\right)_{\Omega} + \left(\pd{q}{x_i},v_i^{f,k+1}\right)_{\Omega} = 0,
			\label{eqn:cont_var_semi_disc_vms1} \\
			\begin{split}
				\text{Poisson:} & \quad 
				2\Lambda^2 \left(\pd{q}{x_i}, \pd{\phi^{c,k+1}}{x_i}\right)_\Omega 
				- 2\Lambda^2 \left(q, n_i\pd{\phi^{c, k+1}}{x_i}\right)_{\partial\Omega}\\
				& \quad 
				-  \left(q,  \sum c^{s,c, k+1} z^{s}\right)_\Omega   
				- \left(q,  \sum c^{s,f, k+1} z^{s}\right)_\Omega= 0, \label{eqn:poisson_eqn_var_semi_disc_vms1}
			\end{split}\\
			\begin{split}
				\text{Nernst-Planck:} 
				& \quad 
				\left(q, \frac{\beta_0 c^{s,c,k+1} + \beta_1 c^{s,c,k} + \beta_2 c^{s,c,k-1}}{\delta t} \right)_\Omega \\ & \quad 
				+ \left(q, \, \tvi^k \pd{c^{s, c, k+1}}{x_i} \right)_\Omega 
				- \left(\pd{q}{x_i}, \, \tvi^k c^{s, f, k+1} \right)_\Omega \\  & \quad
				+ \left(\pd{q}{x_i}, \, z^{s}c^{s, c, k+1}\;\pd{\phi^{c, k+1}}{x_i} \right)_\Omega   
				+ \left(\pd{q}{x_i}, \, z^{s}c^{s, f, k+1}\;\pd{\phi^{c, k+1}}{x_i} \right)_\Omega
				\\  & \quad
				- \left(q, \, n_i z^{s}c^{s, c, k+1}\;\pd{\phi^{c, k+1}}{x_i} \right)_{\partial\Omega} 
				\\  & \quad
				+ \left(\pd{q}{x_i}, \, \pd{c^{s, c, k+1}}{x_i} \right)_\Omega 
				- \left(q, \, n_i \pd{c^{s, c, k+1}}{x_i} \right)_{\partial\Omega}
				= 0.
				\label{eqn:np_eqn_var_semi_disc_vms1}
			\end{split} 
		\end{align}
		$\forall \vec{w} \in \mathscr{P}\vec{H}^{r}(\Omega)$, $q \in \mathscr{P}H^r(\Omega)$\\
		, and $\vec{v}^{f, k+1} \in (\mathscr{I} - \mathscr{P})\vec{H}^{r}(\Omega)$,
		$c^{s,f, k+1}, p^{f, k+1} \in (\mathscr{I} - \mathscr{P})H^r(\Omega)$\\
		, given $\vec{v}^{k},\;\vec{v}^{k-1} \in \vec{H}^r(\Omega)$, and $\phi^{k},\phi^{k-1}, c^{s,k}, c^{s,k} \in H^r(\Omega)$.
	\end{definition}
	Note that the above variational problem is not closed, as the fine-scale velocities, pressure, and ion concentrations are unknown.  
	To close the system, we use the \textit{residual-based approximation} proposed in~\citet{ Bazilevs2007} for the fine-scale components. Defining $\mathcal{R}_{m}$, $\mathcal{R}_{sol}$, and $\mathcal{R}_{c^s}$ as the residuals of the momentum, solenoidality and Nernst-Plank equations, the fine-scale components are:   
	\begin{equation}
		v_i^{f} = -\tau_m \mathcal{R}_m\left(v_i^{c}, p^{c}, c^{s,c}\right), \qquad 
		p^{f} = - \tau_{sol} \mathcal{R}_{sol}(v_i^{c}), \qquad \text{and} \qquad
		c^{s, f} = - \tau_{c^{s}} \mathcal{R}_{c^{s}}\left(c^{s, c}, v_i^{c}, \phi^c\right),
	\end{equation}
	where,
	\begin{align}
		\begin{split}
			\tau_m &= \left( \frac{4}{\Delta t^2}  + v_i^{c}G_{ij}v_j^{c} + C_{I} \left(\frac{1}{Re}\right)^2 G_{ij}G_{ij}\right)^{-1/2},
		\end{split}\\
		\begin{split}
			\tau_{sol} &=  \frac{1}{tr(G_{ij})\tau_m},
		\end{split}\\
		\tau_{c^s} &=  \left( \frac{4}{\Delta t^2}  + v_i^{c}G_{ij}v_j^{c} + C_{I} G_{ij}G_{ij}\right)^{-1/2}.
	\end{align}
	
	Here, we set $C_{I}$ and $C_{\phi}$ for all our simulations to 6, and $G_{ij} = \sum_{k=1}^{3} \partial \xi_k/ \partial x_i \partial \xi_k/ \partial x_j$ is a mesh-based tensor that accounts for the inverse mapping ($\partial \xi_k/ \partial x_i$) between the parametric and the physical domain of the element. For octrees with an equal aspect ratio ($\Delta x = \Delta y = \Delta z = h $), such a transformation simplifies $G_{ij}$ to $(2/h)^2$, where $h$ is the element's size.   
	Subsequently, the final fully-discrete variational form that is solved can be written as follows.   
	\begin{definition}{}{variational_form_sem_disc_vms_rb}
		Let $(\cdot,\cdot)$ be the standard $L^2$ inner product. We state the time-discrete variational problem as follows: find 
		$\vec{v}^{c, k+1}, \in \mathscr{P}\vec{H}^{r}(\Omega)$, and $p^{c, k+1},\;c^{s,c, k+1}(\vec{x}), \phi^{c, k+1} \in \mathscr{P}H^r(\Omega)$
		such that
		\begin{align}
			\begin{split}	
				\text{Momentum Eqns:}& \quad \frac{1}{Sc}\left(w_i,\; \frac{\beta_0v^{c,k+1}_i + \beta_1v^{c,k}_i + \beta_2v^{c,k-1}_i}{\delta t} \right)_{\Omega}\\ 
				& 
				+  \frac{1}{Sc}\left(w_i, \tvj^{k+1}\pd{v_i^{c,k+1}}{x_j} \right)_{\Omega}
				+ \frac{1}{Sc}\left(\pd{w_i}{x_j}, \tvj^{k+1}\tau_m \mathcal{R}_m\left(v_i^{c,k+1}, p^{c,k+1}, c^{s,c, k+1}\right) \right)_{\Omega} \\
				& +  \left(w_i,\pd{p^{c,k+1}}{x_i}\right)_{\Omega} 
				+  \left(\pd{w_i}{x_i},\tau_{sol} \mathcal{R}_{sol}\left(v_i^{c, k+1}\right)\right)_{\Omega}\\
				&+ \left(\pd{w_i}{x_j}, \pd{v_i^{c,k+1}}{x_j}\right)_{\Omega} \\
				&+ \frac{\kappa}{2 \Lambda^2} \left(w_i, \left(\sum c^{s, k+1}z^{s}\right) \pd{\phi^{c, k+1}}{x_i}\right)_{\Omega}= 0, 
				\label{eqn:nav_stokes_var_semi_disc_vms_rb}
			\end{split}\\
			\text{Solenoidality:} & \quad 
			\quad \left(q,\pd{v_i^{c,k+1}}{x_i}\right)_{\Omega} - \left(\pd{q}{x_i}q, \tau_m \mathcal{R}_m\left(v_i^{c,k+1}, p^{c,k+1}, c^{s,c, k+1}\right)\right)_{\Omega} = 0,
			\label{eqn:cont_var_semi_disc_vms_rb} \\
			\begin{split}
				\text{Poisson:} & \quad 
				2\Lambda^2 \left(\pd{q}{x_i}, \pd{\phi^{c,k+1}}{x_i}\right)_\Omega 
				- 2\Lambda^2 \left(q, n_i\pd{\phi^{c, k+1}}{x_i}\right)_{\partial\Omega}\\
				& \quad 
				-  \left(q,  \sum  z^{s} c^{s,c, k+1}\right)_\Omega   
				+ \left(q,  \sum z^{s}\tau_{c^{s}} \mathcal{R}_{c^{s}}\left(v_i^{c, k+1}, \phi^{c,k+1}\right) \right)_\Omega= 0, \label{eqn:poisson_eqn_var_semi_disc_vms_rb}
			\end{split}\\
			\begin{split}
				\text{Nernst-Planck:} 
				& \quad 
				\left(q, \frac{\beta_0 c^{s,c,k+1} + \beta_1 c^{s,c,k} + \beta_2 c^{s,c,k-1}}{\delta t} \right)_\Omega \\ & \quad 
				+ \left(q, \, \tvi^k \pd{c^{s, c, k+1}}{x_i} \right)_\Omega 
				+ \left(\pd{q}{x_i}, \, \tvi^k \tau_{c^{s}} \mathcal{R}_{c^{s}}\left(v_i^{c, k+1}, \phi^{c,k+1}\right) \right)_\Omega \\  & \quad
				+ \left(\pd{q}{x_i}, \, z^{s}c^{s, c, k+1}\;\pd{\phi^{c, k+1}}{x_i} \right)_\Omega   
				- \left(\pd{q}{x_i}, \, z^{s}\tau_{c^{s}} \mathcal{R}_{c^{s}}\left(v_i^{c, k+1}, \phi^{c,k+1}\right)\;\pd{\phi^{c, k+1}}{x_i} \right)_\Omega
				\\  & \quad
				- \left(q, \, n_i z^{s}c^{s, c, k+1}\;\pd{\phi^{c, k+1}}{x_i} \right)_{\partial\Omega}
				\\  & \quad
				+ \left(\pd{q}{x_i}, \, \pd{c^{s, c, k+1}}{x_i} \right)_\Omega 
				- \left(q, \, n_i \pd{c^{s, c, k+1}}{x_i} \right)_{\partial\Omega}
				= 0.
				\label{eqn:np_eqn_var_semi_disc_vms_rb}
			\end{split} 
		\end{align}
		$\forall \vec{w} \in \mathscr{P}\vec{H}^{r}(\Omega)$, $q \in \mathscr{P}H^r(\Omega)$\\
		, and $\vec{v}^{f, k+1} \in (\mathscr{I} - \mathscr{P})\vec{H}^{r}(\Omega)$,
		$c^{s,f, k+1}, p^{f, k+1} \in (\mathscr{I} - \mathscr{P})H^r(\Omega)$\\
		, given $\vec{v}^{k},\;\vec{v}^{k-1} \in \vec{H}^r(\Omega)$, and $\phi^{k},\phi^{k-1}, c^{s,k}, c^{s,k-1} \in H^r(\Omega)$.
	\end{definition}
	
	\section{Solution strategy and design of the numerical framework}
	
	\subsection{Solution strategy}
	\label{sec:sol_strat}
	It is important to note that because we are using a block solution method, the fully discrete Navier-Stokes equations,  \crefrange{eqn:nav_stokes_var_semi_disc_vms_rb}{eqn:cont_var_semi_disc_vms_rb}, and the Poisson Nernst Planck equations, \crefrange{eqn:poisson_eqn_var_semi_disc_vms_rb}{eqn:np_eqn_var_semi_disc_vms_rb}, are solved as two different sub-problems. Note that the NS (\crefrange{eqn:nav_stokes_var_semi_disc_vms_rb}{eqn:cont_var_semi_disc_vms_rb}) is a linear system, whereas, PNP (\crefrange{eqn:poisson_eqn_var_semi_disc_vms_rb}{eqn:np_eqn_var_semi_disc_vms_rb}) is a non-linear system.  We use the Newton method to reduce the non-linear system~(\crefrange{eqn:nav_stokes_var_semi_disc_vms_rb}{eqn:cont_var_semi_disc_vms_rb}) to an iteration of linear problems. Symbolically, we can write this nonlinear algebraic system as
	\begin{equation}
		T_i\left(U^k_1, U^k_2, \dots, U^k_n\right) = 0,
	\end{equation} 
	where $\vec{U}^{k}$ is a finite-dimensional vector that contains all of the degrees of freedom at the discrete-time level $t^k$. In the system of ~\cref{eqn:nav_stokes_var_semi_disc_vms_rb,eqn:cont_var_semi_disc_vms_rb}, $\vec{U}^{k}$ contains ion concentrations and electric potential at time level $k$. Then, the Newton iteration can be posed as, 
	\begin{align}
		A_{ij}^{m,k} \; \delta U_j^{m,k} = -T_i\left(U_1^{m,k}, U_2^{m,k}, 
		\dots, U_n^{m,k} \right), \quad 
		\quad A_{ij}^{m,k} := \pd{}{U_j} F_i\left(U_1^{m,k}, U_2^{m,k}, \dots, U_n^{m,k}\right),
		\label{eqn:newton_linear_system}
	\end{align}
	where  $U_j^{m, k}$ is the vector containing all the degrees of freedom at the $k^{\text{th}}$ time step and at the
	$m^{\text{th}}$ Newton iteration. $\delta U_j^{m, k}$ is the ``variation"  vector that will be used to update the current Newton iteration guess:
	\begin{align}
		U_j^{m+1, k} = U_j^{m, k} + \delta U_j^{m, k}.
	\end{align}
	and $A_{ij}^{m,k}$ is the Jacobian matrix computed analytically by calculating the variations (partial derivatives) of the operators with respect to the degrees of freedom. The iterative procedure begins with an initial guess, which is set to the solution from the previous time step:
	\begin{equation}
		U_i^{0,k} = U_i^{k-1},
	\end{equation}
	and ends once we reach the desired tolerance:
	\begin{equation}
		\|\delta U_j^{m,k}\| \le \text{TOL}.
	\end{equation}
	
\begin{figure}
\centering			
\includegraphics{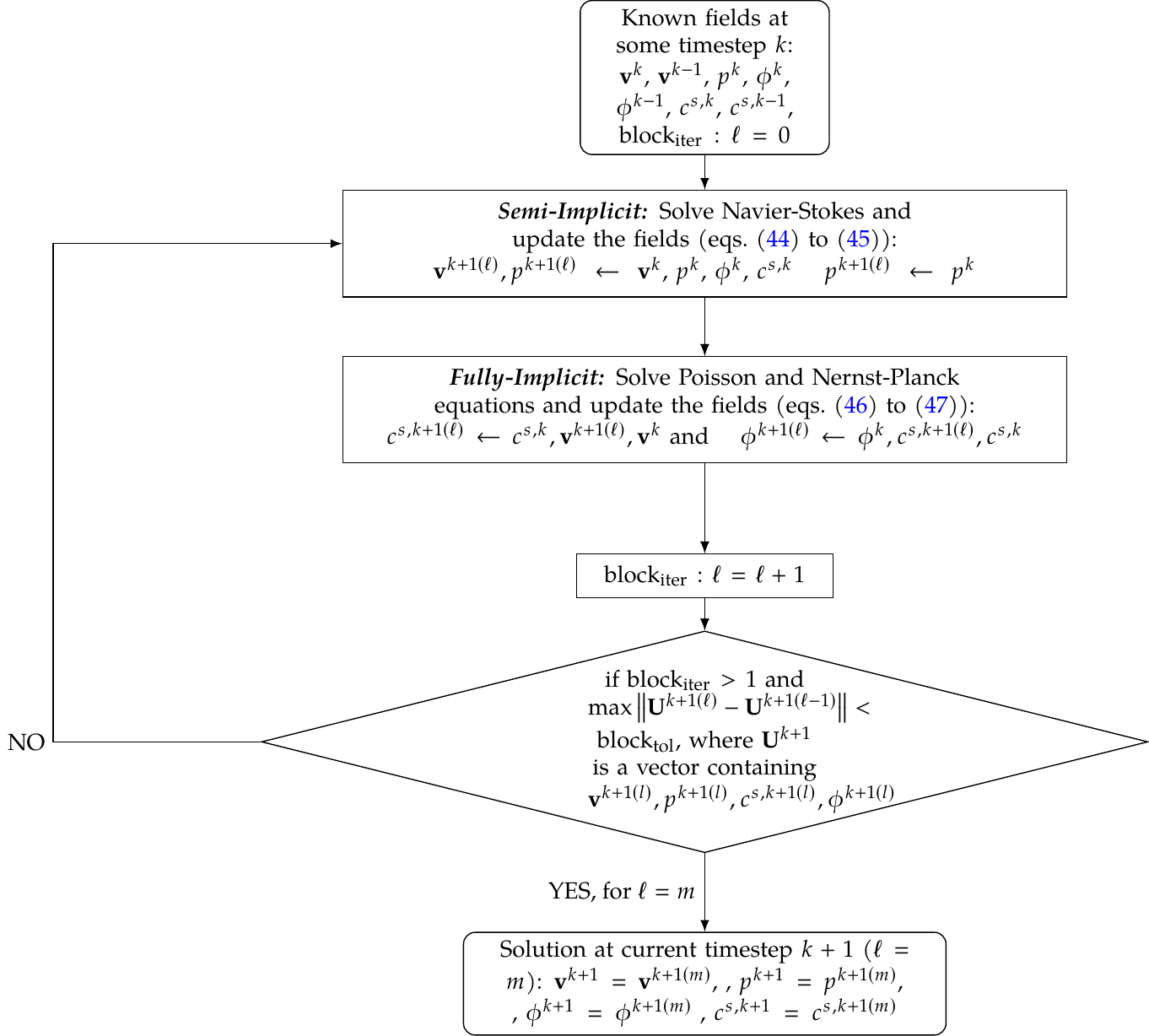}
\caption{Flowchart for the block iteration technique as described
in \cref{sec:sol_strat}
.}
\label{fig:flowchart_block}		
\end{figure}

	The block-iteration strategy for solving the system is given in~\figref{fig:flowchart_block}. The flow chart in~\cref{fig:flowchart_block} shows that the linear semi-implicit Navier Stokes system is solved first from the available data from the previous timestep. Subsequently, the Poisson Nernst-Planck system is solved from the data from the previous timestep and the first block of Navier-Stokes. The block iteration is repeated until a desired convergence is reached. 
	
	Note that the momentum block is linear due to the semi-implicit discretization. The number of iterations of the non-linear Poisson-Nernst Planck (PNP) block depends on the timestep. We generally observe that the non-linear PNP block converges in two to three Newton iterations (each iteration being a linear solve).  
	
	We choose the timestep for accuracy reasons close to the advective CFL of momentum equations. The implicit treatment of higher order operators and the coupling operator allows us to do this (Note that this timestep is still larger compared to the case if the operators were treated explicitly). Then, both PNP and momentum blocks are generally iterated twice. In some cases, like the chaotic electroconvective instability (section 5.2), where the timestep (\num{1e-4}) is chosen to be of the order of smallest element size (0.0009765625). The timestep is small enough that the block iteration converges in just one block of momentum and PNP. The accuracy of quantities of interest (QOI) from our solutions is compared in \figref{fig: AliManiCase}, and it shows very good agreement with high-resolution simulations from~\citet{druzgalski2013direct}. We have added this explanation to the manuscript. 
	
	For chaotic cases where a small timestep is naturally chosen to capture the flow structures, the PNP system relaxes much faster than the momentum system. Therefore, when a small timestep is chosen for momentum (of the order of the smallest mesh size), the convergence of PNP is even faster, and the block iteration converges in one iteration. Numerical experiments reveal that two block iterations produce tight convergence for all examples, with one block iteration often providing good results for most examples.

\subsection{Software design}
The software is designed in a layer-by-layer fashion. The numerical simulation presented in this work is implemented on top of massively parallel framework \textsc{Dendrite-KT}. \textsc{Dendrite-KT} uses ~\dkt{} for adaptive octree mesh generation and ~\petsc{} for linear algebra solver by utilizing \texttt{KSP} (Krylov Subspace) context for solving a linear system of equations and \texttt{SNES} context for Newton iteration. 
~\dkt{} allows allows for targeted, dynamic mesh refinement based on either the location or \textit{aposteriori}  estimates. 

\begin{figure}
\centering
\begin{adjustbox} {width = \textwidth}
\includegraphics{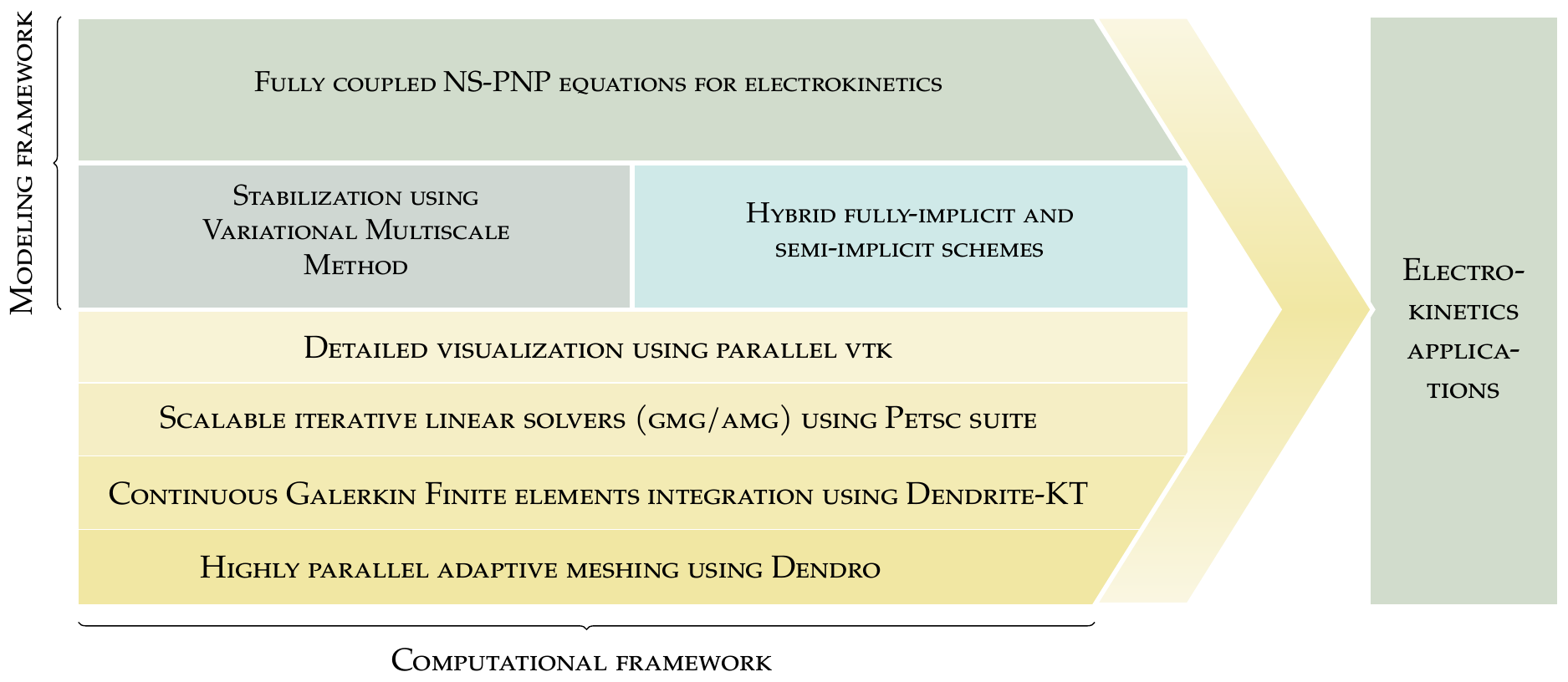}
\end{adjustbox}
\caption{Design schematic of the parallel framework for solving NS-PNP.} 
\label{fig:computational_fram}
\end{figure}

The organization of the overall framework is illustrated in~\figref{fig:computational_fram}. The yellow layers represent the computational framework. The modeling framework - shown as three blocks on top of the computational framework, represents the fully discrete scheme presented  in~\cref{sec:numerical_tecniques}. The non-linear solve for the PNP system (\crefrange{eqn:poisson_eqn_var_semi_disc_vms_rb}{eqn:np_eqn_var_semi_disc_vms_rb}) is solved by coupling the data structures on the octree meshing framework with~\petsc~using the~\texttt{SNES} context. Subsequently, the inner linear solves required for the PNP system and the NS system (\crefrange{eqn:nav_stokes_var_semi_disc_vms_rb}{eqn:cont_var_semi_disc_vms_rb}) are solved using the \texttt{KSP} construct. This allows us to exploit the plethora of preconditioner and linear iterative solvers in~\petsc. The final solution vectors are saved in the parallel~\texttt{VTK} format, which makes visualization using \texttt{ParaView} straightforward.  

\textbf{\dkt~meshing framework}: An ``octant`` represents the fundamental building block of the octree mesh. At the root of the octree, an octant is a cube that encompasses the full domain. In order to construct the octree,  all $p$ compute processes start at the root node. To avoid communication during the refinement stage, we perform redundant computations on all $p$ processes. Starting from the root node, all processes refine (similar to the sequential code) until at least $\mathcal{O}(p)$ octants are produced. Then using a weighted space-filling-curve (SFC) based partitioning, we partition the octants across all processes.  Further, proceeding in a top-down fashion, an octant is refined based on user-specified criteria. The refinement criteria is a user-specified function that takes the coordinates of the octant and the associated values, and returns \texttt{true} or \texttt{false}. Since the refinement happens in an element-local fashion, this step is embarrassingly parallel. In distributed memory, the initial top-down tree construction (which describes the mesh) also enables an efficient partitioning of the domain across an arbitrary number of processes.

The key distinguishing factor for ~\dkt{}, as compared to other octree libraries \citep{BursteddeWilcoxGhattas11,zhang2019amrex,popinet2003gerris}, is the ability to carve out any arbitrarily shaped object based on a simple \Out{} - \In{} classification of a point $\in \mathbb{R}^d$, leading to the generation of \textit{incomplete} octree. An element is said to be \Out{}, if all the nodes of the element falls within the active region of the domain $\Omega$ and \In{} if all the nodes falls within the inactive region. If some of the nodes fall within the active region and some within the inactive, the element is marked \Intercepted{}.  The basic idea of ~\dkt{} is to retain the \Out{} and \Intercepted{} elements and remove the \In{} elements. Such a construction leads to generation of octree that is \textit{incomplete}. An octree is called incomplete if siblings of at least one octant (at any given level) are missing. During the top-down construction, if any given octant falls into the \In{} region of the domain, then the given octant is not further considered a candidate for refinement. This introduces a very efficient way to generate octrees by pruning the tree at a coarser level. Additionally, incomplete octrees provide saving in the total number of degrees of freedom to be solved, in comparison to traditional complete octrees. This translates to additional savings in time-to-solve as well as conditioning of the matrix. We refer interested readers to our previous work ~\citep{saurabh2021scalable,saurabh2022scalable} for further details of incomplete octree construction and scalability performance.

\textbf{{FEM operations:}}
\dkt{} supports both matrix and matrix-free computations. To prevent indirect memory access, \dkt{} does not store any elemental to global map data structure but relies on performing top-down and bottom-up traversals of the mesh tree. The top-down and bottom-up traversals can be thought of as scatter and gather operations, respectively. 

The top-down phase selectively copies nodes from coarser to finer levels until the leaf level is reached. We create buckets for all child subtrees. Looping through the nodes, a node is copied into a bucket if the node is incident on the child subtree corresponding to that bucket. A node that is incident on multiple child subtrees will be duplicated. By recursing on each child subtree and its corresponding bucket of incident nodes, we eventually reach the leaf level. Once the traversal reaches a leaf octant, the elemental nodes have been copied into a contiguous array.
The elemental vector is computed directly without the use of an element-to-node map. The result is stored in a contiguous output buffer the same size as the local elemental input vector.

After all child subtrees have been traversed, the bottom-up phase returns results from a finer to a coarser level.
The parent subtree nodes are once again bucketed to child subtrees,
but instead of the parent values being copied, the values of nodes from each child are accumulated into a parent output array.
For any node that is incident on multiple child subtrees, the values from all node instances are summed to a single value.
The global vector is assembled after the bottom-up phase executes at the root of the octree. A detailed description of this process is provided in our previous work, specifically ~\citet{Ishii2019,saurabh2021scalable,khanwale2023projection}

\textbf{Scalability:}
We deploy the numerical methods on a scalable framework named~\textsc{Dendrite-KT}, which combines the octree-based meshing with continuous Galerkin finite elements. We leverage the highly optimized and scalable implementations of general finite element kernels in~\textsc{Dendrite-KT} built on top of the scalable \dkt{} framework. We have shown massive scalability of~\textsc{Dendrite-KT} for various multiphysics phenomena in~\citep{Ishii2019, saurabh2021industrial, saurabh2022scalable, khanwale2022fully, khanwale2022breakup, khanwale2023projection, saurabh2023scalable}. For example, ~\citet{saurabh2023scalable} recently showed the scalability of the multiphase flow framework called \textsc{Proteus} based on~\textsc{Dendrite-KT} till $\sim 114000$ MPI processes on TACC Frontera.

\section{Results}
\subsection{Method of manufactured solutions for convergence analysis}
\begin{figure}[ht!]
\centering
\centering
\includegraphics{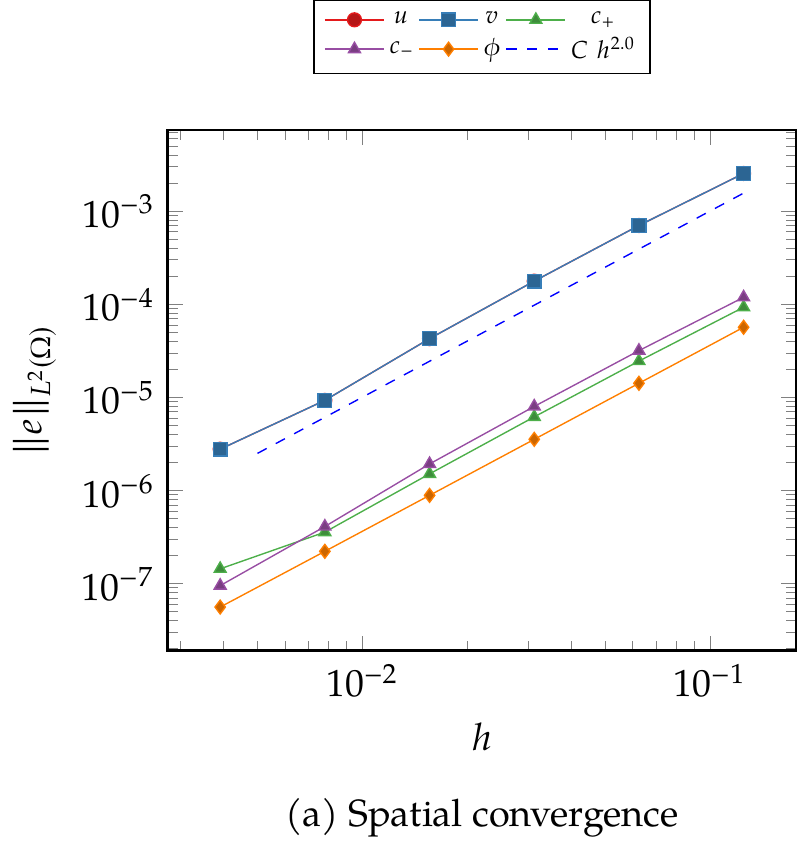}
\includegraphics{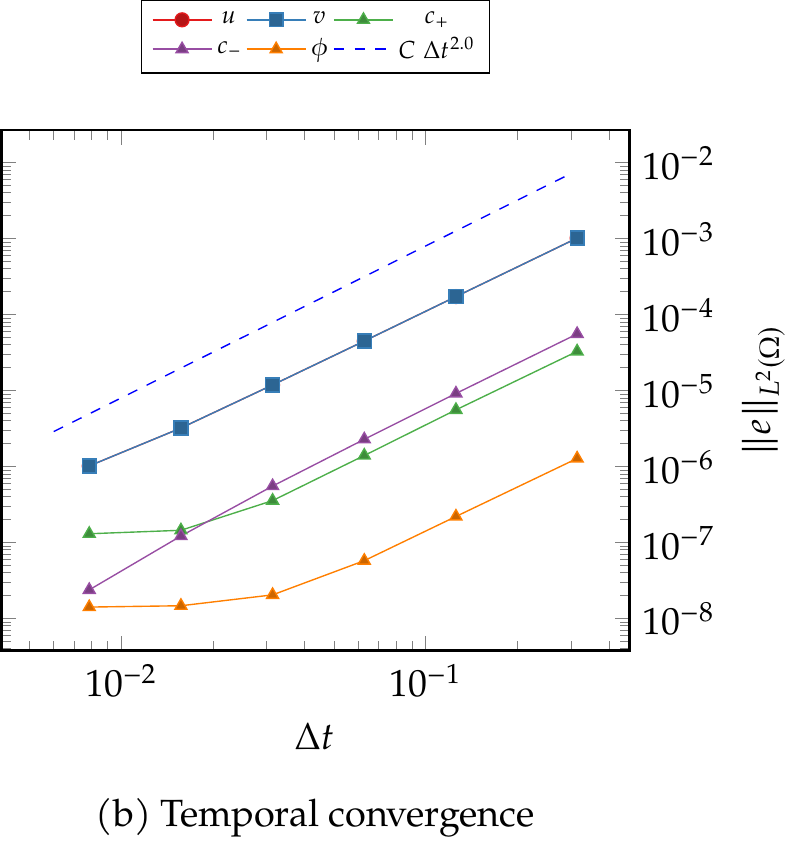}
\caption{Convergence analysis using the method of manufactured solutions. (a) spatial, and (b) temporal convergence of the NS-PNP solver. The dashed line is a reference line of slope 2.}
\label{fig:Convergence}
\end{figure}

We use the method of manufactured solutions to assess the temporal and spatial convergence of our proposed scheme. In this work, we restrict ourselves to a linear basis function ($r = 1$). We construct the appropriate forcing function to satisfy the analytical solution of the form:

\begin{equation}
\begin{split}
u & = \cos(2t) \sin(2 \pi x) \cos(2 \pi y) \\
v & = - \cos(2t) \cos(2 \pi x) \sin(2 \pi y) \\
p & =  \cos(2t) \sin(2 \pi x) \cos(2 \pi y) \\
\phi & = -\cos(2t) \cos(2 \pi x) \sin(2 \pi y) \\
c^{+} & = \cos(2t) \cos(2 \pi x) \sin(2 \pi y) \\
c^{-} & = \cos(2t) \sin(2 \pi x) \cos(2 \pi y)
\end{split}
\label{eqn:manufac_sol}
\end{equation}
We solve the NS-PNP equations on a square domain $[0, 1]^2$, with Dirichlet conditions enforced on all boundaries. We note that the manufactured solution for the fluid velocity in~\cref{eqn:manufac_sol}~is divergence-free. Here, we use a uniform mesh with varying refinement levels for the convergence study. ~\figref{fig:Convergence} shows both the spatial and temporal discretization convergence results. For the spatial convergence analysis, we fixed the timestep $\Delta t$ at 0.0157 and varied the spatial resolution. Next, to study the temporal convergence behavior, we chose the spatial resolution to be level 9 ($\Delta x = \Delta y =  1/2^9$). The $L_2$ error is reported at the final time of $t = \pi$. We observe an expected second-order spatial (~\figref{fig:Convergence}(a)) and temporal convergence (~\figref{fig:Convergence}(a)). At an error of about $\mathcal{O}(10^{-7})$, we see the flattening out of the error curves, especially in temporal convergence studies. This flattening can be attributed to the errors from the spatial resolution limit starting to dominate the $L_2$ error. These results indicate that the implementation exhibits the expected second-order spatial and temporal convergence behavior.

\subsection{Electro-convective instability}\label{sec:alimani}
We next demonstrate electrokinetic transport near an ion-selective surface, which is a good benchmarking problem because 1) the problem can be found in a wide range of applications, 2) the length scale of the EDL near the surface requires an extremely small mesh size, which makes the problem a good test case for the proposed octree-based mesh, and 3) The chaotic nature of the flow triggers structures over a wide range of scales in all field quantities, necessitating numerical stabilization from VMS.
We compare our results against the benchmark results of~\citep{druzgalski2013direct}.  

\begin{figure}
\centering
%
%
%
\includegraphics{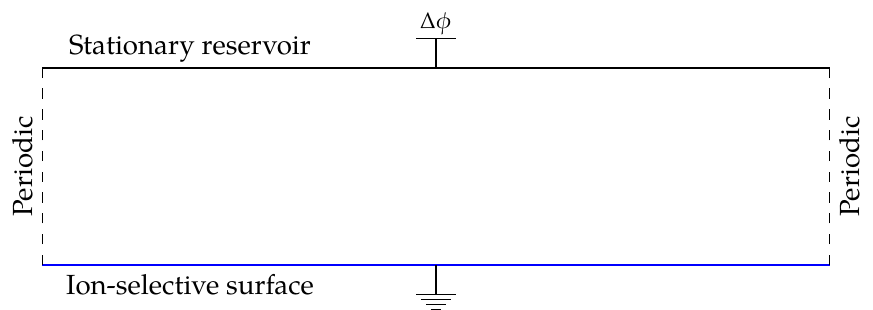}
\caption{Schematic diagram (not on scale) for the electro-convective instability. }
\label{fig:electroconvective-schematic}
\end{figure}

The simulation was performed for the transport of monovalent binary electrolytes in a rectangular domain whose dimension is 8 $\times$ 1. The top boundary is a reservoir at which the potential and the concentration of both species are fixed (Dirichlet). The ion-selective surface is located at the bottom boundary, where the potential is grounded. The cation concentration was $c^{+}=2$, and no anion flux was allowed across the ion-selective surface. A periodic boundary condition is applied to the side walls. \figref{fig:electroconvective-schematic} shows the schematic of the domain, and \tabref{tab:BC-AliMani} lists the boundary conditions for the different variables. We consider two cases, one when a moderate potential difference of $\Delta \phi = 20$ is imposed across the domain and another when a large potential difference of $\Delta \phi = 120$ is imposed.  We anticipate significant electro-convective instabilities to occur in the second case. 


\begin{table}[htb!]
\centering
\renewcommand{\arraystretch}{1.2}
\begin{tabular}{|ccccc|}
\hline
\multicolumn{2}{|c|}{} & \multicolumn{1}{c|}{\textbf{Variable}} &\multicolumn{1}{c|}{\textbf{Type}} & \textbf{Value} \\ \hline
\hline
\multicolumn{1}{|c|}{\multirow{5}{*}{NS}} & \multicolumn{1}{c|}{\multirow{2}{*}{Top Wall}} & \multicolumn{1}{c|}{u} & \multicolumn{1}{c|}{Dirichlet} & 0 \\ \cline{3-5} 
\multicolumn{1}{|c|}{} & \multicolumn{1}{c|}{} & \multicolumn{1}{c|}{v} & \multicolumn{1}{c|}{Dirichlet} & 0 \\ \cline{2-5} 
\multicolumn{1}{|c|}{} & \multicolumn{1}{c|}{\multirow{2}{*}{Bottom Wall}} & \multicolumn{1}{c|}{u} & \multicolumn{1}{c|}{Dirichlet} & 0 \\ \cline{3-5} 
\multicolumn{1}{|c|}{} & \multicolumn{1}{c|}{} & \multicolumn{1}{c|}{v} & \multicolumn{1}{c|}{Dirichlet} & 0 \\ \cline{2-5} 
\multicolumn{1}{|c|}{} & \multicolumn{1}{l|}{Side walls} & \multicolumn{3}{c|}{Periodic} \\ \hline
\hline
\multicolumn{1}{|c|}{\multirow{7}{*}{PNP}} & \multicolumn{1}{c|}{\multirow{3}{*}{Top Wall}} & \multicolumn{1}{c|}{$\phi$} & \multicolumn{1}{c|}{Dirichlet} & 20  or 120\\ \cline{3-5} 
\multicolumn{1}{|c|}{} & \multicolumn{1}{c|}{} & \multicolumn{1}{c|}{$c_{+}$} & \multicolumn{1}{c|}{Dirichlet} & 1 \\ \cline{3-5} 
\multicolumn{1}{|c|}{} & \multicolumn{1}{c|}{} & \multicolumn{1}{c|}{$c_{-}$} & \multicolumn{1}{c|}{Dirichlet} & 1 \\ \cline{2-5} 
\multicolumn{1}{|c|}{} & \multicolumn{1}{c|}{\multirow{3}{*}{Bottom Wall}} & \multicolumn{1}{c|}{$\phi$} & \multicolumn{1}{c|}{Dirichlet} & 0 \\ \cline{3-5} 
\multicolumn{1}{|c|}{} & \multicolumn{1}{c|}{} & \multicolumn{1}{c|}{$c_{+}$} & \multicolumn{1}{c|}{Dirichlet} & 2 \\ \cline{3-5} 
\multicolumn{1}{|c|}{} & \multicolumn{1}{c|}{} & \multicolumn{1}{c|}{$c_{-}$} & \multicolumn{1}{c|}{Zero flux} & 0 \\ \cline{2-5} 
\multicolumn{1}{|c|}{} & \multicolumn{1}{l|}{Side walls} & \multicolumn{3}{c|}{Periodic} \\ \hline
\end{tabular}
\caption{Boundary condition for the electro-convective instabilty}
\label{tab:BC-AliMani}
\end{table}

\begin{table}[h]
\centering
\begin{tabular}{|c|c|c|l|l|c|}
\hline
$\mathbf{\Delta \phi}$  & \textbf{Mesh size}                                         & \textbf{\# of MPI tasks}                                        & \textbf{\# of timestep} & \textbf{Final time}    & \textbf{Wall time} \\ \hline
$20$ & \multirow{2}{*}{\begin{tabular}[c]{@{}c@{}} $\sim$ 1M\\ (8.3 M equivalent)\end{tabular}} & \multirow{2}{*}{\begin{tabular}[c]{@{}c@{}}224 CPU processes\\ (4 nodes on TACC Frontera)\end{tabular}}  & 1000 & 10 & 5 hours    \\ \cline{1-1} \cline{4-6} 
$120$ &                                                                                  &                                                                                               & 12000    & 1.2 & 4 days     \\ \hline
\end{tabular}
\caption{Mesh details and time taken for the electro-convective instability.}
\label{table: mesh_AliMani}
\end{table}

Similar to the previous example, finer refinement was used near the ion-selective surface to resolve steep gradients in the boundary layer. The bottom part of the domain, extending from y = 0 to y = 0.1, consists of elements refined to level 13 ($\Delta x = \Delta y =  8/2^{13}$), while the rest of the region consists of elements refined to level 10 ($\Delta x = \Delta y =  8/2^{10}$). The setup resulted in a domain with approximately 1M cells. In contrast, an equivalent uniform mesh at the finest resolution would have $\sim~10$M cells. The initial condition with fully developed EDL was generated by solving only the PNP equations, then randomly perturbing the resultant concentrations locally by 1\% to initiate the instability before solving the fully coupled NS-PNP equations.
The timestep was $10^{-2}$ for the case with potential difference, $\Delta \phi = 20$, and a time step of $10^{-4}$ for the case of potential difference $\Delta \phi = 120$. We ran the case of $\Delta \phi = 20$ until about 10 seconds when it reached steady-state (approximately 1000 timesteps). For the case of $\Delta \phi = 120$ we ran the simulation for 1.2 seconds, which encompasses the time to reach statistical stationarity and the time used to collect data (approximately 12000 timesteps). \tabref{table: mesh_AliMani} briefs the mesh details along with the number of MPI tasks and the wall time elapsed for both the cases.

We first compare the temporal evolution of the charge density ($c^{+} - c^{-}$) and velocity fields under the imposition of the potential drops. This is shown in \figref{fig:plot20V} for the $\Delta \phi = 20$ case and \figref{fig:plot120V} for the $\Delta \phi = 120$ case. Once a potential drop is applied, the ion concentration profile near the ion-selective surface forms an electric double layer (EDL) and diffusion layer (DL). Depending on the level of the potential drop, the extended space charge layer (ESC) may form between EDL and DL~\citep{druzgalski2013direct}. When $\Delta \phi=20$, the ESC is not formed, and the magnitude of convection is insufficient to trigger instabilities. A series of counter-rotating vortices are observed in \figref{fig:plot20V}.
~On the other hand, an ESC was observed in the initial stage of $\Delta \phi=120$. The interplay between ESC and convection triggered electrohydrodynamic instability, as reported in~\citet{druzgalski2013direct}. As the simulation proceeds, the development of electrohydrodynamic instability becomes even more evident. We note that the average velocity at $\Delta \phi = 120$ was greater than that of $\Delta \phi = 20$ case by two orders of magnitude.

\begin{figure}
\begin{subfigure}{0.45\textwidth}
\includegraphics{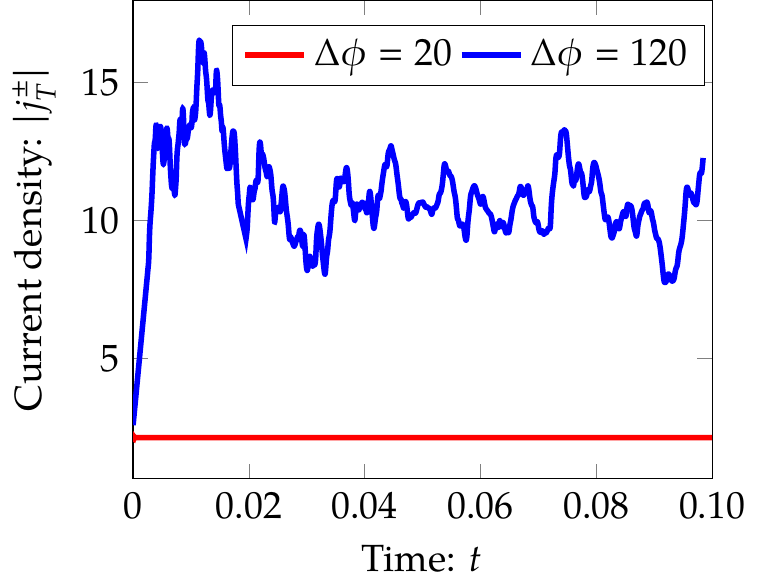}
\caption{Variation of current with time}
\label{fig:current}
\end{subfigure}
\begin{subfigure}{0.45\textwidth}
\includegraphics{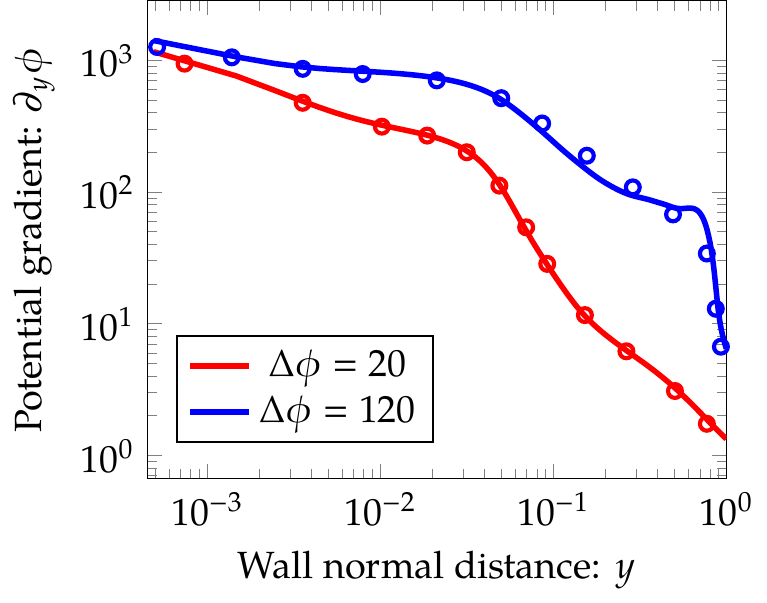}
\caption{Variation of potential gradient with $y$}
\label{fig:ptVar}
\end{subfigure}
\par\bigskip
\begin{subfigure}{0.45\textwidth}
\includegraphics{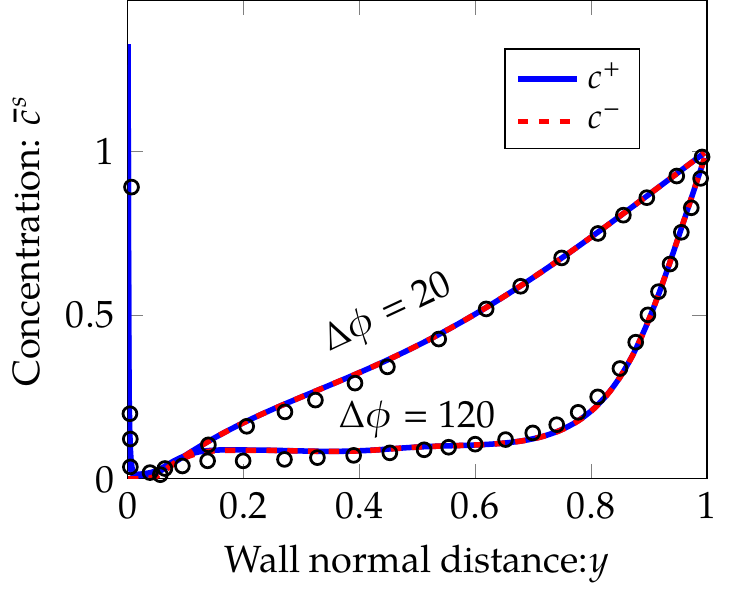}
\caption{Variation of concentration with $y$}
\label{fig:cVar}
\end{subfigure}
\begin{subfigure}{0.45\textwidth}
\includegraphics{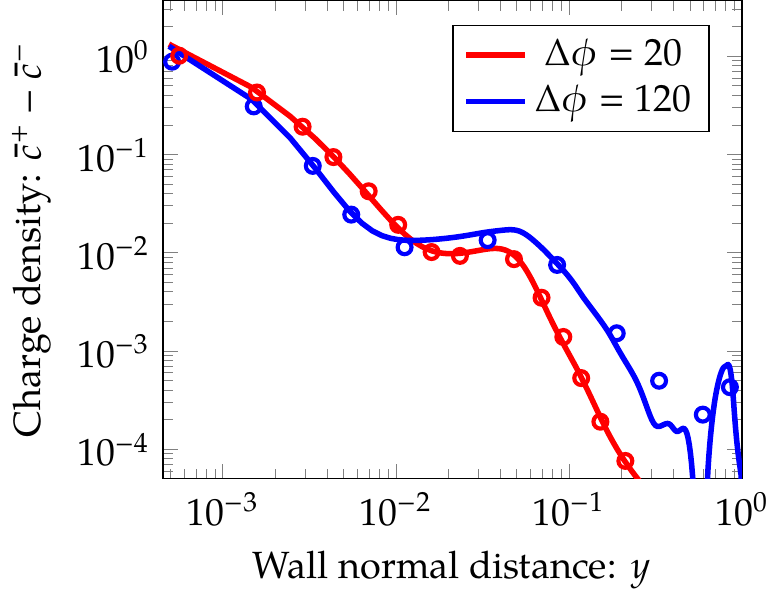}
\caption{Variation of charge density  with $y$}
\label{fig:cdVar}
\end{subfigure}
\caption{ \textit{Electro-osmotic instability:} Variation of current with time (in s) for different $\Delta \phi$ is shown in ~\figref{fig:current}.  Time and X averaged variation of potential gradient ($\partial_y\phi$) (\figref{fig:ptVar}), concentration (\figref{fig:cVar}) and charge density  (\figref{fig:cdVar}) with wall normal distance $y$ for $\Delta \phi = $ 20  and 120. The dotted markers in the plot corresponds to the previous DNS results by ~\citet{druzgalski2013direct}}
\label{fig: AliManiCase}
\end{figure}

We next extract quantitative properties for comparison. The current measured on the boundary is an often used property, as it is the easiest property to measure experimentally.~\figref{fig: AliManiCase} (a) shows the temporal variation of the magnitude of the net species flux, $|j_{T}^{\pm}|$, through the top boundary, calculated as 
\begin{equation}
j_{T}^{\pm}= - \frac{1}{L_{T}}\sum_{l=1}^{N_{T}}\left[\left(q, \, \hat{n}_T \left(\pd{c^{+, t}}{x_i} + z^{+}c^{+, t}\;\pd{\phi^{t}}{x_i}\right) \right)_{\Gamma_{l}}-\left(q, \, \hat{n}_T \left(\pd{c^{-, t}}{x_i} + z^{-}c^{-, t}\;\pd{\phi^{t}}{x_i}\right) \right)_{\Gamma_{l}}\right],
\end{equation}
where ${L_{T}}$, $N_{T}$, $\hat{n}_T$, and $\Gamma=\partial{\Omega_{\Gamma}}$ are the length of the top boundary, the number of the top boundary elements, the normal vector to the top boundary, and the top boundary, respectively. 
We see that, in the absence of electro-convective instabilities, the current (for the $\Delta \phi = 20$ case) was uniform over time, as expected in the absence of ESC and not enough convection to perturb established EDL. In contrast, significant current fluctuations were observed for the $\Delta \phi = 120$ case due to the instability.

We finally plot $x-$averaged field values for the potential (\figref{fig: AliManiCase} (b)), concentration (\figref{fig: AliManiCase} (c)) and charge density (\figref{fig: AliManiCase} (d)). Note the  logarithmic axis in these figures to ensure that the variation close to the ion-selective surface is emphasized. These figures also allow detailed comparison with previous DNS simulations of the same case by ~\citet{druzgalski2013direct}, which shows excellent agreement. 
Note that our numerical scheme allows three orders of magnitude larger time steps as compared to the benchmarking problem by ~\citet{druzgalski2013direct}, with identical results.



\newcommand \figwidthCase{20}
\begin{figure}[htb!]
\begin{tabular}{|l|l|l|l|}
\hline
& 
\parbox[c]{\figwidthCase em}{
\includegraphics[width=\linewidth,trim=100 1200 100 100,clip]{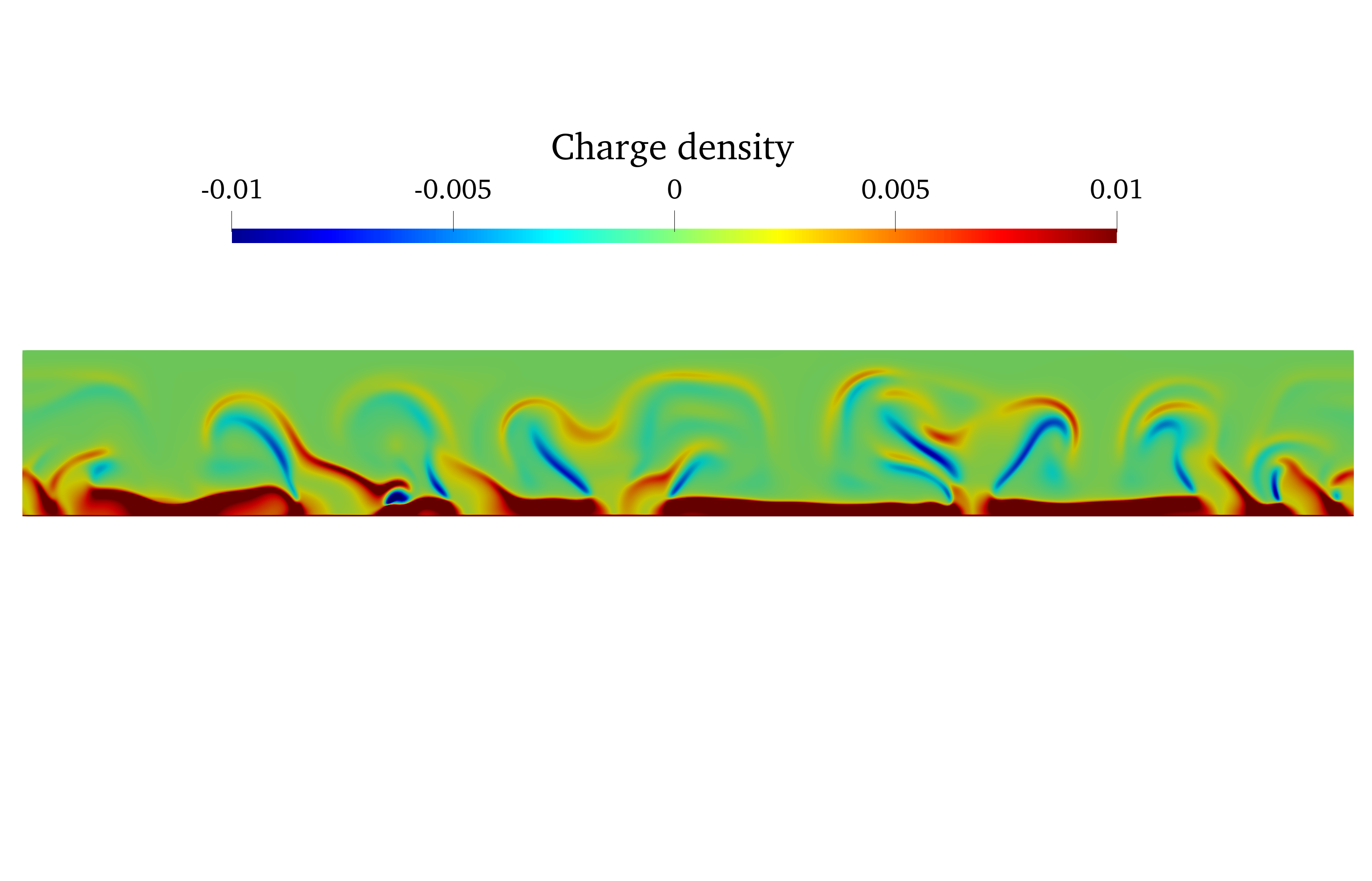}}
& 
\parbox[c]{\figwidthCase em}{
\includegraphics[width=\linewidth,trim=100 1300 110 50,clip]{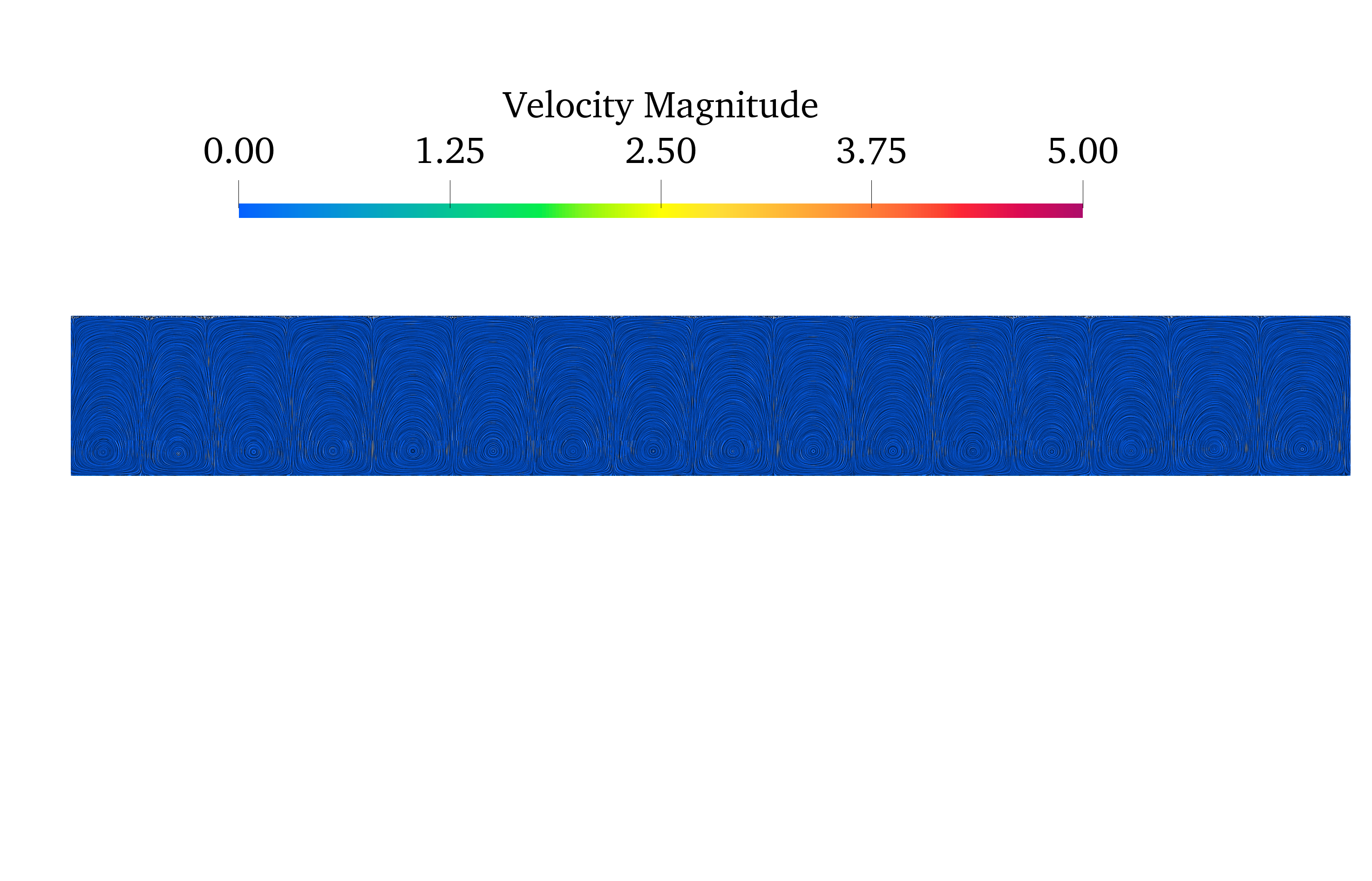}}
\\ \hline
\rotatebox[origin=c]{90}{t = 0.01} 
& 
\parbox[c]{\figwidthCase em}{
\includegraphics[width=\linewidth,trim=100 550 100 500,clip]{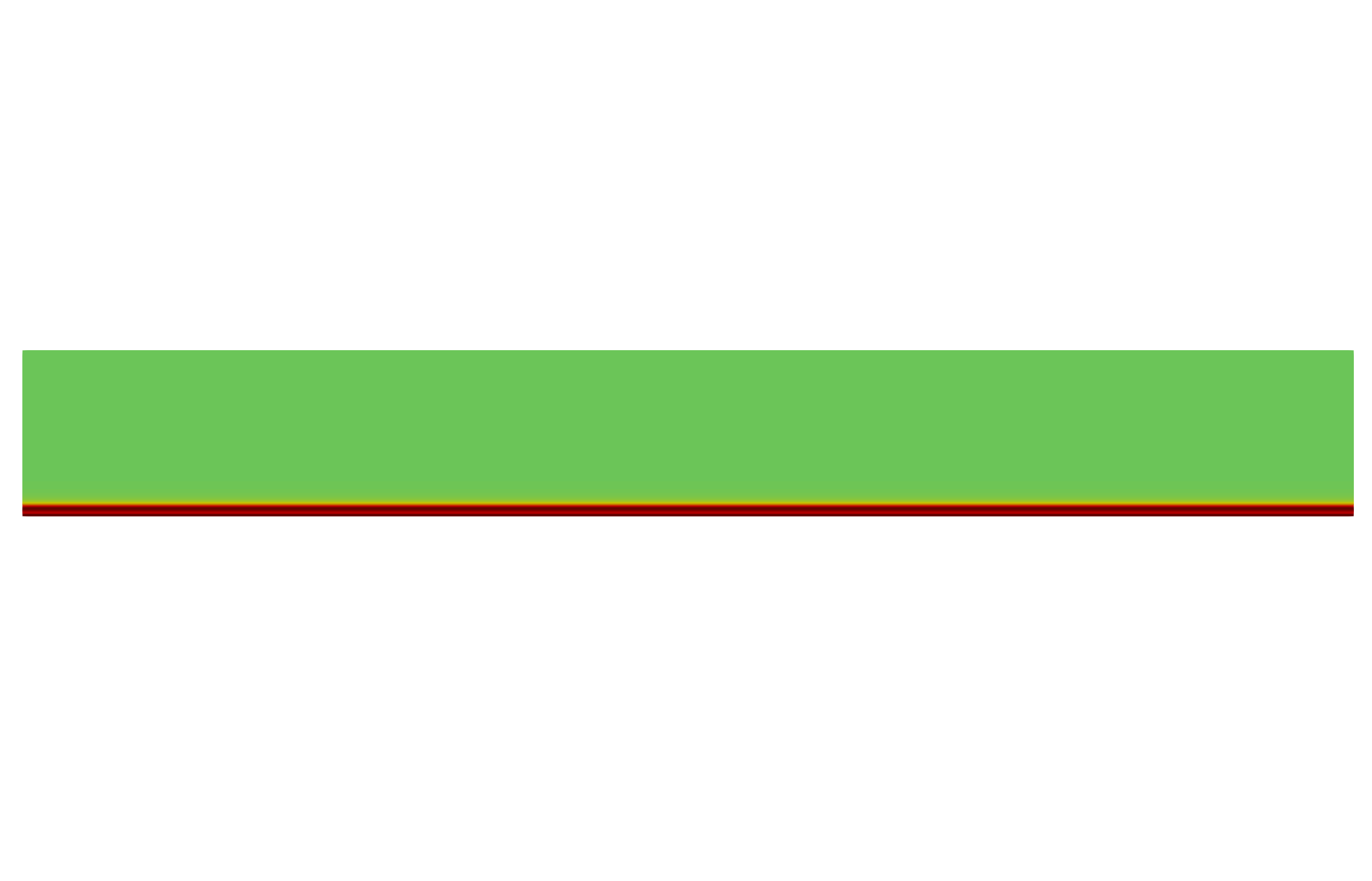}}
&
\parbox[c]{\figwidthCase em}{
\includegraphics[width=\linewidth,trim=50 300 50 200,clip]{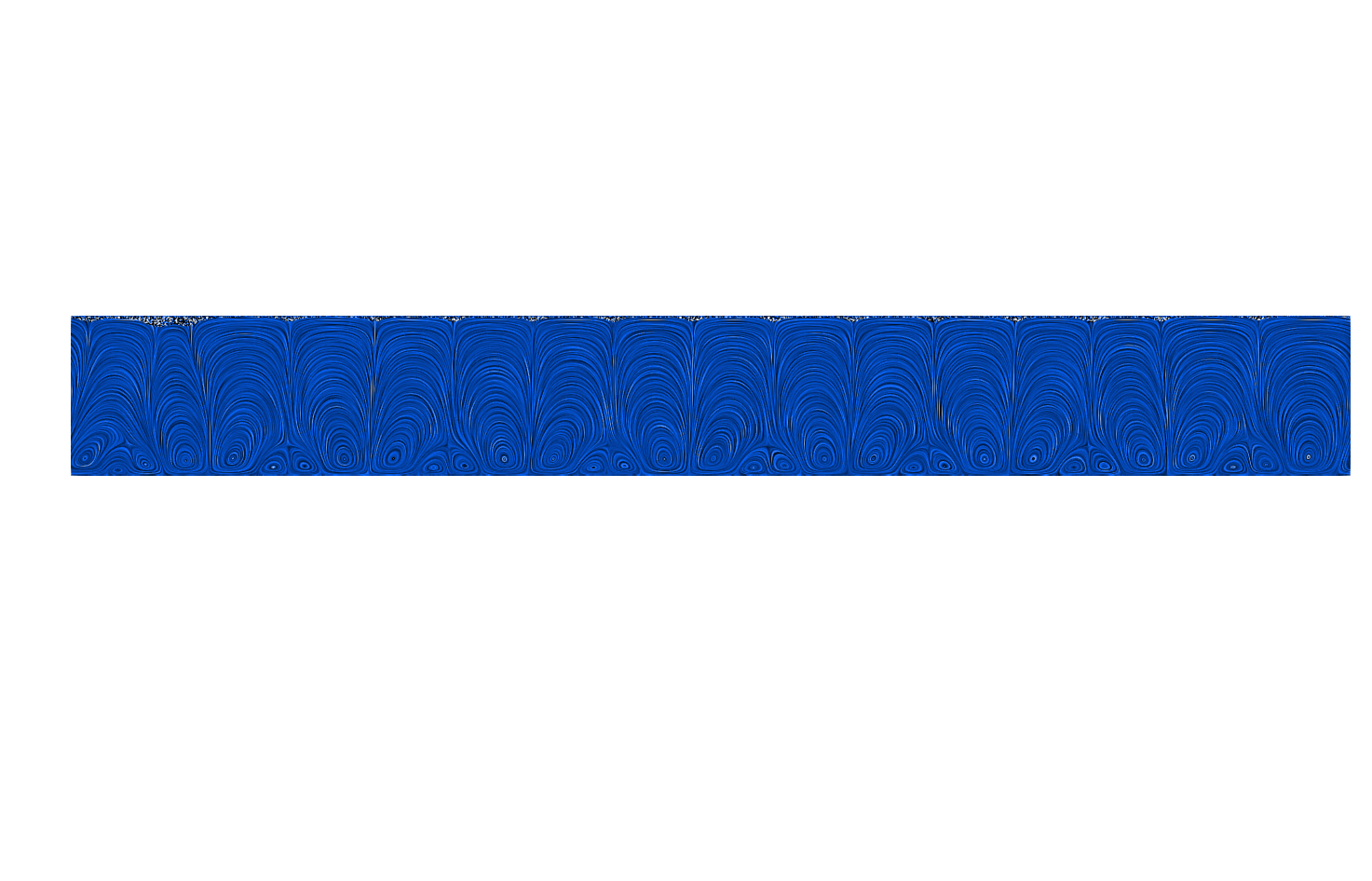}}    

\\ \hline
\rotatebox[origin=c]{90}{t = 0.1} 
&
\parbox[c]{\figwidthCase em}{
\includegraphics[width=\linewidth,trim=100 550 100 500,clip]{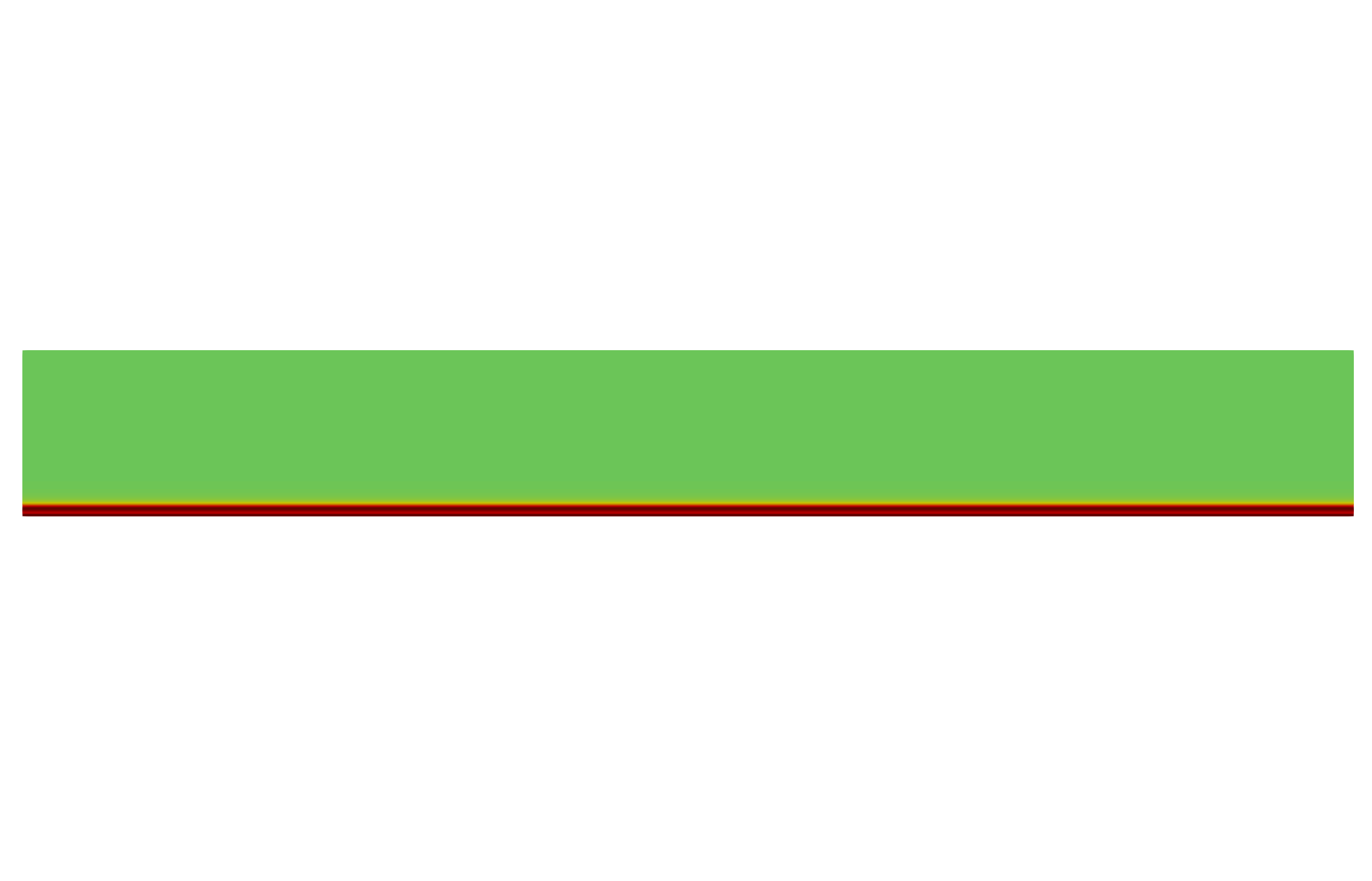}} 
&  
\parbox[c]{\figwidthCase em}{
\includegraphics[width=\linewidth,trim=50 300 50 200,clip]{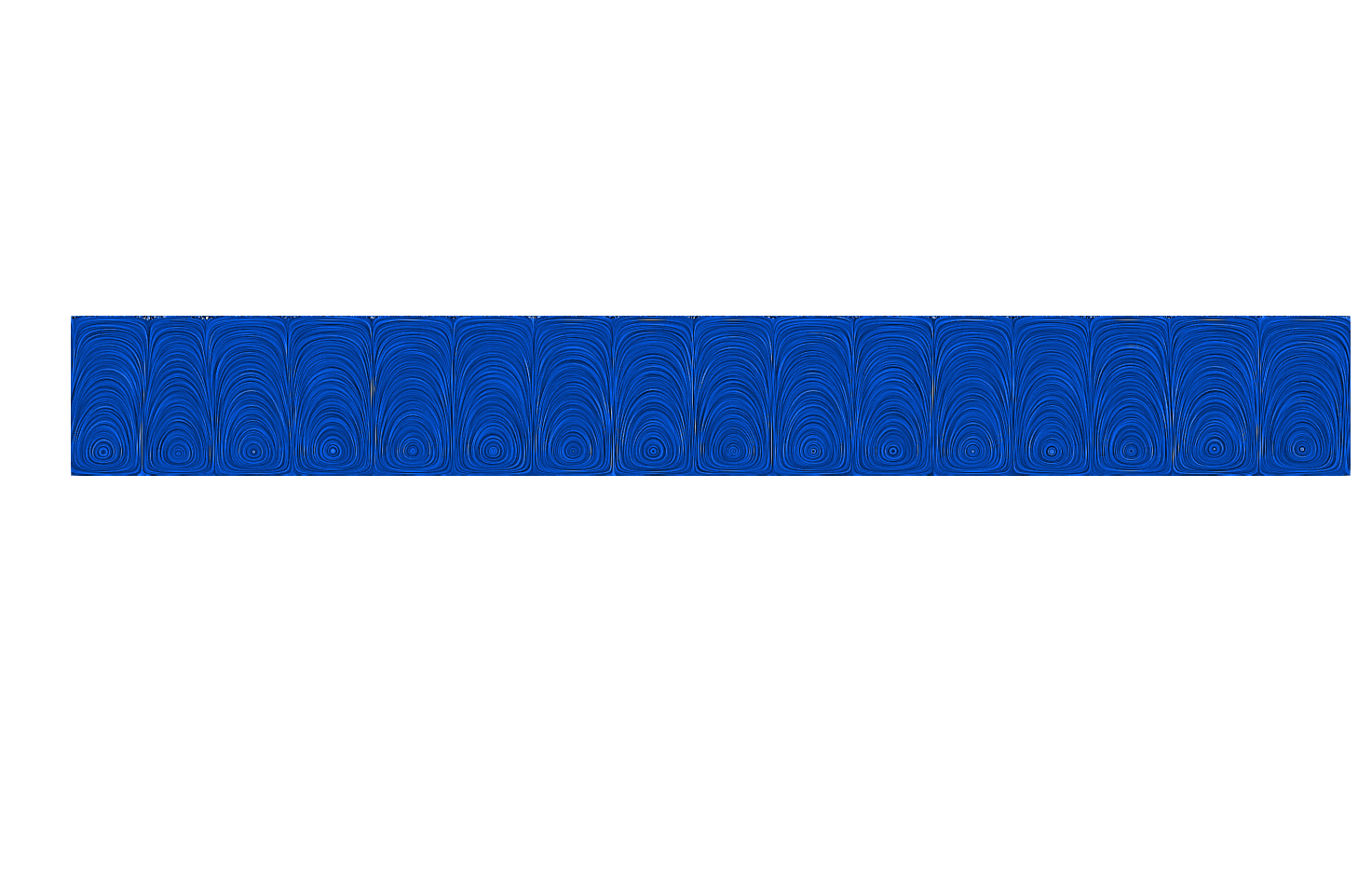}}

\\ \hline
\rotatebox[origin=c]{90}{t = 1.0} 
&
\parbox[c]{\figwidthCase em}{
\includegraphics[width=\linewidth,trim=100 550 100 500,clip]{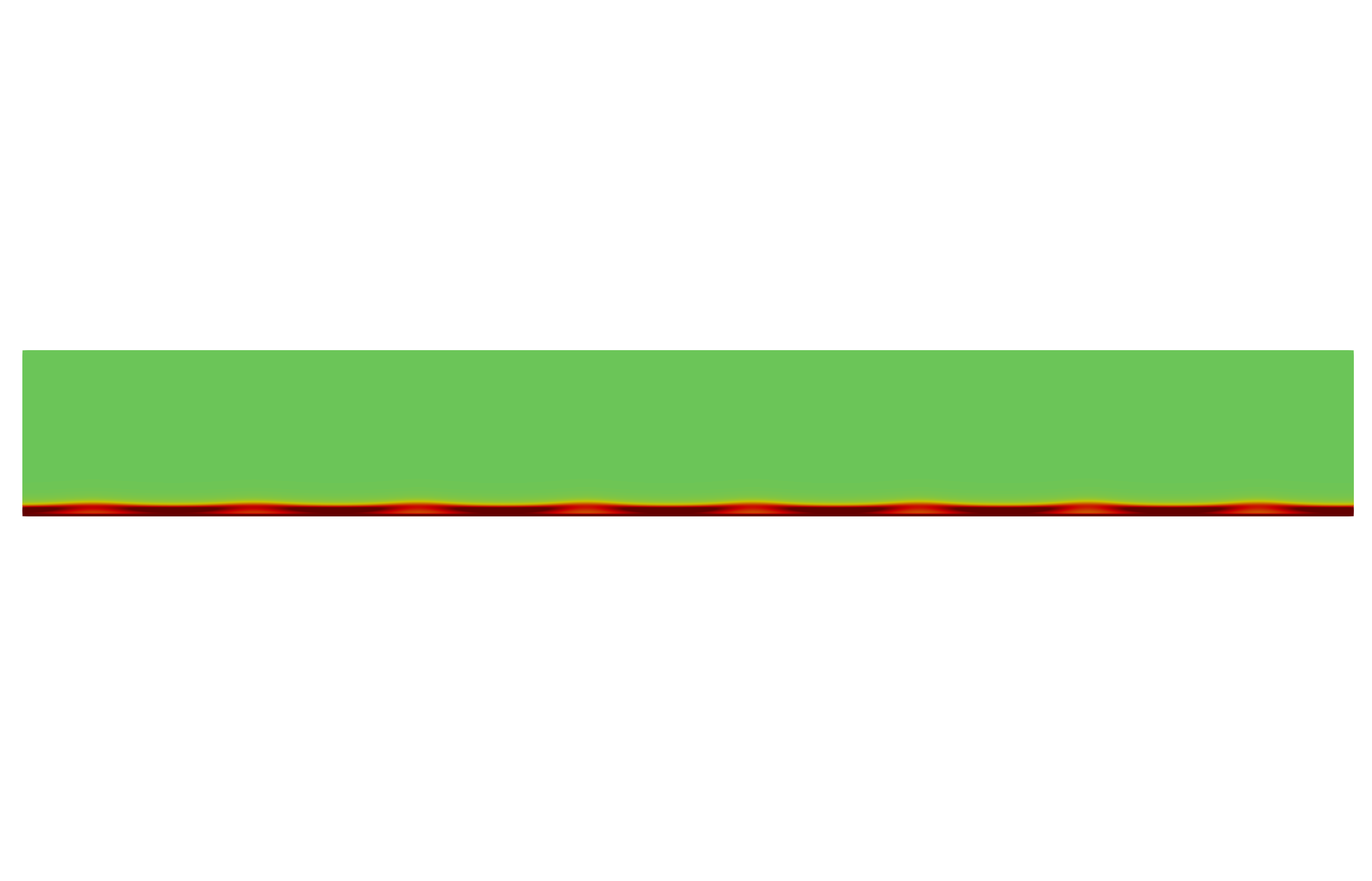}}  
&
\parbox[c]{\figwidthCase em}{
\includegraphics[width=\linewidth,trim=50 300 50 200,clip]{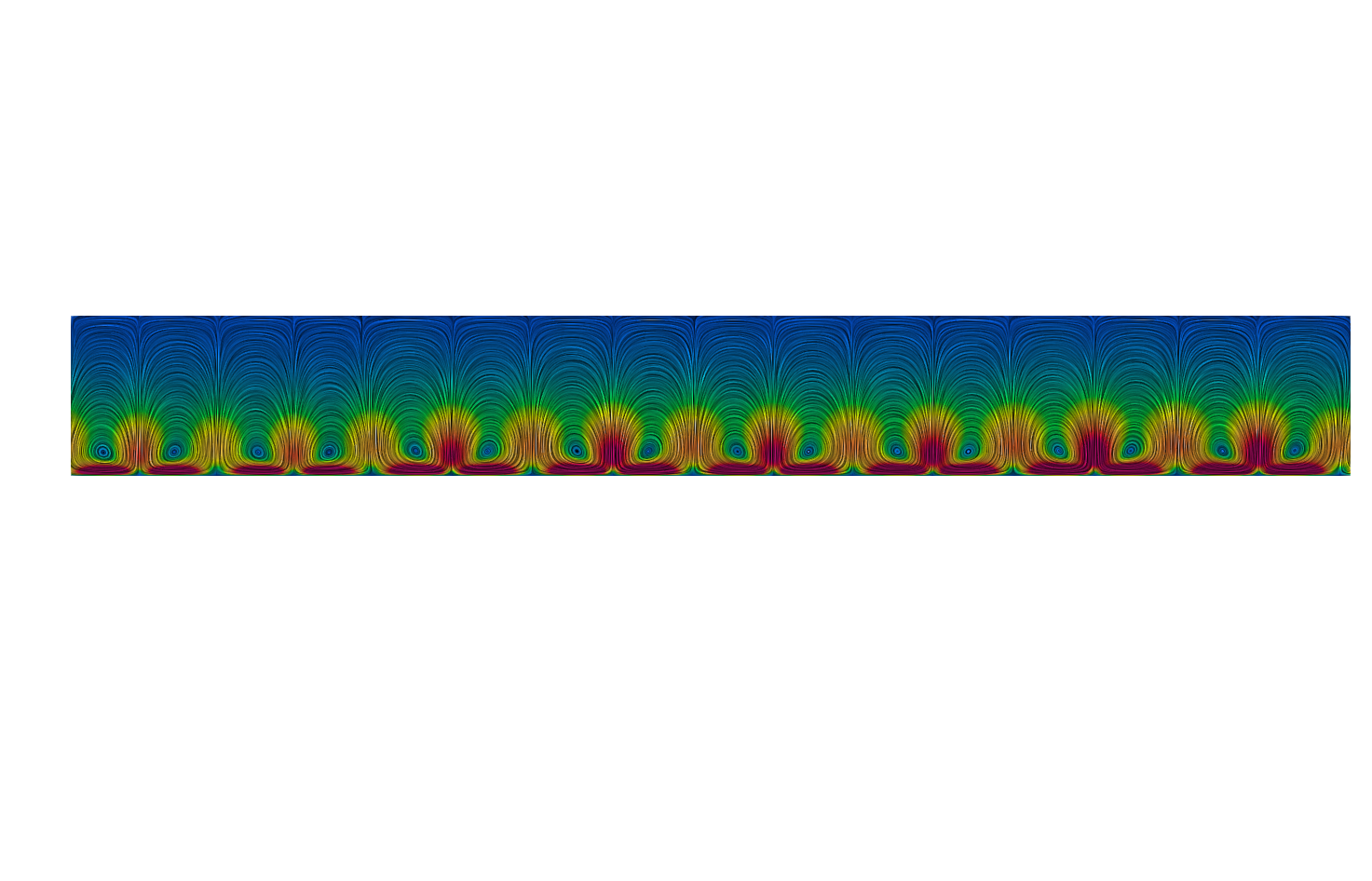}}

\\ \hline
\rotatebox[origin=c]{90}{t = 5.0} & 
\parbox[c]{\figwidthCase em}{
\includegraphics[width=\linewidth,trim=100 550 100 500,clip]{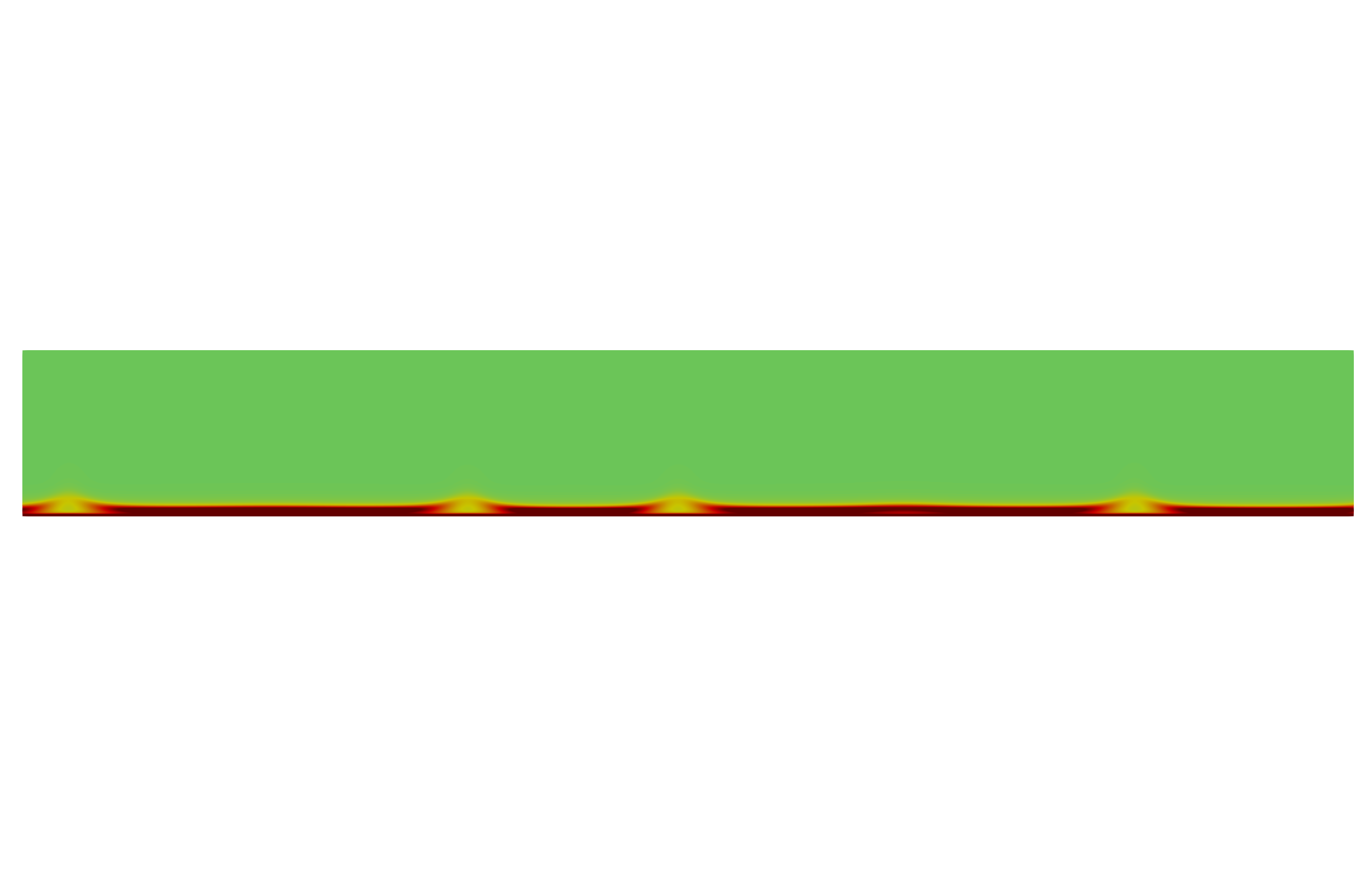}} 
&
\parbox[c]{\figwidthCase em}{
\includegraphics[width=\linewidth,trim=50 300 50 200,clip]{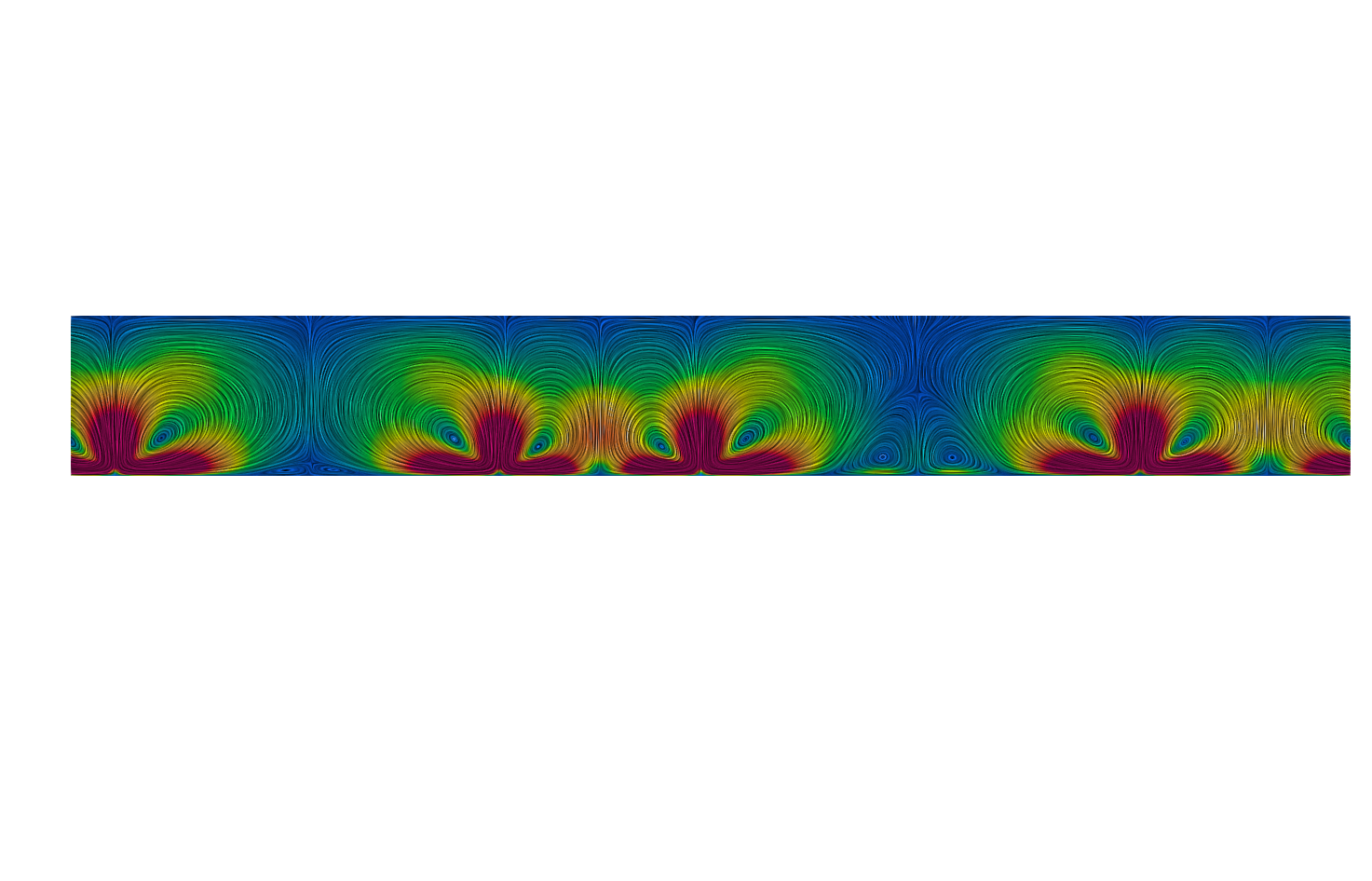}}  
\\ \hline
\end{tabular}
\caption{Contour plots of charge density and velocity magnitudes  at different times for $\Delta \phi = $ 20. }
\label{fig:plot20V}
\end{figure}


\begin{figure}[htb!]
\begin{tabular}{|l|l|l|l|}
\hline
& 
\parbox[c]{\figwidthCase em}{
\includegraphics[width=\linewidth,trim=100 1200 100 100,clip]{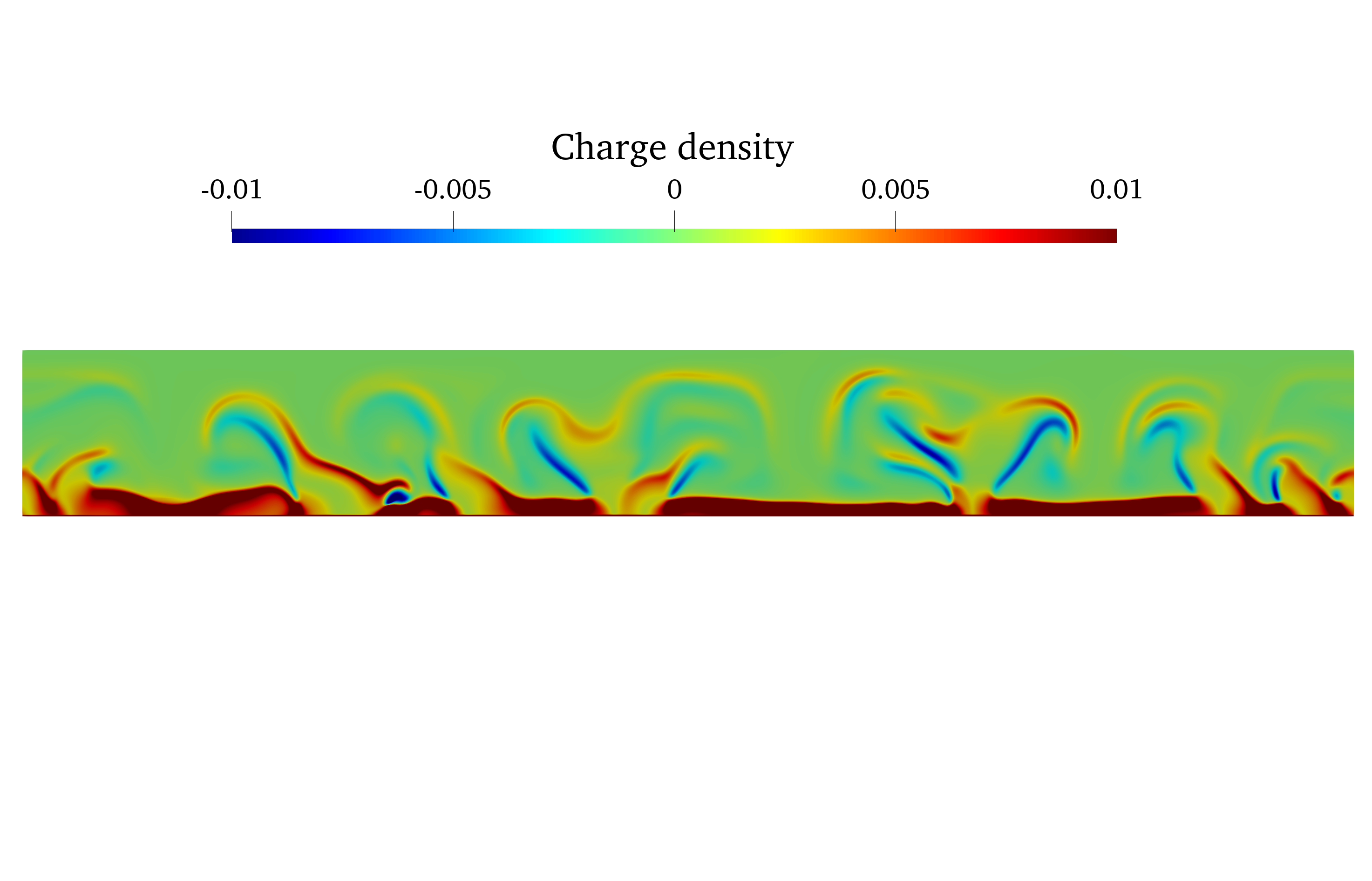}}
& 
\parbox[c]{\figwidthCase em}{
\includegraphics[width=\linewidth,trim=100 1200 100 100,clip]{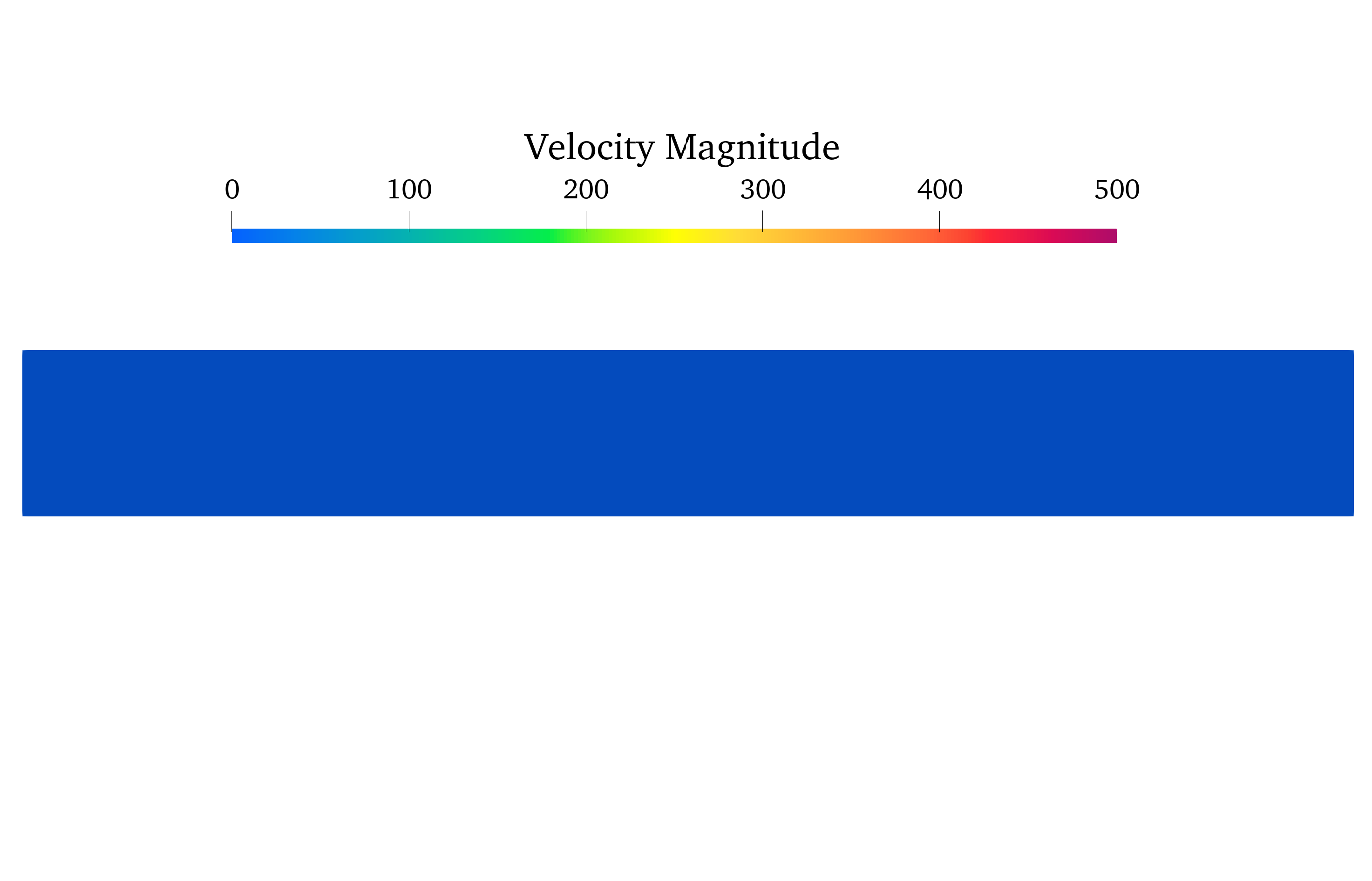}}
\\ \hline
\rotatebox[origin=c]{90}{t = 0.001} 
& \parbox[c]{\figwidthCase em}{
\includegraphics[width=\linewidth,trim=100 550 100 500,clip]{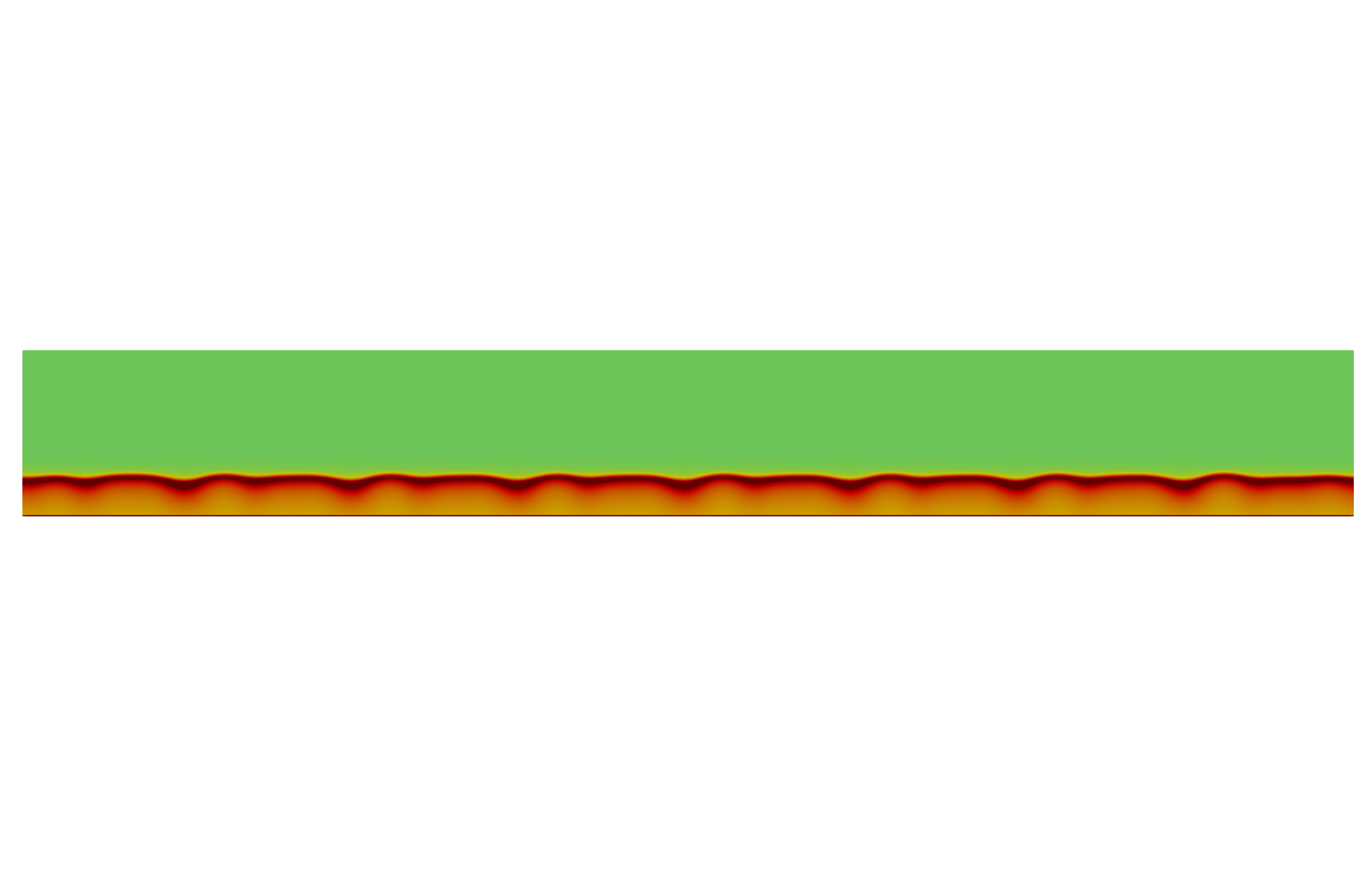}}
& \parbox[c]{\figwidthCase em}{
\includegraphics[width=\linewidth,trim=50 275 50 250,clip]{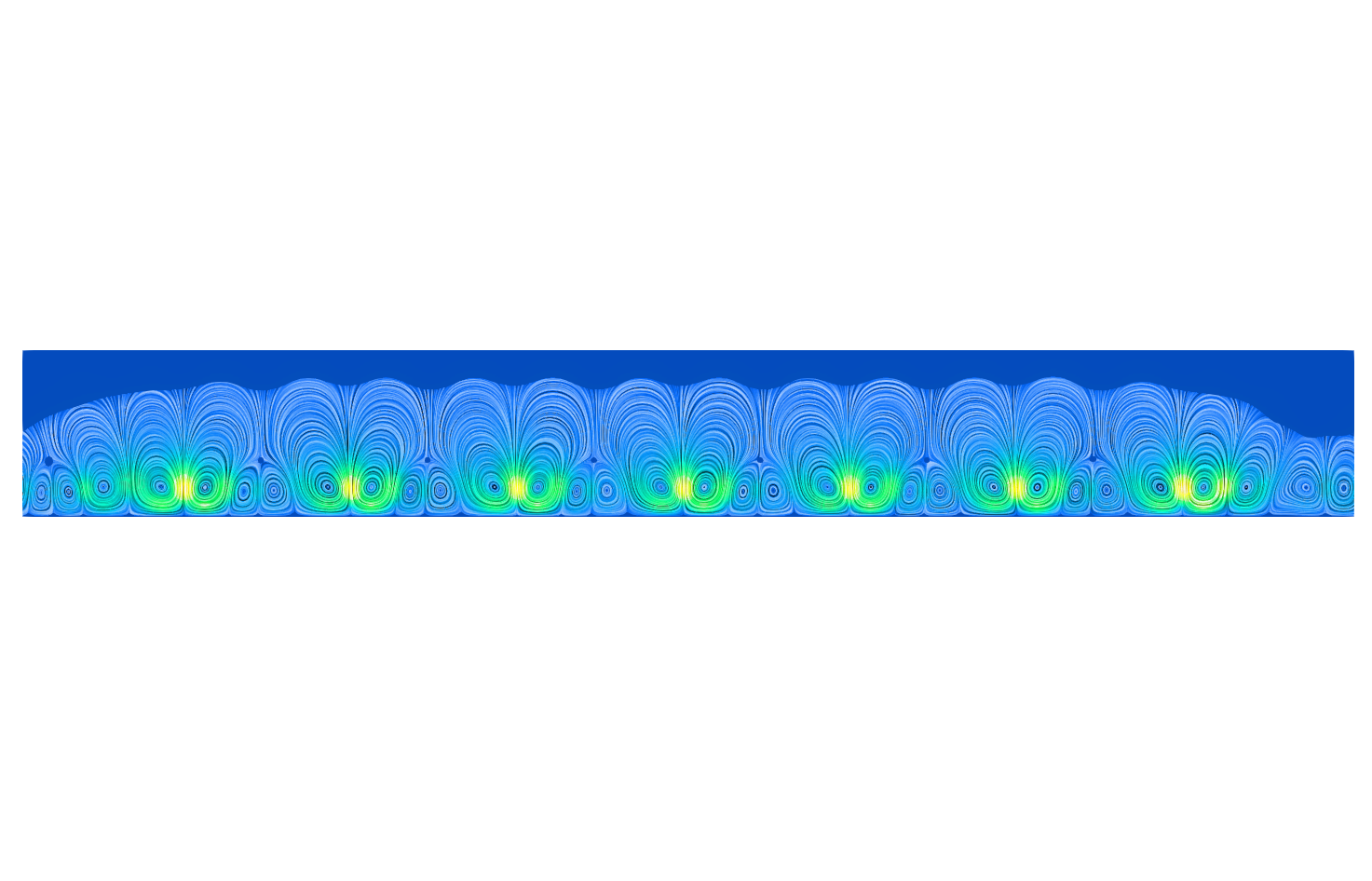}}   

\\ \hline
\rotatebox[origin=c]{90}{t = 0.005} & \parbox[c]{\figwidthCase em}{
\includegraphics[width=\linewidth,trim=100 550 100 500,clip]{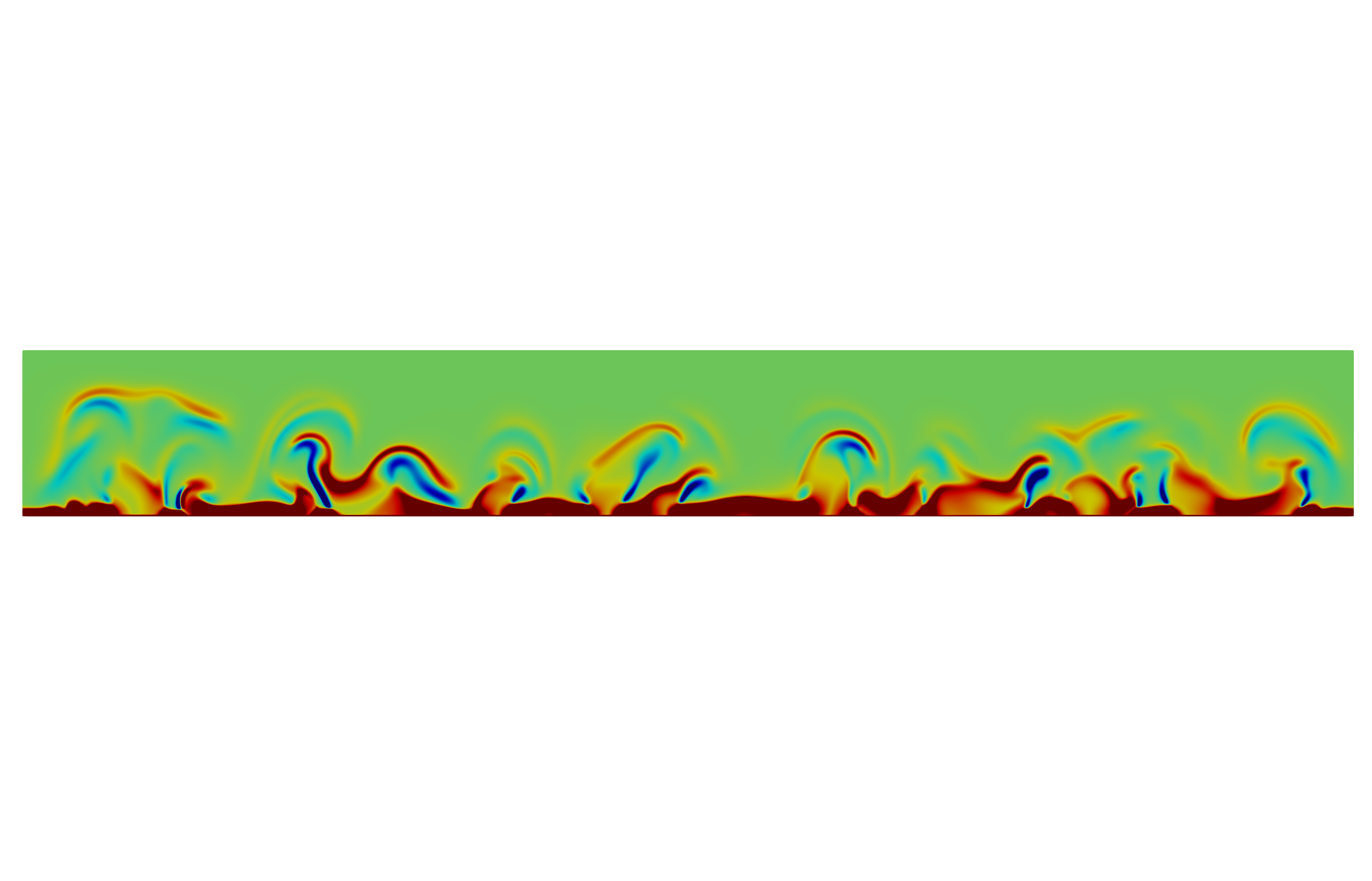}} 
& 
\parbox[c]{\figwidthCase em}{
\includegraphics[width=\linewidth,trim=50 275 50 250,clip]{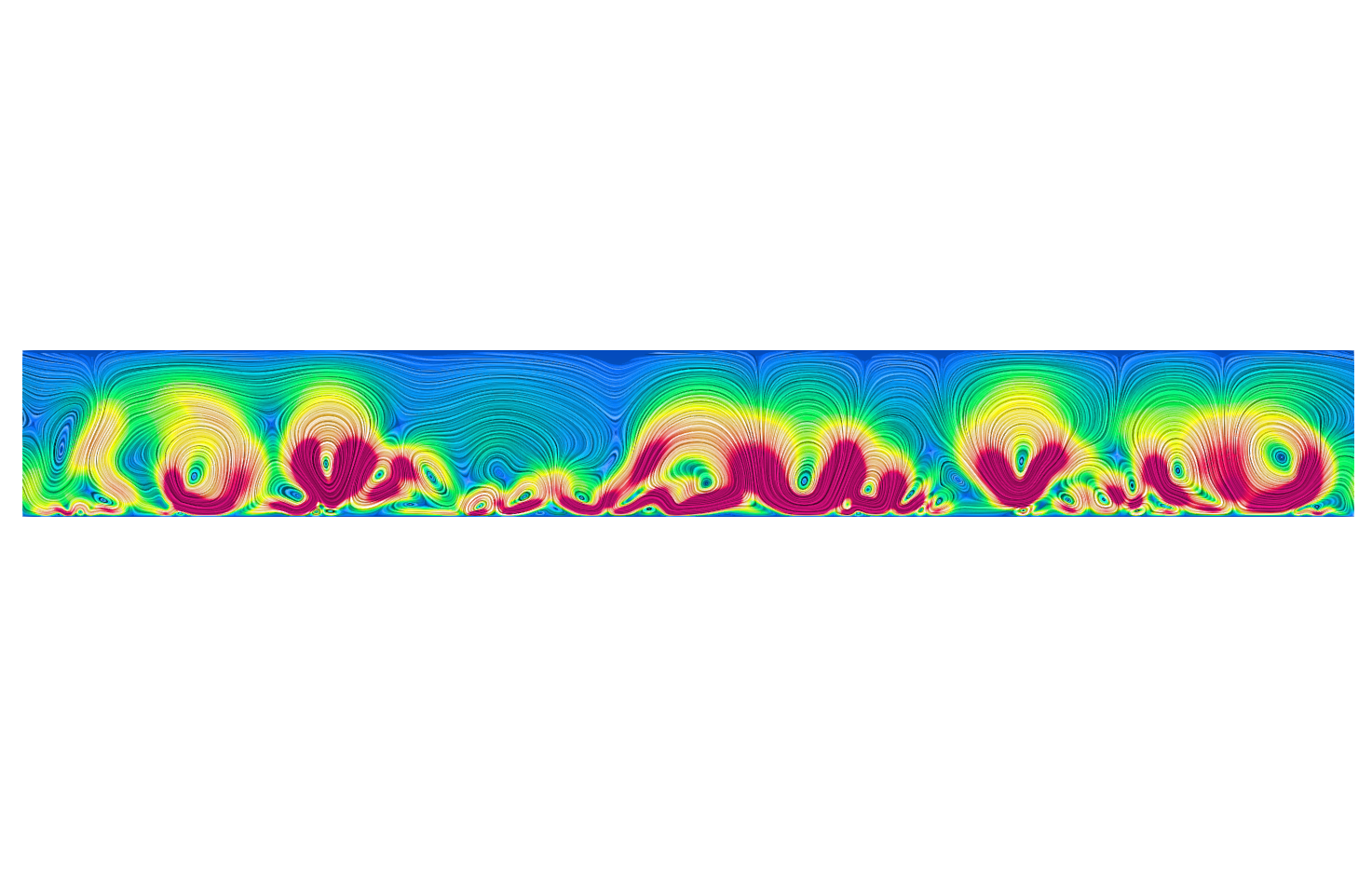}}

\\ \hline
\rotatebox[origin=c]{90}{t = 0.01} & \parbox[c]{\figwidthCase em}{
\includegraphics[width=\linewidth,trim=100 550 100 500,clip]{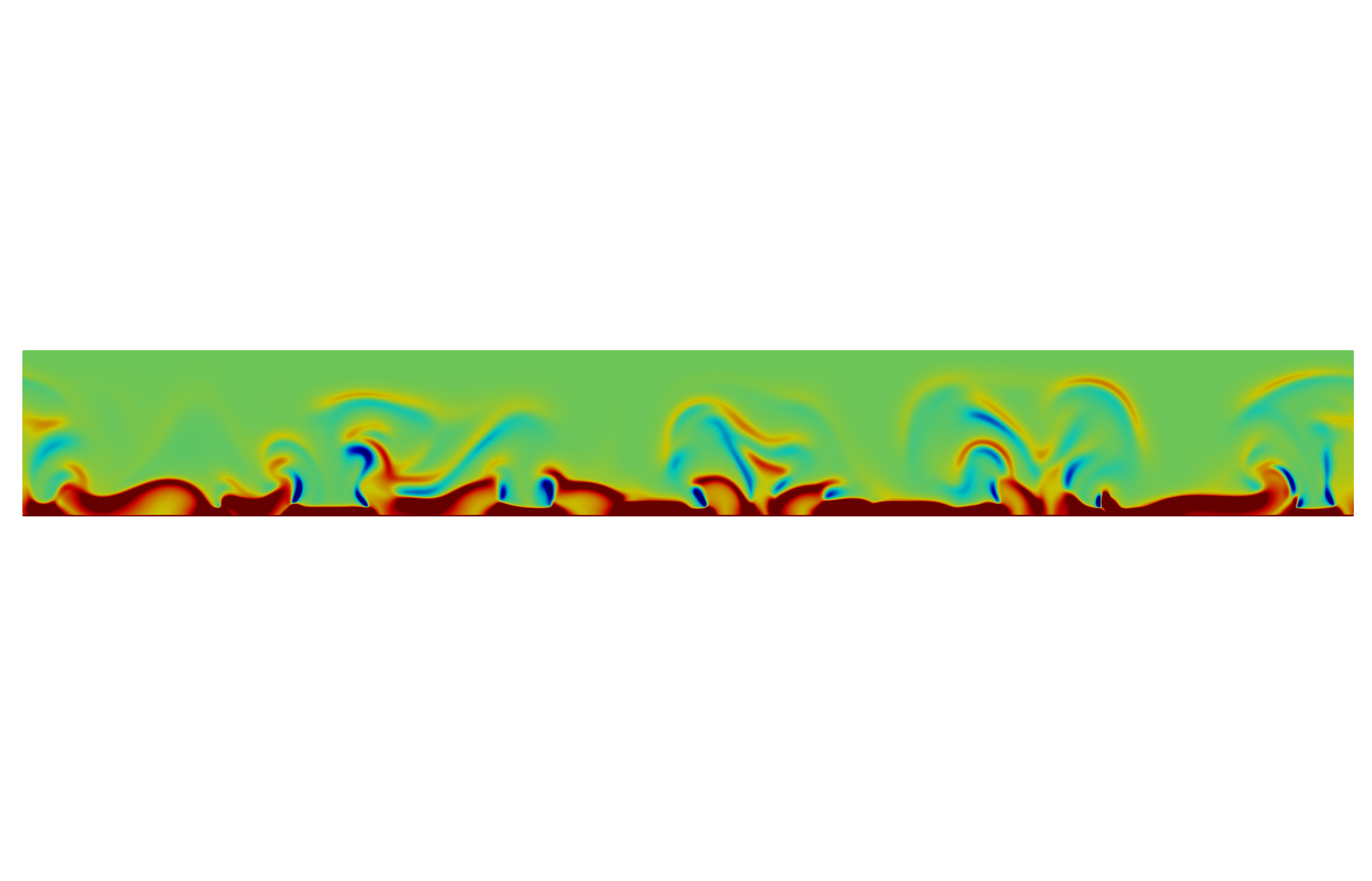}}  
& 
\parbox[c]{\figwidthCase em}{
\includegraphics[width=\linewidth,trim=50 275 50 250,clip]{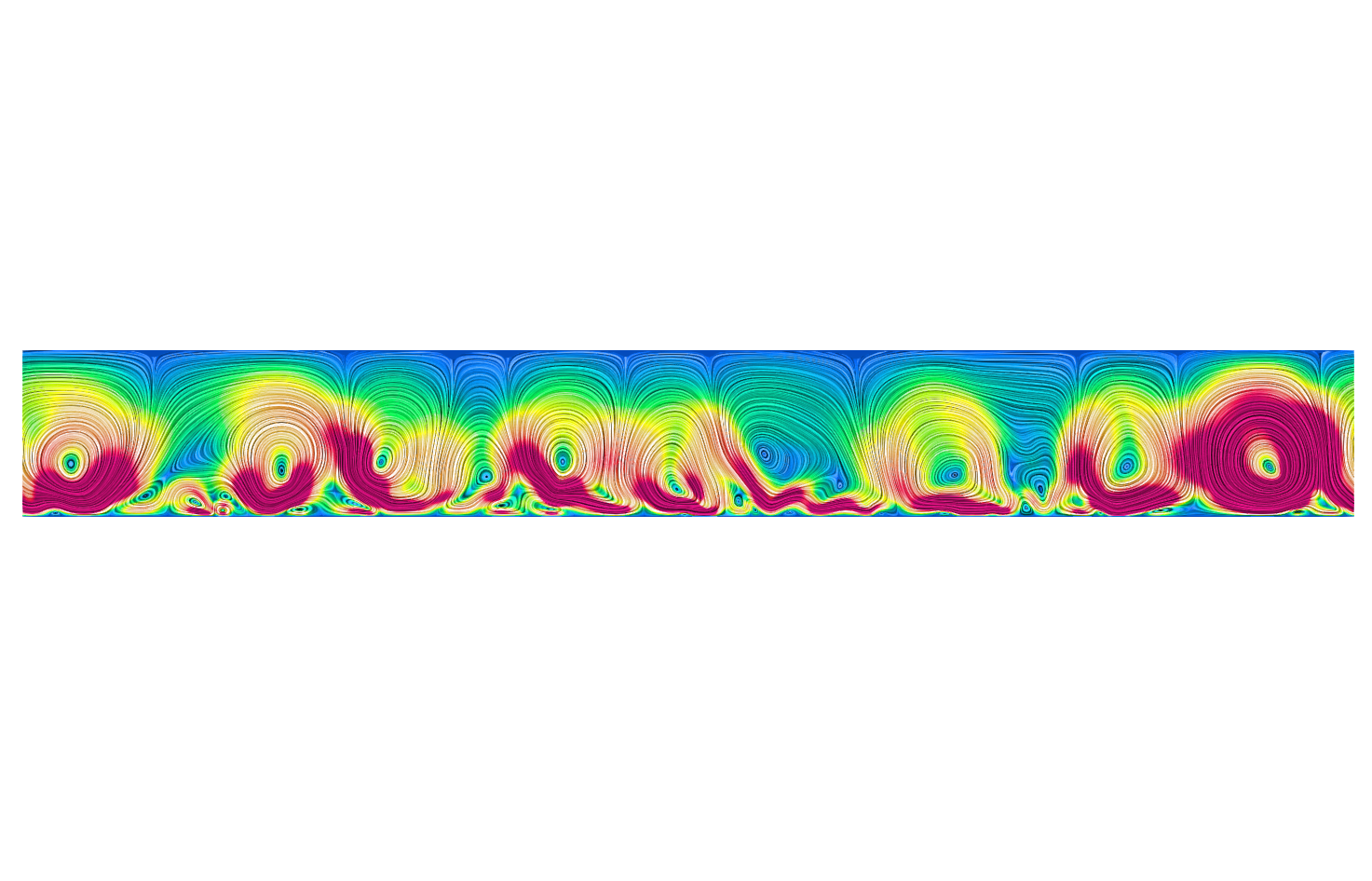}}

\\ \hline
\rotatebox[origin=c]{90}{t = 0.05} & 
\parbox[c]{\figwidthCase em}{
\includegraphics[width=\linewidth,trim=100 550 100 500,clip]{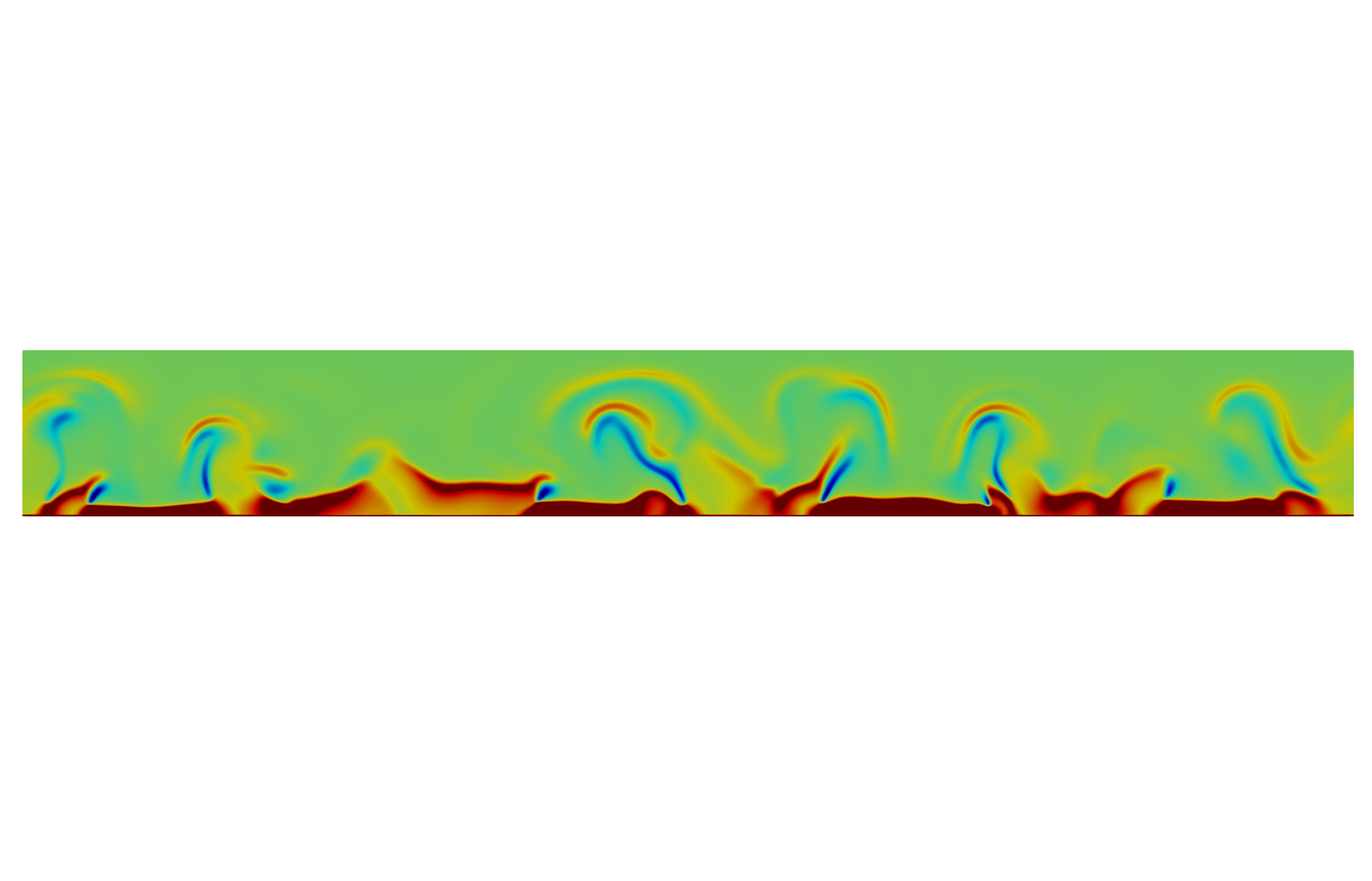}} 
&  
\parbox[c]{\figwidthCase em}{
\includegraphics[width=\linewidth,trim=50 275 50 250,clip]{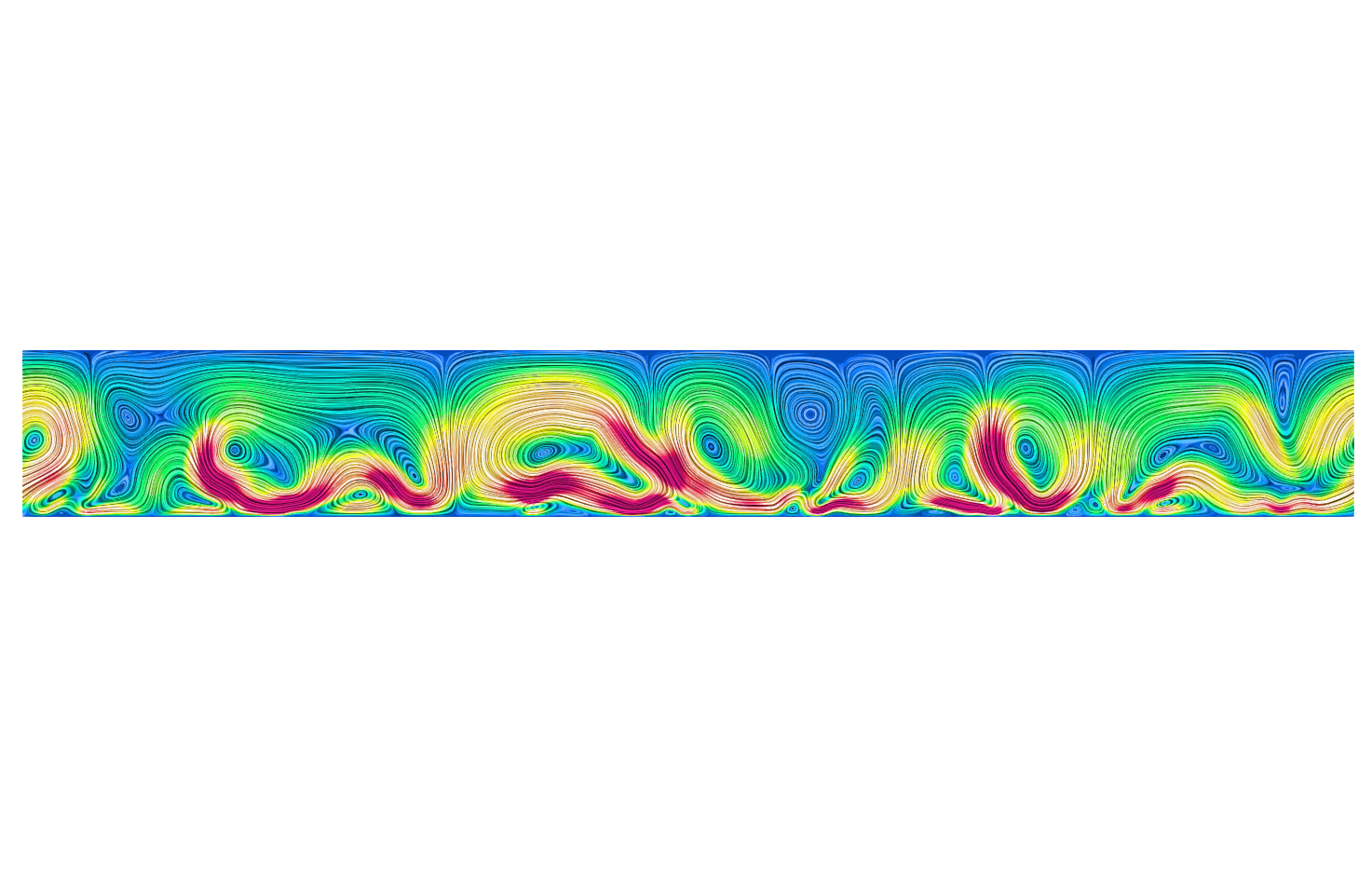}}

\\ \hline
\rotatebox[origin=c]{90}{t = 0.1} & 
\parbox[c]{\figwidthCase em}{
\includegraphics[width=\linewidth,trim=100 550 100 500,clip]{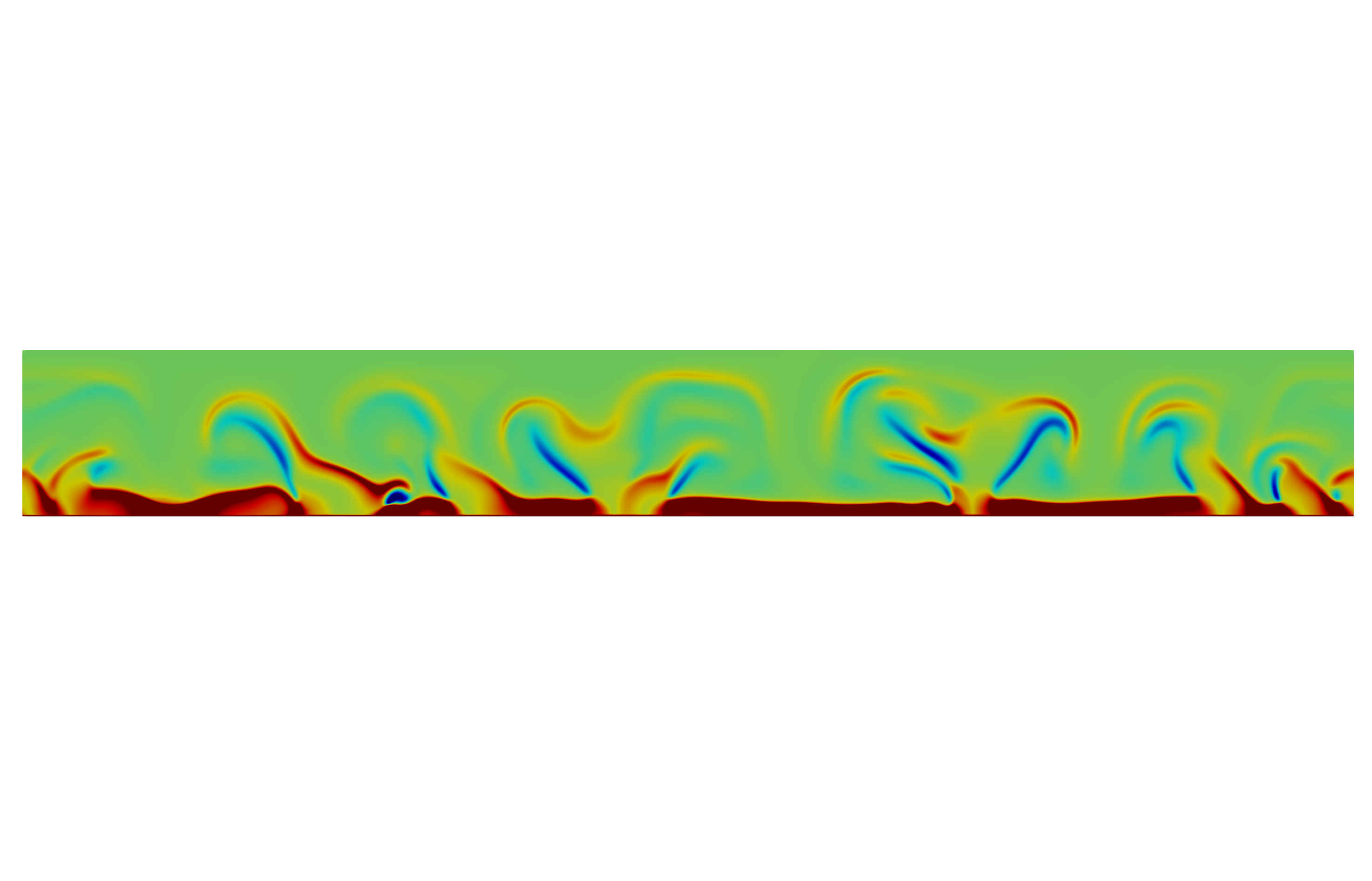}} 
&  
\parbox[c]{\figwidthCase em}{
\includegraphics[width=\linewidth,trim=50 275 50 250,clip]{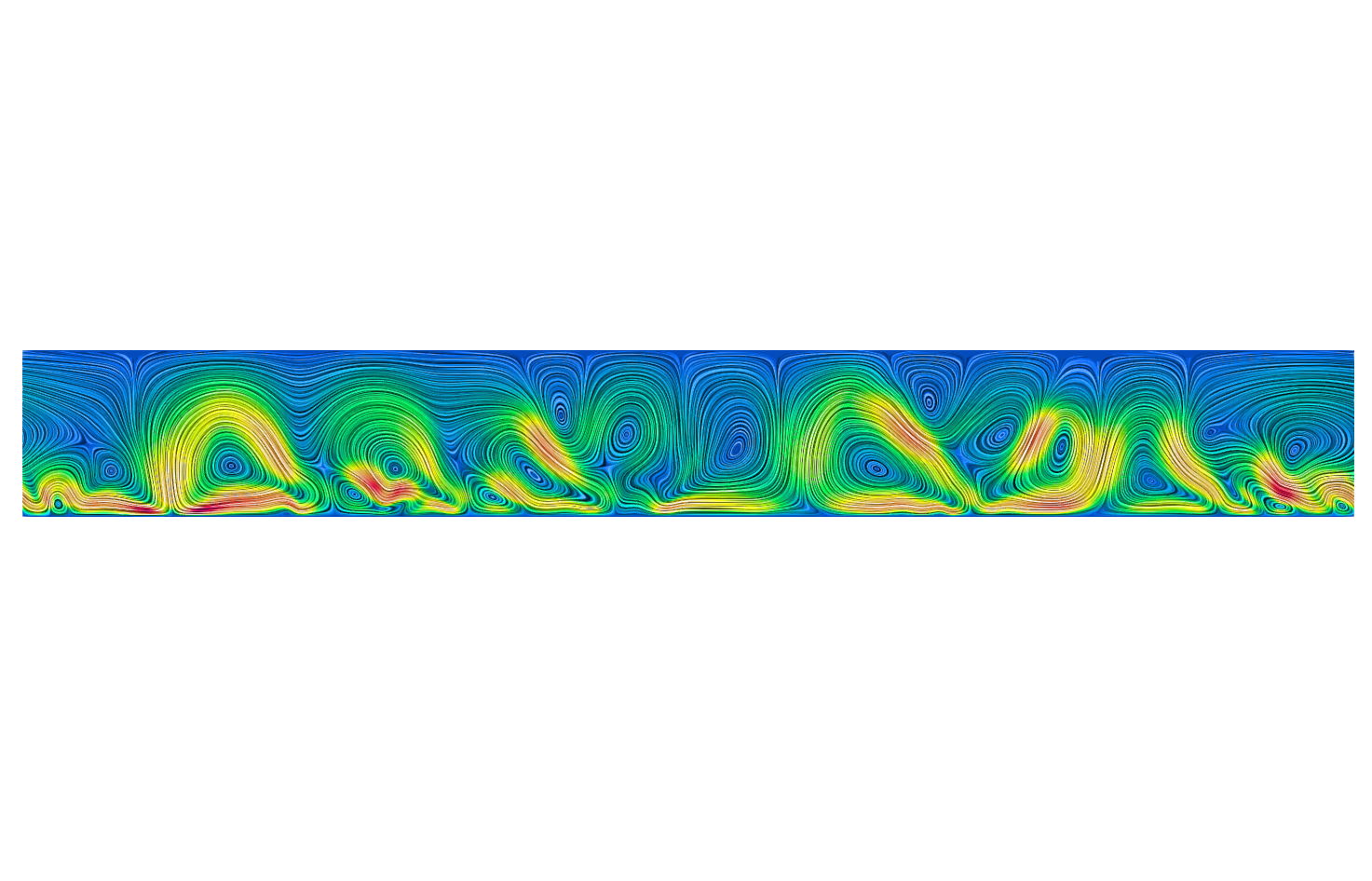}}   
\\ \hline
\end{tabular}
\caption{Contour plots of charge density and velocity magnitudes  at different times for $\Delta \phi = $ 120. }
\label{fig:plot120V}
\end{figure}

{}

\newpage 
\subsection{Packed bed ion concentration polarization: Complex geometries, small EDLs, and upstream vortex formation}\label{sec:microbeads}
Our final example illustrates a simulation of ion concentration polarization (ICP) -- a phenomenon used in a wide array of electrochemical unit operations -- in a realistic, state-of-art engineering device. The device consists of a packed bed of metallic beads that act as a set of electrodes embedded within a microfluidic device. This device was recently studied experimentally by \citet{berzina2022out} to circumvent two challenges faced by conventional ICP occurring at a planar electrode, detailed in the next paragraph. The device configuration offers a challenging canonical simulation example due to (a) the complex geometries of the electrodes involved, (b) the steep gradients (especially near these complex geometries), and small EDL's produced by the operating conditions, (c) the ensuing channel scale vortical structures that have to be accurately resolved, and (d) the long time horizon over which the current measurements stabilize. 


There has been a growing interest in ICP due to its ability to enrich and separate charged species. The key idea of ICP is to manipulate target analytes by leveraging a high electric field in the ion depletion zone (IDZ) created by selective charge transport. Previous studies have successfully demonstrated a variety of applications, including seawater desalination~\citep{kim2010direct}, hemodialysis~\citep{berzina2018electrokinetic}, analyte preconcentration~\citep{ko2012nanofluidic}, and manipulation of neutral species following their partition into ionic micelles ~\citep{berzina2019continuous}, all using planar ICP. However, two challenges preclude the scaling up of ICP from enabling commercial deployment. First, an IDZ formed at a planar membrane or electrode may not extend across the full channel cross-section under the flow rate employed for focusing, thereby allowing the analyte to “leak” past the IDZ. Second, within the IDZ, large fluid vortices lead to mixing, which decreases the efficiency of analyte enrichment and worsens with increased channel dimensions. A promising approach has been to move away from a planar electrode (along the bottom of the channel) to using engineered, porous electrodes across the channel cross-section. Recently, \citet{berzina2022out} introduced 3-D packed bed electrodes that successfully control unwanted instabilities and extend high electric field area for the enrichment of bioanalytes before detection. In a follow-up study by \citet{berzina2023electrokinetic}, an additional bed of packed bioconjugated beads were incorporated for DNA capture. We deploy our computational framework on this configuration to explore the impact of the packed bed electrodes on the current generation, concentration distribution, electric field extent, and vorticity patterns.

\begin{figure}
\centering
\begin{subfigure}{.75\textwidth}
\centering
\includegraphics[width = \linewidth]{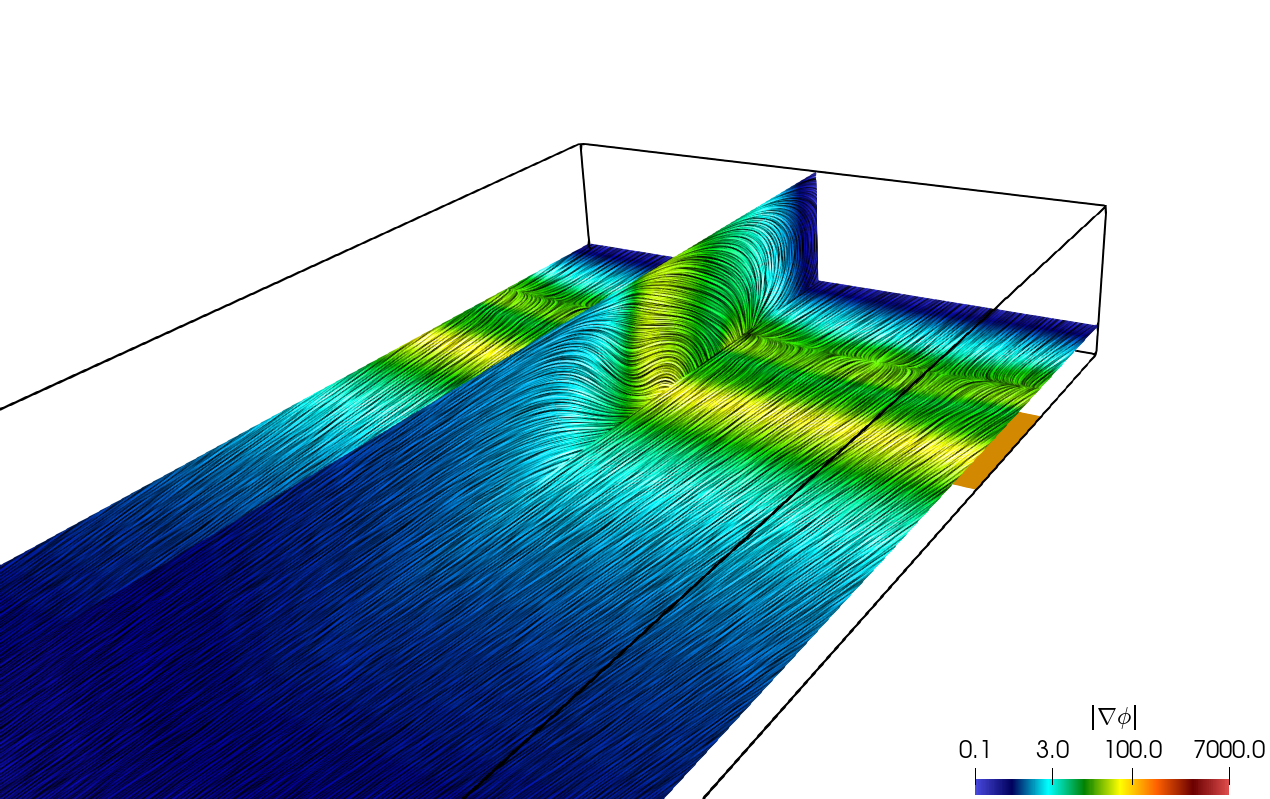}
\subcaption{}
\end{subfigure}
\begin{subfigure}{.75\textwidth}
\centering
\includegraphics[width = \linewidth]{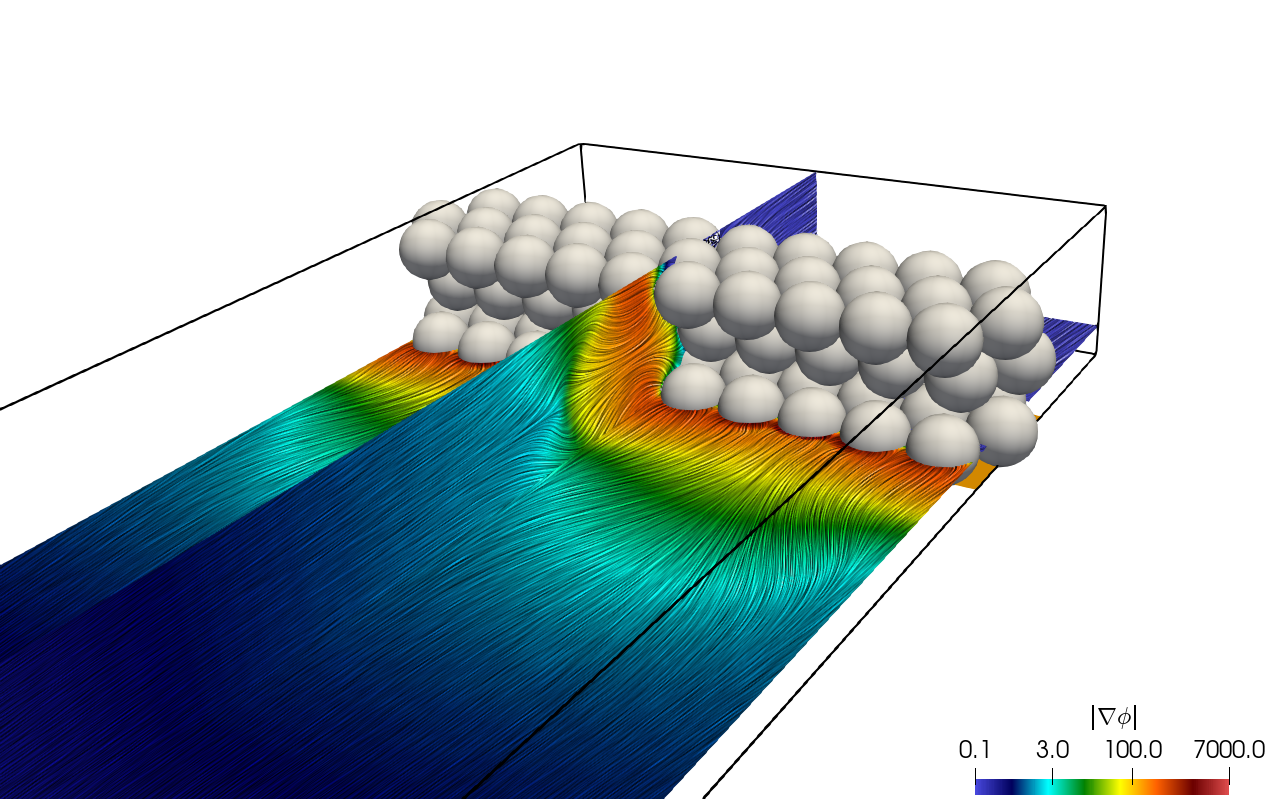}
\subcaption{}
\end{subfigure}
\caption{The magnitude of the electric field (a) only with the plate electrode, and (b) with added 3-D electrode bed.}
\label{fig:sphere_efield}
\end{figure}

\begin{figure}
\centering
\begin{subfigure}{.75\textwidth}
\centering
\includegraphics[width = \linewidth]{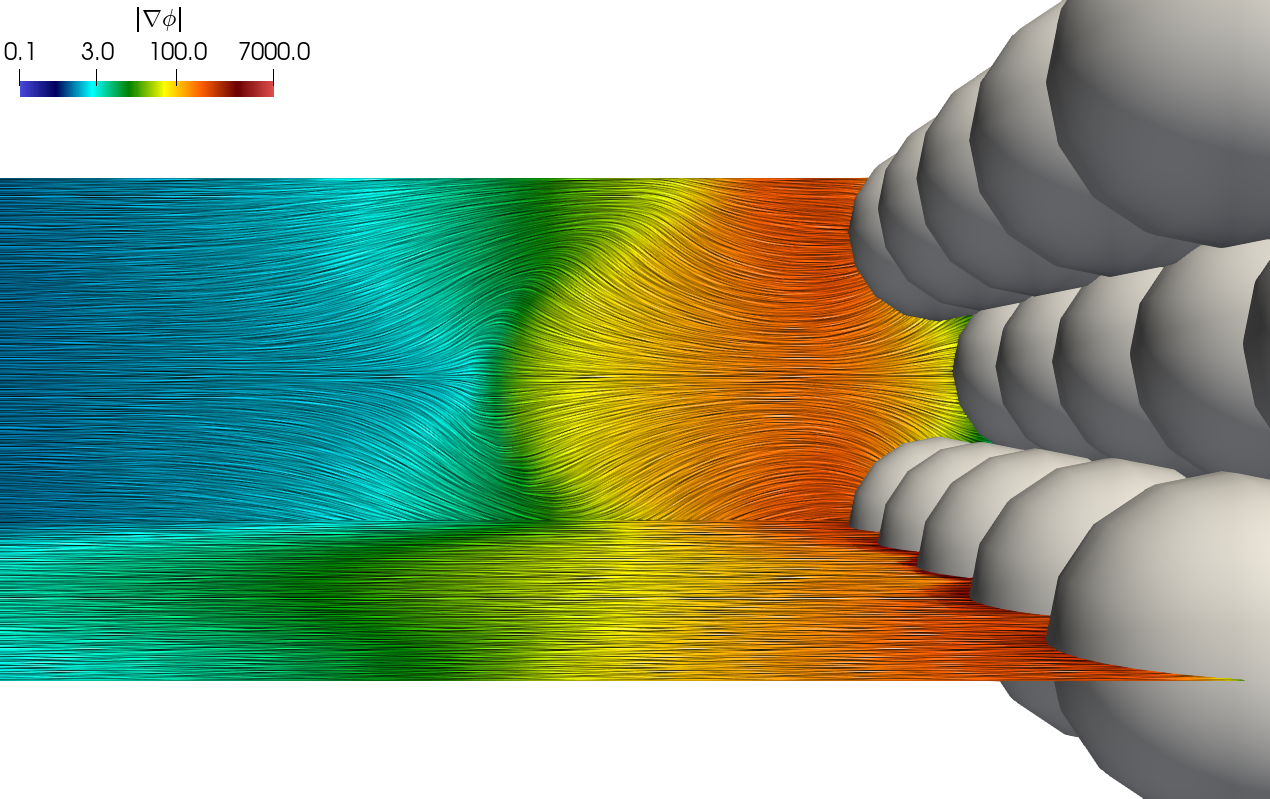}
\subcaption{}
\end{subfigure}
\begin{subfigure}{.75\textwidth}
\centering
\includegraphics[width = \linewidth]{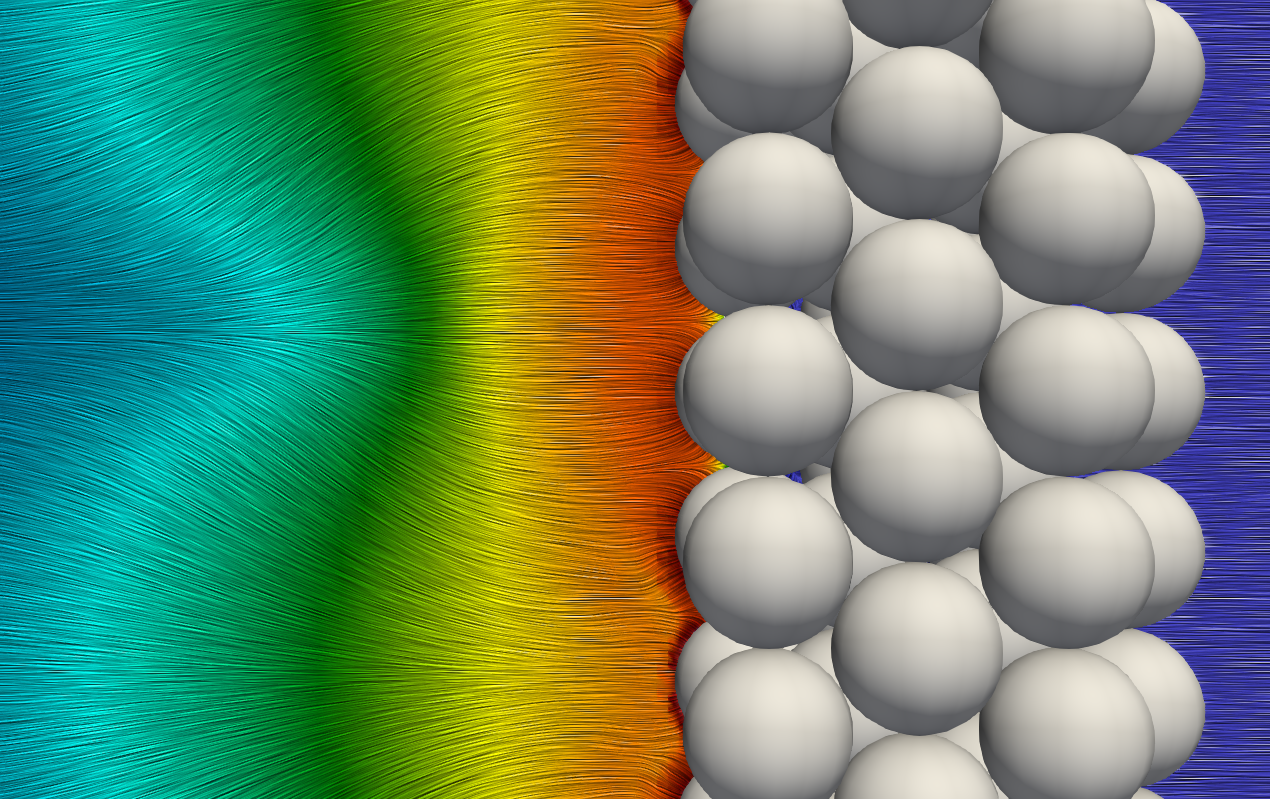}
\subcaption{}
\end{subfigure}
\caption{The magnified view of the electric field near the packed bed from the side (a) and top (b)}
\label{fig:sphere_electricMag}
\end{figure}

\begin{figure}
\centering
\begin{subfigure}{.75\textwidth}
\centering
\includegraphics[width = \linewidth]{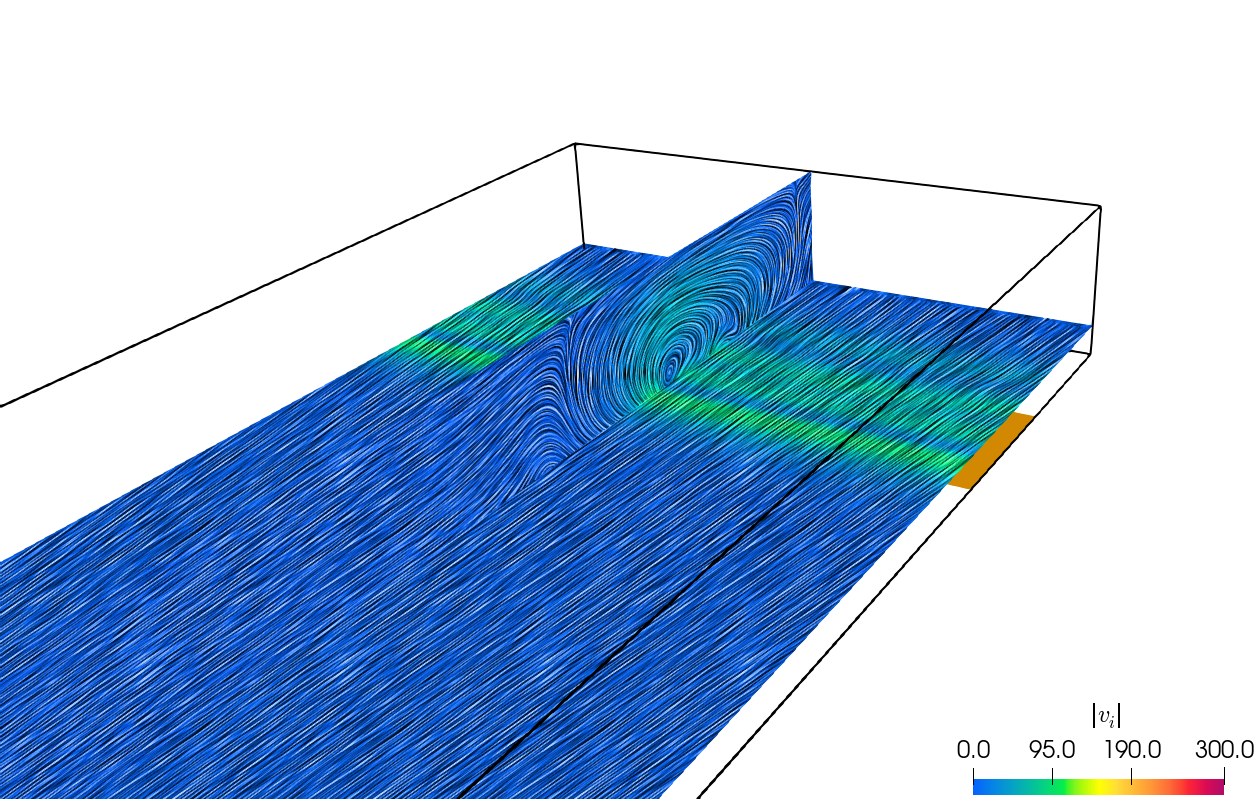}
\subcaption{}
\end{subfigure}
\begin{subfigure}{.75\textwidth}
\centering
\includegraphics[width = \linewidth]{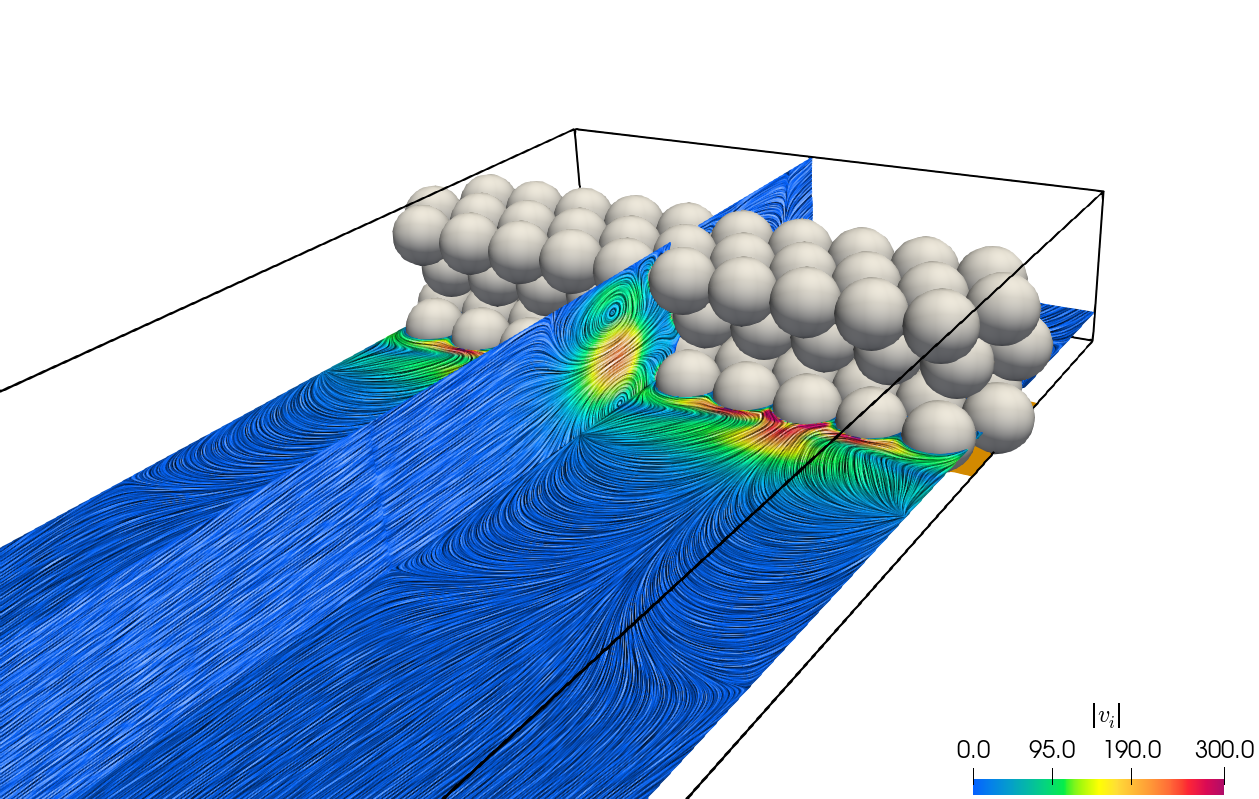}
\subcaption{}
\end{subfigure}
\caption{The flow streamline and magnitude around the electrodes. (a) only with the plate electrode, and (b) with added 3-D electrode bed.}
\label{fig:sphere_flow}
\end{figure}

The domain consists of a cuboidal geometry with packed bed electrodes -- made up of spherical beads -- downstream of the channel inlet. Fluid is pushed from left to right, and a potential difference is maintained. The length, width, and height of the simulation domain are $32 \times 1 \times 1$. We only simulate a part of the width considering symmetry. On the floor of the packed bead bed, a planar electrode is located. The geometry around the packed beds is shown in~\figref{fig:sphere_efield} - ~\figref{fig:sphere_flow}.

\begin{figure}
\centering
\begin{subfigure}{0.3\textwidth}
\includegraphics[width=\textwidth]{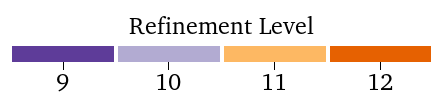}
\label{fig:colourmap}
\end{subfigure}\\
\begin{subfigure}{0.5\textwidth}
\includegraphics[width=\textwidth]{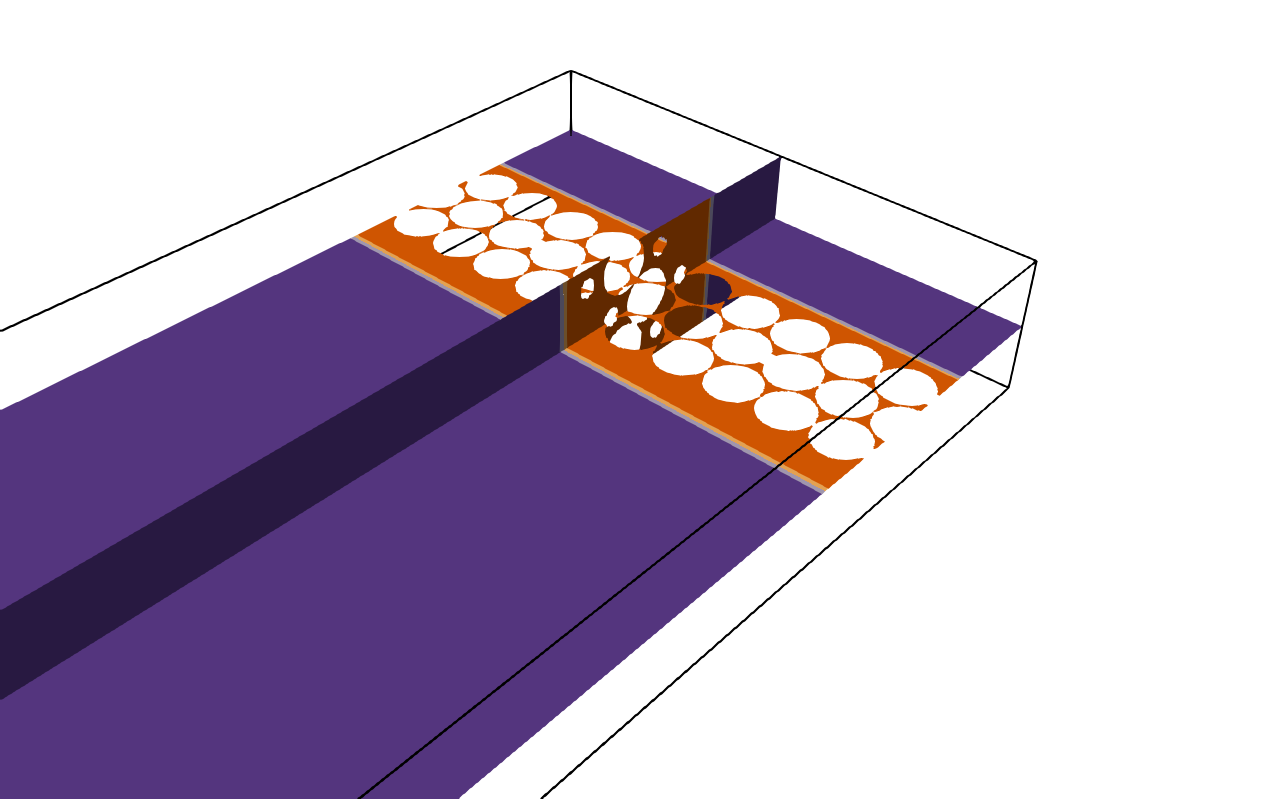}
\caption{}
\label{fig:iso}
\end{subfigure}\\
\begin{subfigure}{0.4\textwidth}
\includegraphics[width=\textwidth]{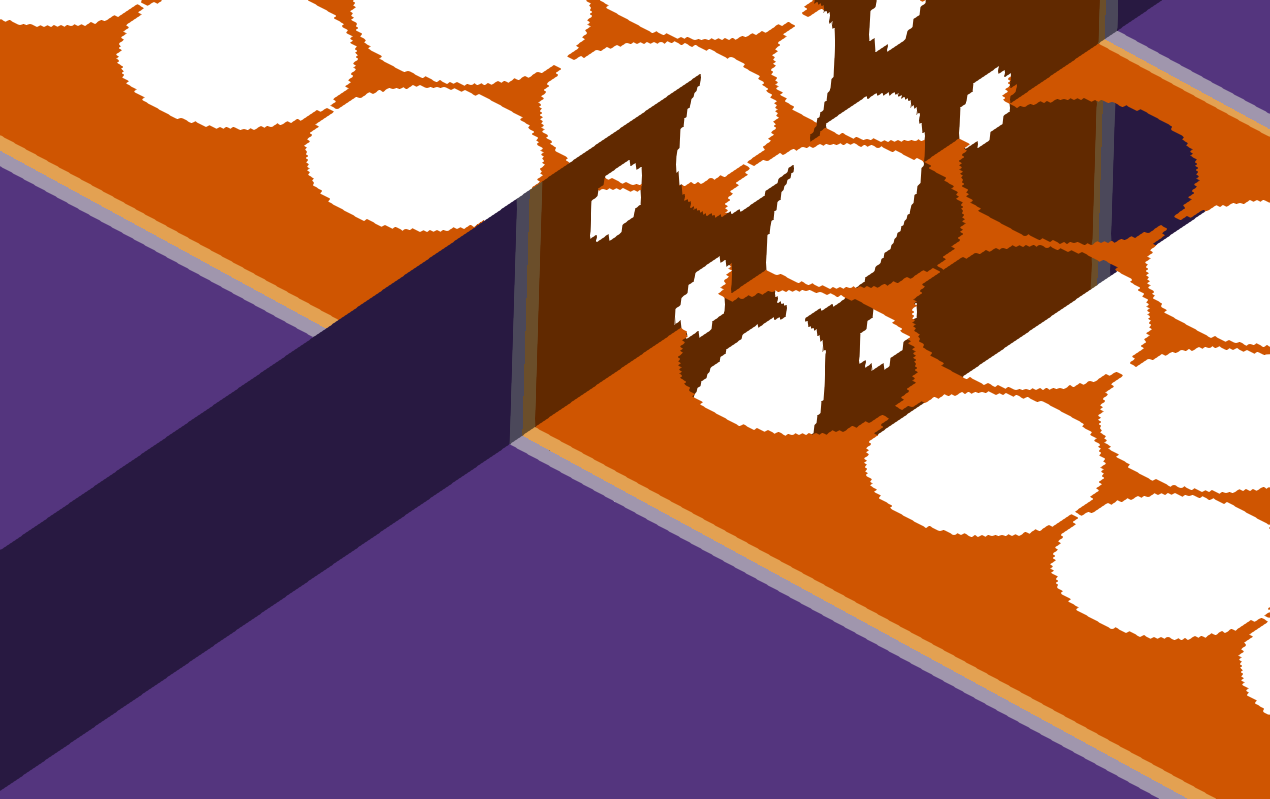}
\caption{}
\label{fig:isoZoom}
\end{subfigure}
\begin{subfigure}{0.4\textwidth}
\includegraphics[width=\textwidth]{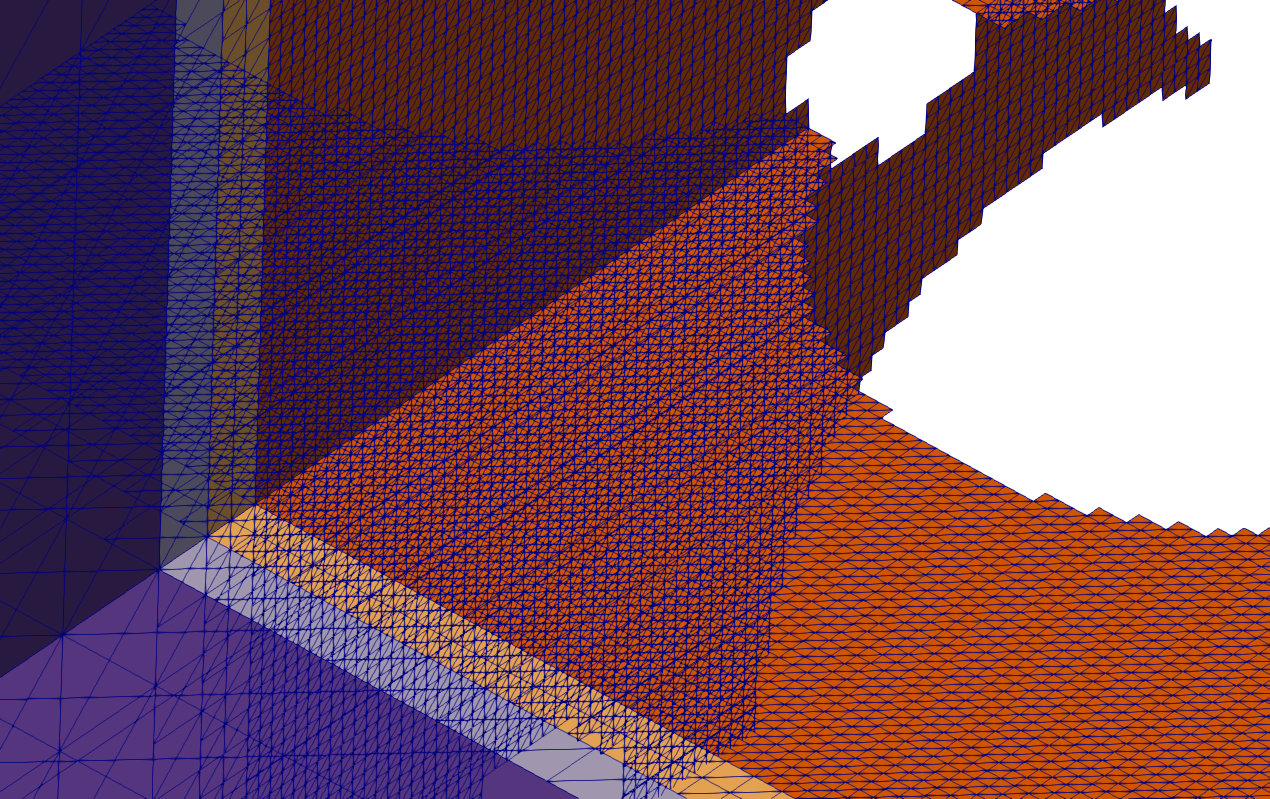}
\caption{}
\label{fig:isoZoomMesh}
\end{subfigure}
\begin{subfigure}{0.4\textwidth}
\includegraphics[width=\textwidth]{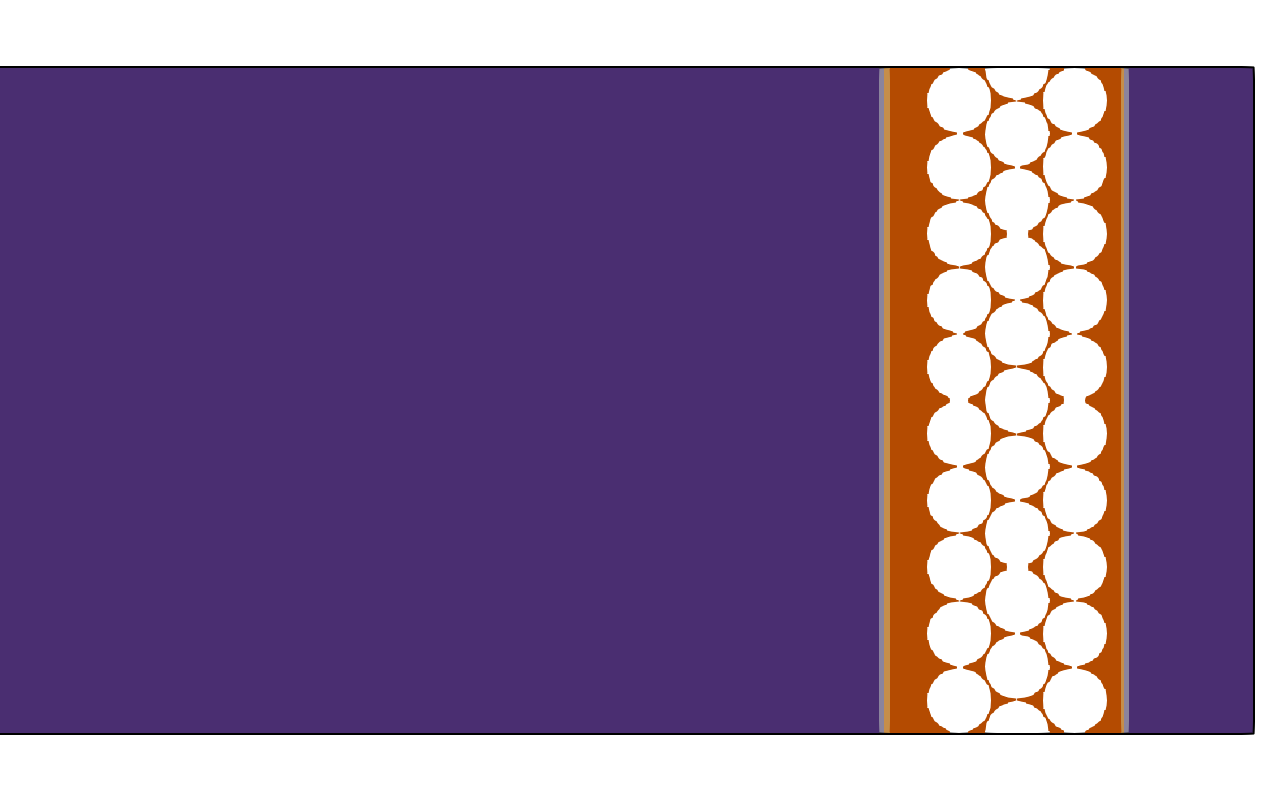}
\caption{}
\label{fig:top}
\end{subfigure}
\begin{subfigure}{0.4\textwidth}
\includegraphics[width=\textwidth]{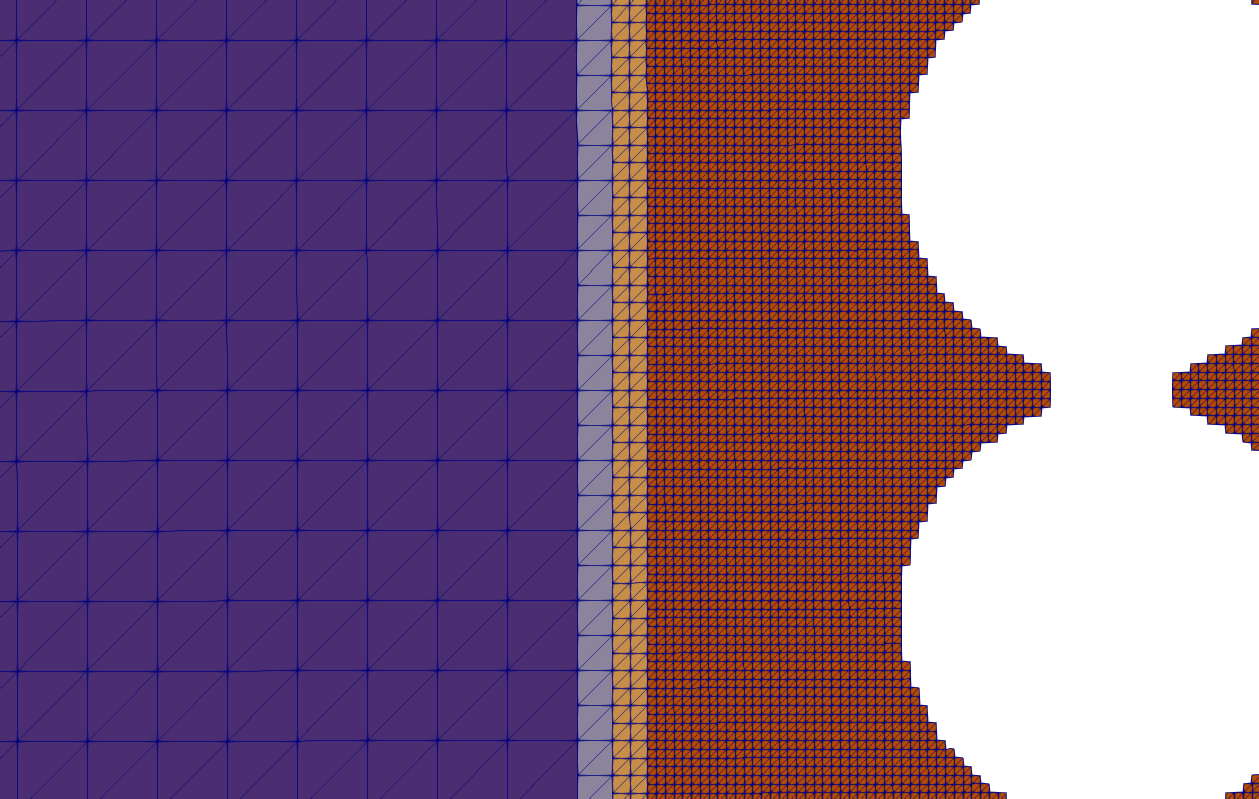}
\caption{}
\label{fig:topZoomMesh}
\end{subfigure}
\caption{Mesh refinement of pack bed ion concentration polarization.}
\label{fig:zoomedMesh}
\end{figure}

A challenging aspect in such simulations is the need to construct an adaptively refined mesh of this complex geometry. Here, the octree-based adaptive meshing feature in our framework allows us to efficiently construct a high-quality mesh of the system. \figref{fig:zoomedMesh}~illustrates the mesh refinement near the packed bed on a vertical plane. The refinement level ranges from 9 (bulk) to 13 (electrodes), which corresponds to the element sizes ranging from $32/2^{9} \sim 0.0625$ and $32/2^{13} \sim 0.0039$, which captures physics in a wide range of length scales from the EDL around the "rough" electrode surface to the fluid flow at bulk. The non-dimensional Debye length, in this case, is $0.0036$. 
The details of boundary conditions are listed in \tabref{tab:BC-Sphere}. As before, the initial conditions were obtained by solving only the PNP equations until the EDL was fully developed. Then, the fully coupled NS-PNP equations were solved.  
\begin{table}[htb!]
\centering
\renewcommand{\arraystretch}{1.2}
\begin{tabular}{|ccccc|}
\hline
\multicolumn{3}{|c|}{} & \multicolumn{1}{c|}{\textbf{Type}} & \textbf{Value} \\ \hline \hline
\multicolumn{1}{|c|}{\multirow{13}{*}{NS}} & \multicolumn{1}{c|}{\multirow{4}{*}{\begin{tabular}[c]{@{}c@{}}Packed Bed Electrode, \\ Top \& Bottom Walls\\ (no - slip)\end{tabular}}} & \multicolumn{1}{c|}{u} & \multicolumn{1}{c|}{Dirichlet} & 0 \\ \cline{3-5} 
\multicolumn{1}{|c|}{} & \multicolumn{1}{c|}{} & \multicolumn{1}{c|}{v} & \multicolumn{1}{c|}{Dirichlet} & 0 \\ \cline{3-5} 
\multicolumn{1}{|c|}{} & \multicolumn{1}{c|}{} & \multicolumn{1}{c|}{w} & \multicolumn{1}{c|}{Dirichlet} & 0 \\ \cline{3-5} 
\multicolumn{1}{|c|}{} & \multicolumn{1}{c|}{} & \multicolumn{1}{c|}{p} & \multicolumn{1}{c|}{Neumann} & 0 \\ \cline{2-5} 
\multicolumn{1}{|c|}{} & \multicolumn{1}{c|}{\multirow{4}{*}{\begin{tabular}[c]{@{}c@{}}Inlet\end{tabular}}} & \multicolumn{1}{c|}{u} & \multicolumn{1}{c|}{Dirichlet} & 0.371 \\ \cline{3-5} 
\multicolumn{1}{|c|}{} & \multicolumn{1}{c|}{} & \multicolumn{1}{c|}{v} & \multicolumn{1}{c|}{Dirichlet} & 0 \\ \cline{3-5} 
\multicolumn{1}{|c|}{} & \multicolumn{1}{c|}{} & \multicolumn{1}{c|}{w} & \multicolumn{1}{c|}{Dirichlet} & 0 \\ \cline{3-5} 
\multicolumn{1}{|c|}{} & \multicolumn{1}{c|}{} & \multicolumn{1}{c|}{p} & \multicolumn{1}{c|}{Neumann} & 0 \\ \cline{2-5} 
\multicolumn{1}{|c|}{} & \multicolumn{1}{c|}{\multirow{4}{*}{\begin{tabular}[c]{@{}c@{}}Outlet\end{tabular}}} & \multicolumn{1}{c|}{u} & \multicolumn{1}{c|}{Neumann} & 0 \\ \cline{3-5} 
\multicolumn{1}{|c|}{} & \multicolumn{1}{c|}{} & \multicolumn{1}{c|}{v} & \multicolumn{1}{c|}{Neumann} & 0 \\ \cline{3-5} 
\multicolumn{1}{|c|}{} & \multicolumn{1}{c|}{} & \multicolumn{1}{c|}{w} & \multicolumn{1}{c|}{Neumann} & 0 \\ \cline{3-5} 
\multicolumn{1}{|c|}{} & \multicolumn{1}{c|}{} & \multicolumn{1}{c|}{p} & \multicolumn{1}{c|}{Dirchlet} & 0 \\ \cline{2-5} 
\multicolumn{1}{|c|}{} & \multicolumn{1}{c|}{Side Walls} & \multicolumn{3}{c|}{Periodic} \\ \hline \hline 
\multicolumn{1}{|c|}{\multirow{13}{*}{PNP}} & \multicolumn{1}{c|}{\multirow{3}{*}{\begin{tabular}[c]{@{}c@{}}Packed Bed \& Planar\\ Electrodes \end{tabular}}} & \multicolumn{1}{c|}{$\phi$} & \multicolumn{1}{c|}{Dirichlet} & 0 \\ \cline{3-5} 
\multicolumn{1}{|c|}{} & \multicolumn{1}{c|}{} & \multicolumn{1}{c|}{$c_{+}$} & \multicolumn{1}{c|}{Dirichlet} & 0 \\ \cline{3-5} 
\multicolumn{1}{|c|}{} & \multicolumn{1}{c|}{} & \multicolumn{1}{c|}{$c_{-}$} & \multicolumn{1}{c|}{Zero Flux} & 0 \\ \cline{2-5} 
\multicolumn{1}{|c|}{} & \multicolumn{1}{c|}{\multirow{3}{*}{\begin{tabular}[c]{@{}c@{}}Inlet\end{tabular}}} & \multicolumn{1}{c|}{$\phi$} & \multicolumn{1}{c|}{Dirichlet} & 277.1 \\ \cline{3-5} 
\multicolumn{1}{|c|}{} & \multicolumn{1}{c|}{} & \multicolumn{1}{c|}{$c^{+}$} & \multicolumn{1}{c|}{Dirichlet} & 1 \\ \cline{3-5} 
\multicolumn{1}{|c|}{} & \multicolumn{1}{c|}{} & \multicolumn{1}{c|}{$c^{-}$} & \multicolumn{1}{c|}{Dirichlet} & 1 \\ \cline{2-5} 
\multicolumn{1}{|c|}{} & \multicolumn{1}{c|}{\multirow{3}{*}{\begin{tabular}[c]{@{}c@{}}Outlet\end{tabular}}} & \multicolumn{1}{c|}{$\phi$} & \multicolumn{1}{c|}{Neumann} & 0 \\ \cline{3-5} 
\multicolumn{1}{|c|}{} & \multicolumn{1}{c|}{} & \multicolumn{1}{c|}{$c^{+}$} & \multicolumn{1}{c|}{Neumann} & 0 \\ \cline{3-5} 
\multicolumn{1}{|c|}{} & \multicolumn{1}{c|}{} & \multicolumn{1}{c|}{$c^{-}$} & \multicolumn{1}{c|}{Neumann} & 0 \\ \cline{2-5} 
\multicolumn{1}{|c|}{} & \multicolumn{1}{c|}{\multirow{3}{*}{\begin{tabular}[c]{@{}c@{}}Top \& Bottom Walls \\ (except electrode) \end{tabular}}} & \multicolumn{1}{c|}{$\phi$} & \multicolumn{1}{c|}{Neumann} & 0 \\ \cline{3-5} 
\multicolumn{1}{|c|}{} & \multicolumn{1}{c|}{} & \multicolumn{1}{c|}{$c^{+}$} & \multicolumn{1}{c|}{Zero Flux} & 0 \\ \cline{3-5} 
\multicolumn{1}{|c|}{} & \multicolumn{1}{c|}{} & \multicolumn{1}{c|}{$c^{-}$} & \multicolumn{1}{c|}{Zero Flux} & 0 \\ \cline{2-5} 
\multicolumn{1}{|c|}{} & \multicolumn{1}{c|}{Side Walls} & \multicolumn{3}{c|}{Periodic} \\ \hline
\end{tabular}
\caption{Boundary condition for 3-D electrode packed bed ICP}
\label{tab:BC-Sphere}
\end{table}

\tabref{table: mesh_sphere} briefs the mesh information and the total time taken. It is worth noting the importance of adaptivity. The ability to selectively refine only the regions of interest resulted mesh comprised of around 1.6 M elements and a total wall time of approximately two days. In the absence of any adaptivity, a resultant mesh with a refinement level of 13 would result in approximately 536 million elements. Such a drastic increase in the number of elements would render the problem expensive in a reasonable amount of time.

\begin{table}[h]
\centering
\begin{tabular}{|c|c|c|}
\hline
\textbf{Mesh size} & \textbf{Number of MPI task} & \textbf{Wall time} \\ \hline
\multirow{2}{*}{\begin{tabular}[c]{@{}c@{}}1.6 M\\ (536 M equivalent)\end{tabular}} & \multirow{2}{*}{\begin{tabular}[c]{@{}c@{}}544 CPU cores\\ (8 KNL nodes on Stampede2)\end{tabular}} & \multirow{2}{*}{2 days} \\
&  &  \\ \hline
\end{tabular}
\caption{Mesh details and time taken for the Packed bed ion concentration polarization.}
\label{table: mesh_sphere}
\end{table}

\newpage

\begin{figure}
\centering
\begin{subfigure}{.75\textwidth}
\centering
\includegraphics[width = \linewidth]{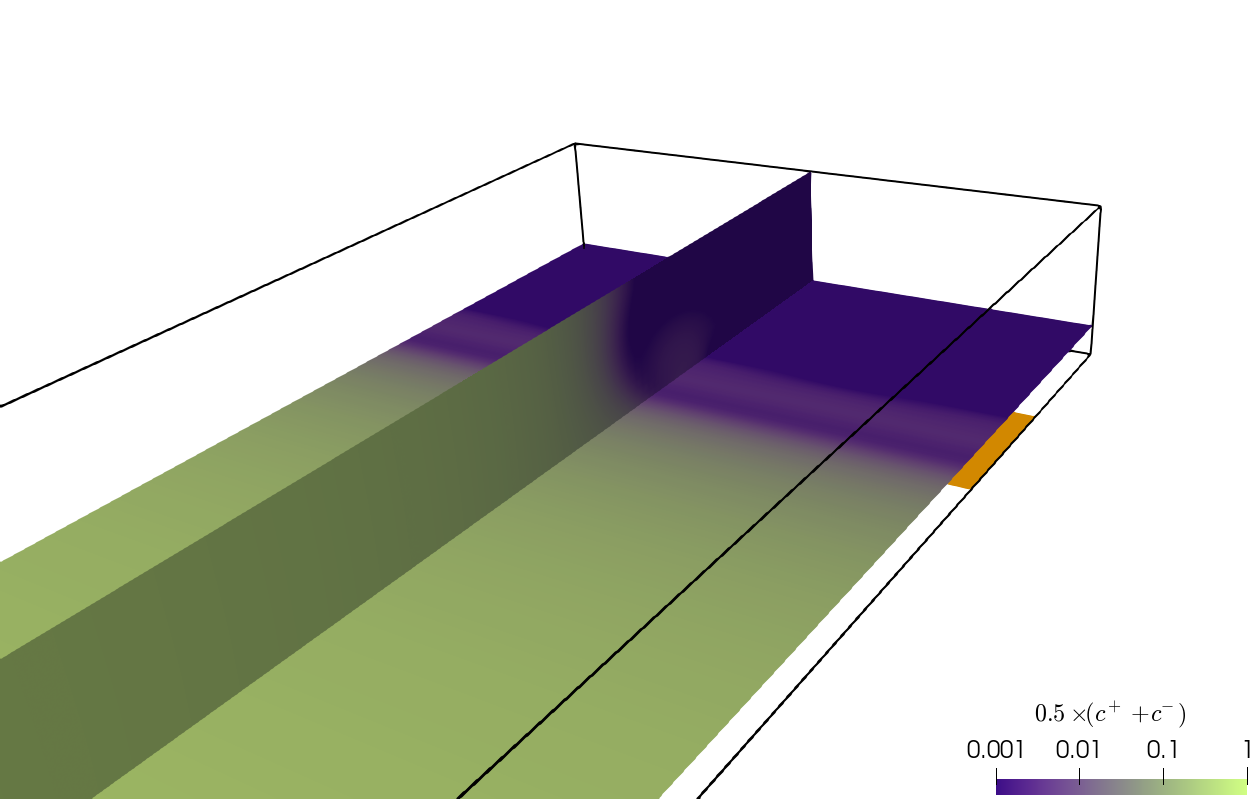}
\subcaption{}
\end{subfigure}
\begin{subfigure}{.75\textwidth}
\centering
\includegraphics[width = \linewidth]{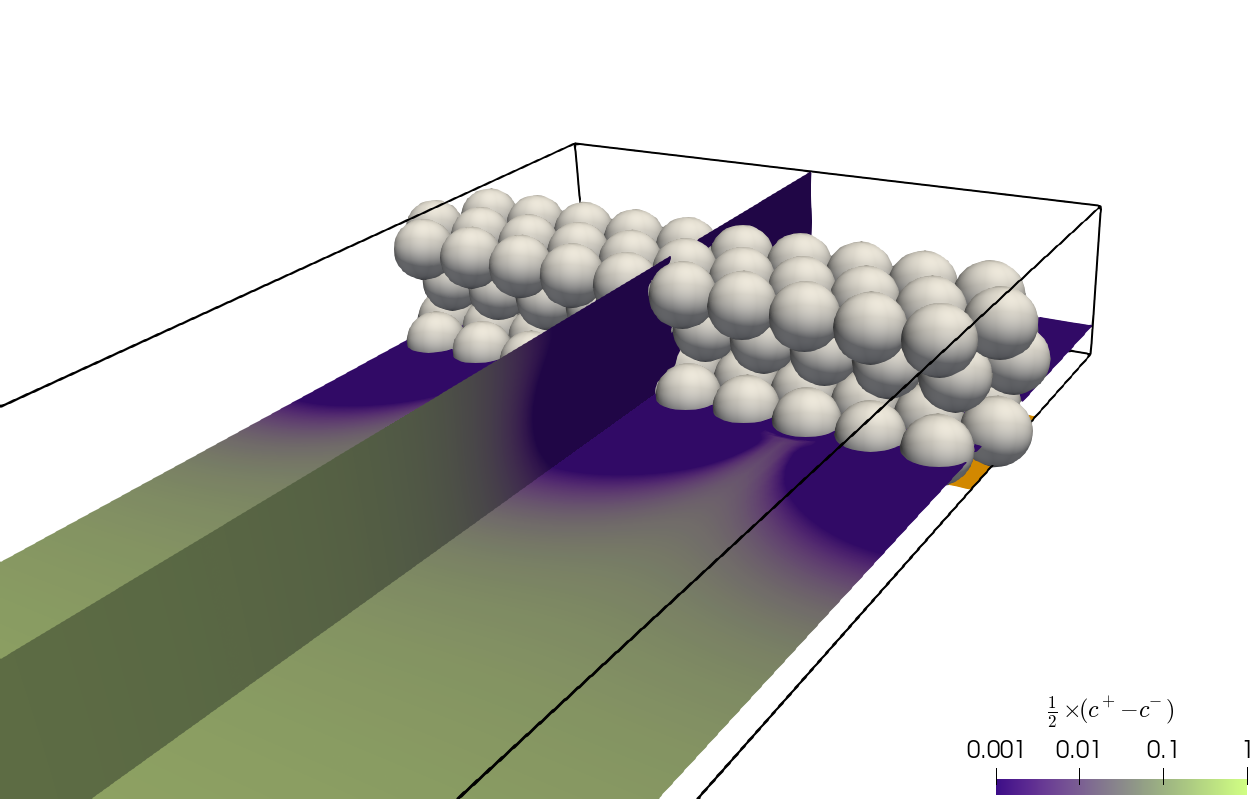}
\subcaption{}
\end{subfigure}
\caption{The ionic strength (average concentration in this case) (a) only with the plate electrode, and (b) with added 3-D electrode bed.}
\label{fig:sphere_ionicstrength}
\end{figure}

We compare the spatial variations in the electric field, charge densities, and flow fields between two cases: a microchannel with a packed bead electrode as described above versus a microchannel with only a planar electrode. We first compare the magnitude and extent of the electric field in ~\figref{fig:sphere_efield}. We see that the presence of the microbeads causes the local enhancement of the electric field to extend throughout the channel cross-section rather than being concentrated only in the bottom region. We note that these observations are qualitatively consistent with the ones reported in the experimental study of~\citet{berzina2022out}. Figure.~\ref{fig:sphere_electricMag} shows a magnified view of the magnitude of the electric field near the bead bed.

Next, we compare the flow fields between the two cases in~\figref{fig:sphere_flow}. Note the set of vortices that span the channel width, with each vortex pair nearly spanning the channel height. This is in contrast to a single large vortex in the planar electrode case. The increased region over which the electric field acts induces a larger velocity field upstream of the packed bed. However, experimental evidence suggests that this does not cause any electroconvective instabilities. The absence of electroconvective instabilities is confirmed by visualizing the ionic strength in ~\figref{fig:sphere_ionicstrength}. 

\section{Conclusions and future work}
\label{sec:conclusions}

Direct numerical simulations of the coupled NS-PNP equations are challenging due to the 1) wide ranges of time and length scales exhibited by the multiscale physics, 2) The chaotic nature of the flow triggers structures over a wide range of scales in all field quantities,  
and 3) geometric complexities that escalate the computational cost. In this work, we presented a finite element-based numerical framework that addresses those challenges by incorporating 1) Variational multiscale pressure stabilization, 2) hybridized semi-implicit and fully implicit time integration schemes, and 3) highly parallelized adaptive meshing. 

We demonstrate the utility of the framework using a series of numerical examples. The temporal and spatial convergence tests were performed to verify the order of accuracy of the numerical framework. Next, we performed simulations capturing electrohydrodynamic instabilities near an ion-selective surface. These simulations showed excellent agreement with previously reported direct numerical simulations. Moreover, the proposed numerical scheme enabled using time steps three orders of magnitude larger than the benchmarking study, allowing significant savings in computational costs. Finally, we deployed the framework for simulating a practical microfluidics application for analyte preconcentration, which involves a complex geometry. The adaptive meshing and numerical approach are able to capture the impact of the complex geometry of 3-D electrodes. The simulation confirmed the effect of the 3-D electrodes on suppressing unwanted vortices and extending the high electric field area, and this phenomenon was also observed in the experimental work. 

In summary, our proposed numerical framework addresses the challenges involved in detailed simulations of electrokinetic phenomena for practical applications. The numerical framework offers robustness for simulations with numerical instabilities, the ability to take large time steps to reach long-time horizon simulations, and adaptive octree-based meshing to resolve geometric complexities, all packaged within a massively scalable software stack. 

\section{Acknowledgements}
The authors acknowledge XSEDE grant number TG-CTS110007 for computing time on TACC Stampede2.  The authors also acknowledge computing allocation through a DD award on TACC Frontera and Iowa State University computing resources. BG, KS, SK, and MAK were funded in part by NSF grants 1935255, and 1855902.

\appendix

\section{Details of solver selection for the numerical experiments}
\label{sec:app_linear_solve}
For the cases presented in \cref{sec:alimani,sec:microbeads} we use the BiCGStab linear solver (a Krylov space solver) with additive Schwarz-based preconditioning.  For better reproduction, the command line options we provide {\sc petsc} are given below which include some commands used for printing some norms as well. The {\sc petsc} options for the Navier-Stokes equations are:
\begin{lstlisting}
solver_options_ns = {
ksp_atol = 1e-10
ksp_rtol = 1e-08
ksp_stol = 1e-14
ksp_max_it = 2000
ksp_type = "bcgs"
pc_type = "asm"
# residual monitoring
ksp_monitor = ""
ksp_converged_reason = ""
};
\end{lstlisting}
and for Poisson-Nernst-Planck equations, we use the ~\textsc{Petsc}~\textsc{SNES} for the Newton iteration. We used the BiCGStab linear solver (a Krylov space solver) with additive Schwarz-based preconditioning for the inner linear solves of the Newton iteration.  
\begin{lstlisting}
solver_options_pnp = {
ksp_atol = 1e-8
ksp_rtol = 1e-8
ksp_stol = 1e-8
#multigrid
ksp_type = "fgmres"
pc_type = "gamg"
pc_gamg_asm_use_agg = True
mg_levels_ksp_type = "gmres"
mg_levels_pc_type = "sor"
mg_levels_ksp_max_it = 10
};
\end{lstlisting}
We used default tolerance for ~\textsc{SNES} convergence.
\bibliographystyle{elsarticle-num-names}

\bibliography{total_reference,hari,references_new}



\end{document}